\documentclass{article}
\usepackage{amssymb,latexsym}

\def\eqnreset{\setcounter{equation}{0}}
\def\eqsection#1{\section{#1}\eqnreset}

\def\refplg{2}
\def\refplo{2.2}

\def\refbplf{3}
\def\refbplg{3.5}

\def\refplt{4}
\def\refbplc{4.3}

\def\refpli{5}

\def\refplw{6}
\def\refplx{6.1}
\def\refply{6.3}
\def\refplz{6.4}
\def\refpls{6.6}
\def\refbpla{6.8}
\def\refbplb{6.10}

\def\refpld{8.7}

\def\refbplm{10.1}

\def\refplk{12}
\def\refpll{12.1}
\def\refpln{12.6}

\def\refplp{13.1}
\def\refplq{13.4}
\def\refplc{13.6}
\def\refplb{13.7}
\def\refplu{13.14}
\def\refbpln{13.15}

\def\refplm{16.1}
\def\refplj{16.2}

\def\refbpli{17.6}

\def\refbplk{18.3}
\def\refbplj{18.6}
\def\refplv{18.8}

\def\refple{19}

\def\refplf{20}

\def\refbpll{22.1}
\def\refbplh{22.2}

\def\refplh{23}
\def\refpla{23.1}
\def\refbpld{23.3}
\def\refplr{23.4}

\newtheorem{Thm}{Theorem}[section]
\newtheorem{Defi}[Thm]{Definition}
\newtheorem{Cor}[Thm]{Corollary}
\newtheorem{Lemma}[Thm]{Lemma}
\newtheorem{Prop}[Thm]{Proposition}
\newtheorem{Rem}[Thm]{Remark}
\newtheorem{Conj}[Thm]{Conjecture}
\newtheorem{Prelim}[Thm]{Preliminary}

\newenvironment{thm}[0]{\begin{Thm}\noindent}%
{\end{Thm}}
\newenvironment{defi}[0]{\begin{Defi}\noindent\rm}%
{\end{Defi}}
\newenvironment{cor}[0]{\begin{Cor}\noindent}%
{\end{Cor}}
\newenvironment{lemma}[0]{\begin{Lemma}\noindent}%
{\end{Lemma}}
\newenvironment{prop}[0]{\begin{Prop}\noindent}%
{\end{Prop}}
\newenvironment{rem}[0]{\begin{Rem}\noindent\rm}%
{\end{Rem}}
{\end{Conj}}
{\end{Prelim}}

\def\medno{\medbreak\noindent}

\def\qed{~\hfill$\square$\medbreak}
\def\proof{\par\noindent{\bf Proof:}{\ }{\ }}

\def\naam#1{\label{#1}}
\def\refer#1{\ref{#1}}

\def\bib#1{\cite{#1}}

\def\text#1{\;\;\;\;{\rm \hbox{#1}}\;\;\;\;}
\def\qquad{\quad\quad}

\def\itema{\vspace{-1mm}\item[{\rm (a)}]}
\def\itemb{\item[{\rm (b)}]}
\def\itemc{\item[{\rm (c)}]}

\def\msy#1{{\mathbb #1}}
\def\C{{\msy C}}
\def\N{{\msy N}}
\def\Z{{\msy Z}}
\def\R{{\msy R}}
\def\D{{\msy D}}

\def\ga{\alpha}

\def\gd{\delta}

\def\gf{\varphi}
\def\gg{\gamma}

\def\gl{\lambda}

\def\gs{\sigma}

\def\gD{\Delta}

\def\gS{\Sigma}

\def\got#1{\mathfrak #1}
\def\fa{{\got a}}
\def\fb{{\got b}}
\def\fg{{\got g}}
\def\fh{{\got h}}
\def\fj{{\got j}}
\def\fk{{\got k}}
\def\fl{{\got l}}
\def\fm{{\got m}}
\def\fn{{\got n}}

\def\fp{{\got p}}
\def\fq{{\got q}}

\def\implies{\Rightarrow}
\def\to{\rightarrow}
\def\Re{{\rm Re}\,}

\def\inp#1#2{\langle#1\,,\,#2\rangle}
\def\hinp#1#2{\langle#1\,|\,#2\rangle}
\def\Ad{{\rm Ad}}
\def\End{{\rm End}}
\def\Hom{{\rm Hom}}

\def\after{\,{\scriptstyle\circ}\,}

\def\pr{{\rm pr}}

\def\ik{{\rm k}}

\def\iq{{\rm q}}
\def\ip{{\rm p}}
\def\iC{{\scriptscriptstyle \C}}

\def\cA{{\cal A}}

\def\cC{{\cal C}}

\def\cE{{\cal E}}
\def\cF{{\cal F}}
\def\cH{{\cal H}}

\def\cM{{\cal M}}

\def\cP{{\cal P}}

\def\cS{{\cal S}}

\def\cV{{\cal V}}
\def\cW{{\cal W}}

\mathcode`:="603A
\def\col{\,:\,}

\def\fC{{\mathfrak C}}
\def\fH{{\mathfrak H}}
\def\fc{{\mathfrak c}}
\def\fj{{\mathfrak j}}

\def\fZ{{\mathfrak Z}}

\def\ep{{\rm ep}}
\def\ev{{\rm ev}}
\def\hyp{{\rm hyp}}

\def\rmU{{\rm U}}
\def\glob{{\rm glob}}

\def\indm{{\rm m}}
\def\inda{{\rm a}}

\def\bs{\backslash}

\def\asmid{\mid}

\def\cC{{\cal C}}
\def\Hyp{{\cal H}}
\def\Hypvee{\Hyp^\vee}
\def\dvee{d^\vee}

\def\tF{\widetilde F}
\def\derF{\widetilde F}

\def\nj{j^\circ}

\def\DrP{\gD_r(P)}
\def\DrQ{\Delta_r(Q)}
\def\gDr{\gD_r}

\def\gSr{\Sigma_r}

\def\spX{{\rm X}}
\def\spXp{\spX_+}

\def\spXQv{\spX_{Q,v}}

\def\spXoneQv{\spX_{1Q, v}}

\def\spXQvp{\spX_{Q,v,+}}

\def\cAtwo{\cA_2}
\def\cAtwoQ{\cA_{2,Q}}

\def\LtwosubA{{\mathfrak L}^2_{0\cA}}

\def\Cep{C^{\rm ep}}

\def\Mer{\cM}

\def\Ci{C^\infty}

\def\Cci{C_c^\infty}

\def\DX{{\msy D}(\spX)}

\def\Ad{{\rm Ad}}
\def\Wave{{\cal J}}

\def\Fou{\cF}

\def\FOU{{\mathfrak F}}

\def\Exp{{\rm Exp}\,}

\def\faq{\fa_\iq}

\def\faqd{\fa_{\iq}^*}
\def\faqdc{\fa_{\iq\iC}^*}

\def\faQq{\fa_{Q\iq}}
\def\faQqdc{\fa_{Q\iq\iC}^*}
\def\faQqd{\fa_{Q\iq}^*}
\def\faQqp{\fa_{Q\iq}^+}

\def\faPqdc{\fa_{P\iq\iC}^*}

\def\staQq{{}^*\fa_{Q\iq}}

\def\staQqdc{{}^*\fa_{Q\iq\iC}^*}
\def\staQqd{{}^*\fa_{Q\iq}^*}
\def\faPq{\fa_{P\iq}}

\def\fbdc{\fb_{\iC}^*}

\def\fgc{\fg_\iC}
\def\fgdc{\fg^*_\iC}
\def\fkQp{\fk{\scriptstyle (Q')}}
\def\fkQ{\fk{\scriptstyle (Q)}}
\def\fmQgs{\fm_{Q\gs}}
\def\barfn{\bar \fn}

\def\AQq{A_{Q\iq}}
\def\Aqp{A_\iq^+}
\def\Aq{A_\iq}
\def\AQqp{A_{Q\iq}^+}

\def\stAQq{{}^*A_{Q\iq}}

\def\MQgs{M_{Q\gs}}
\def\MoneQpv{M_{1Q,v}'}
\def\MoneQp{M_{1Q}'}
\def\MoneQ{M_{1Q}}

\def\HoneQ{H_{1Q}}

\def\minparabs{\cP_\gs^{\rm min}}
\def\allparabs{\cP_\gs}

\def\parone{\cP^1_\gs}
\def\repparabs{{{\bf P}_\gs}}

\def\inn{{\rm n}}
\def\ik{{\rm k}}

\def\DQmaps{{\rm D}_Q}

\def\nE{{E^\circ}}
\def\dnE{{}^\backprime E^\circ}

\def\nEtypes{E^\circ_\types}

\def\nEtypesp{E^\circ_\typesp}
\def\nEtypespp{E^\circ_\typespp}

\def\nC{C^\circ}

\def\WKH{W_{K \cap H}}
\def\NKaq{N_K(\faq)}
\def\NKQaq{N_{K_Q}(\faq)}

\def\QcW{{}^Q\cW}

\def\gL{\Lambda}
\def\dotvar{\,\cdot\,}

\def\pr{{\rm pr}}

\def\embeds{\hookrightarrow}
\def\Cartan{\theta}

\def\rmi{{\rm i}}

\def\lspX{l_\spX}

\def\image{{\rm im}\,}
\def\dega{{\rm deg}_a}

\def\reg{{\rm reg}\,}
\def\der#1{\widetilde{#1}}

\def\taup{\tau'}
\def\tautypesp{{\tau_{\typesp}}}
\def\dK{{\widehat K}}
\def\Vtau{V_\tau}
\def\types{\vartheta}
\def\Vtypes{{\bf V}_\types}
\def\tautypes{\tau_\types}
\def\typesp{{\vartheta'}}
\def\typespp{{\vartheta''}}
\def\itypes{{\rm i}_{\typesp, \types}}
\def\Ptypes{P_{\types,\typesp}}
\def\Ptypespp{P_{\typesp,\typespp}}

\def\sphisotypes{\varsigma_\types}
\def\sphisotypesp{\varsigma_\typesp}
\def\sphiso{\varsigma}

\def\typeso{\vartheta_0}
\def\Vtypesp{{\bf V}_{\typesp}}

\def\Vtaup{V_{\tau'}}
\def\dertau{\widetilde\tau}
\def\typespvee{{{\types'}^\vee}}
\def\typesvee{{\types^\vee}}
\def\typesovee{{\types_0^\vee}}

\def\Ind{{\rm Ind}}
\def\Hilb{{\mathfrak H}}

\def\barVxiomega{\bar V(\xi_\omega)}

\def\Ltwosub{{\mathfrak L}^2_0}

\def\Hxii{\cH_\xi^\infty}
\def\fH{{\mathfrak H}}
\def\fL{{\mathfrak L}}

\def\Hxi{{\cal H}_\xi}
\def\Hxiinfty{{\cal H}_\xi^\infty}
\def\discserP{\spX^\wedge_{P,*,ds}}
\def\barVxi{\bar V (\xi)}
\def\barVxiv{\bar V (\xi, v)}
\def\Mps{\widehat M_{\rm ps}}
\def\discserGtau{\discserG(\tau)}
\def\discserG{\spX^\wedge_{ds}}
\def\discserQ{\spX^\wedge_{Q,*,ds}}
\def\Ltwo{{\mathfrak L}^2}

\def\CirG{C_r^\infty(G)}

\def\naam{\label}

\def\refer{\ref}

\def\bib#1{\cite{#1}}

\begin{document}
\title{The Plancherel decomposition\\
for a reductive symmetric space \\{II. Representation theory}}
\author{E.~P.~van den Ban and H.~Schlichtkrull}
\date{}
\maketitle
\begin{abstract}
We obtain the Plancherel decomposition for a reductive symmetric space in the sense
of representation theory. Our starting
point is the Plancherel formula for spherical Schwartz functions, obtained in part I.
The formula for Schwartz functions involves
Eisenstein integrals obtained by a residual calculus.
In the present paper we identify these integrals as matrix coefficients
of the generalized principal series.
\end{abstract}
\tableofcontents
\eqsection{Introduction} \naam{s: Introduction} In this paper we
establish the Plancherel decomposition for a reductive symmetric
space $\spX = G/H,$ in the sense of representation theory. Here
$G$ is a real reductive group of Harish-Chandra's class and $H$ is
an open subgroup of the group $G^\gs$ of fixed points for an
involution $\gs$ of $G.$ This paper is a continuation of the paper
\bib{BSpl1} in the sense that we derive the Plancherel
decomposition from its main result \bib{BSpl1},Thm.~\refpla, the
Plancherel formula for the space $\cC(\spX\col \tau)$ of
$\tau$-spherical Schwartz functions on $\spX.$  Here $(\tau,
\Vtau)$ is a finite dimensional unitary representation of $K,$ a
$\gs$-stable maximal compact subgroup of $G.$ At the end of the
paper, we make a detailed comparison of our results with those of
P.\ Delorme \bib{Dpl}.

The results of this paper were found and announced in the fall of
1995, when both authors were visitors of the Mittag--Leffler
Institute in Djursholm, Sweden. At the same time  Delorme
announced a proof of the Plancherel theorem. For more historical
comments, we refer the reader to the introduction of \bib{BSpl1}.

Before giving a detailed outline of the results of this paper, we shall first
give some background and describe the main result of \bib{BSpl1},
which serves as the basis for this paper. The space $\spX$ carries an invariant
measure $dx;$ accordingly the regular representation $L$ of $G$ in $L^2(\spX)$ is unitary.
The Plancherel decomposition amounts to an explicit decomposition of $L$
as a direct integral of irreducible unitary representations of $G.$
These representations will turn out to be discrete series representations
of $\spX$ and generalized principal series representations of the form
\begin{equation}
\naam{e: princ ser rep}
\pi_{Q, \xi, \nu} =
\Ind_Q^G(\xi\otimes \nu \otimes 1),
\end{equation}
with $Q = M_QA_QN_Q$ a $\gs\Cartan$-stable parabolic
subgroup of $G$ with the indicated Langlands decomposition,
$\xi$ a discrete series representation of the symmetric
space $\spX_Q: = M_Q/M_Q\cap H,$ and $\nu$ a unitary character of $A_Q/A_Q\cap H.$
To keep the exposition simple, we assume here, and in the rest of the introduction,
that the number of open $H$-orbits on $Q\bs G$ is one.
In general, there are finitely many open
orbits, parametrized by a set $\QcW$ of representatives,
and then  $\xi$
should be taken from the discrete series of the spaces
$\spX_{Q,v}:= M_Q / M_Q \cap v H v^{-1},$ for $v \in \QcW.$

Let $\Cartan$ be the Cartan involution associated with $K;$ it
commutes with $\gs.$ Let $\faq$ be a maximal abelian subspace of
the intersection of the $-1$ eigenspaces for $\Cartan$ and $\gs$
in $\fg,$ the Lie algebra of $G.$ We denote by $\allparabs$ the
collection of $\Cartan\gs$-stable parabolic subgroups of $G$
containing $\Aq:= \exp \faq.$ For  $Q \in \allparabs$ we put
$\faQq: = \fa_Q\cap \faq.$ In \bib{BSpl1} we defined  a spherical
Fourier transform $\Fou_Q$ in terms of a so called normalized
Eisenstein integral
$$
 \nE(Q\col \nu) = E^\circ_\tau(Q\col \nu) .
$$
The Eisenstein integral is a $\DX$-finite and
$1 \otimes \tau$-spherical function in
$C^\infty(\spX)\otimes \Hom(\cAtwoQ, \Vtau),$
depending meromorphically on the parameter $\nu \in
\faQqdc.$
Here $\cAtwoQ = \cAtwoQ(\tau)$ is defined as  the space
of Schwartz functions $\spX_Q \to \Vtau$ that are $\tau_Q:=
\tau|_{K\cap M_Q}$-spherical and behave finitely under the algebra
$\D(\spX_Q)$ of invariant differential operators on $\spX_Q.$ The
space $\cAtwoQ$ is finite dimensional, and inherits the Hilbert
structure from the bigger space  $L^2(\spX_Q\col \tau_Q).$ Without
the simplifying assumption, $\cAtwoQ$ is defined as a finite
direct sum of similar function spaces for $\spX_{Q,v},$ as $v \in
\QcW.$

Let $P_0$ be a fixed minimal element of $\allparabs.$
Then  the
Eisenstein integral $\nE(P_0\col \gl)$ is essentially obtained as a
(sum of) matrix
coefficient(s) of a $K$-finite vector with an $H$-fixed distribution
vector of a $\gs$-minimal
principal series representation of the form (\refer{e: princ ser rep})
with $Q = P_0,$
see \bib{Bps2} and
\bib{BSft}.

In contrast, for non-minimal  $Q \in
\allparabs$ the Eisenstein integral $\nE(Q\col \nu)$ is obtained
from $\nE(P_0\col \gl)$ by means of a residual calculus in the
variable $\gl \in \faqdc,$  see \bib{BSpl1}, Eqn.~(\refpld) and
Lemmas 13.16 and 13.10. In particular, for such $Q$ it is a priori
not clear that the normalized  Eisenstein integral
$\nE(Q\col \nu)$ is a matrix coefficient of the  generalized
principal series representation (\refer{e: princ ser rep}).

In terms of the Eisenstein integral, the spherical
Fourier transform is defined
by the formula
$$
\Fou_Q f(\nu) = \int_\spX \nE(Q \col -\bar\nu \col x)^* f(x)\; dx \in \cAtwoQ,
$$
for $f\in \cC(\spX\col \tau)$ and $\nu \in i\faQqd;$ see
\bib{BSpl1}, \S~\refple. The star indicates that the adjoint of an
endomorphism in $\Hom(\cAtwoQ, \Vtau)$ is taken. The transform
$\Fou_Q$ is a continuous linear map from $\cC(\spX\col \tau)$ into
the space $\cS(i\faQqd) \otimes \cAtwoQ$ of Euclidean Schwartz
functions on $i\faQqd$ with values in the finite dimensional
Hilbert space $\cAtwoQ.$ The wave packet transform $\Wave_Q$ is
defined as the adjoint of the Fourier transform with respect to
the natural $L^2$-type inner products on the spaces involved; see
\bib{BSpl1}, \S~\refplf. It is a continuous linear map $\cS(i\faQqd)
\otimes \cAtwoQ \to \cC(\spX\col \tau),$ given by the formula
$$
\Wave_Q \gf (x) = \int_{i\faQqd} \nE(Q\col \nu\col x)\, \gf(\nu) \; d\nu,
$$
for $\gf \in \cS(i\faQqd) \otimes \cAtwoQ$ and $x \in \spX.$
Here $d\nu$ is Lebesgue measure on $i\faQqd,$ suitably normalized.

Two parabolic subgroups $P,Q \in \allparabs$ are called associated
if their $\gs$-split components $\faPq$ and $\faQq$
are conjugate under the Weyl group $W$ of the root system of $\faq$ in $\fg.$
The notion of associatedness defines an equivalence relation $\sim$ on $\allparabs.$
Let $\repparabs$ be a choice of representatives in
$\allparabs$ for the classes in $\allparabs/\sim.$
Then the Plancherel formula for functions in $\cC(\spX\col \tau)$ takes the form
$$
f = \sum_{Q\in \repparabs} [W\col W^*_Q]\, \Wave_Q\Fou_Q f, \qquad (f \in \cC(\spX\col \tau)),
$$
with $W^*_Q$  the normalizer in $W$ of $\faQq.$
The operator $[W\col W^*_Q] \Wave_Q\Fou_Q$ is a continuous projection operator onto a closed subspace
$\cC_Q(\spX\col \tau)$ of $\cC(\spX\col \tau).$ Moreover,
$$
\cC(\spX\col \tau) = \oplus_{Q\in \repparabs} \;\;\;\cC_Q(\spX\col \tau),
$$
with orthogonal summands. It follows from the above that ${[W\col W^*_Q]}^{1/2}
\Fou_Q$ extends  to a partial isometry
from $L^2(\spX\col \tau)$ to $L^2(i \faQqd)\otimes \cAtwoQ.$  Its adjoint extends
${[W\col W^*_Q]}^{1/2}\Wave_Q$ to a partial isometry in the opposite direction.

In the present paper, we build the Plancherel decomposition for
$(L, L^2(\spX))$ from the above results for all $\tau.$
For this
it is crucial
to relate the Eisenstein integral $\nE(Q\col \nu)$
to the
generalized principal series representations
$\pi_{Q,\xi, -\nu}.$

In \bib{D1n}, Delorme has {\it defined}  a normalized Eisenstein integral
$\dnE(Q\col \nu)$ essentially as a matrix coefficient of the
generalized principal series.
One way to establish the
wanted relationship of $\nE(Q\col \nu)$  with the generalized principal
series would thus be to prove the following identity
of meromorphic functions in the variable $\nu \in \faQqdc:$
\begin{equation}
\naam{e: identity of Eis} \nE(Q\col \nu) = \dnE(Q\col -\nu).
\end{equation}
In view of the vanishing theorem of \bib{BSanfam},
the Eisenstein integral $\nE(Q\col \nu)$
can be uniquely characterized in terms of
its annihilating ideal in $\DX$ and its asymptotic behavior
towards infinity on $\spX;$ see \bib{BSpl1}, Def.~\refplb{} and
Prop.~\refplc{}. The identity (\refer{e: identity of Eis}) would
follow if not only the Eisenstein integral on the left-hand side but
also the Eisenstein integral on the right-hand side satisfied
these characterizing conditions.
For the latter to be true one needs
that, for $\psi \in \cAtwoQ,$ the family $\dnE(Q\col \psi): \nu
\mapsto \dnE(Q\col - \nu)\psi$ belongs to the space $\cE_Q^{\rm
hyp}(\spX\col \tau)$ of \bib{BSpl1}, Prop.~\refplc{}. For this in
turn, the full set of  exponents
of $\dnE(Q\col \psi)$
in its asymptotic expansion along $P_0$
must be of a certain form;
see \bib{BSpl1}, Def.~\refply{}.
We have not been able to deduce this type of information from Delorme's
work.  Nevertheless, by following a different strategy we have been able
to establish (\refer{e: identity of Eis}), but only at the end of the paper,
in Corollary 11.21, after a relation of our Eisenstein integrals
with the principal series has been established.

More precisely, the mentioned characterization
of the Eisenstein integral $\nE(Q\col \nu)$
is used to construct certain embeddings of
$(\fg,K)$-modules
\begin{equation}
\naam{e: embedding J}
\pi_{Q,\xi, - \nu} \embeds (L, \Ci(\spX)).
\end{equation}
The  existence of these embeddings, on the level of $(\fg, K)$-modules,
is sufficient to establish the Plancherel decomposition in the sense
of representation theory, Theorem \refer{t: plancherel with dir int}. Further details
will be given at a later stage in this introduction.

At the end of the paper we invoke the
automatic continuity theorem, Theorem \refer{t: automatic
continuity}, due to W.~Casselman and N.R.~Wallach, see \bib{Cas}
and \bib{Wal2}, to show that the embedding (\refer{e: embedding J})
extends to a $G$-homomorphism. This implies that our Eisenstein integrals are
essentially generalized matrix coefficients of $K$-finite and
$H$-fixed distribution vectors of principal series
representations. {}From this information combined with
results of \bib{CDdvn}, the identity
(\refer{e: identity of Eis}) can then be established.

After this motivation, we shall now give an outline of the paper,
in particular describing how the Eisenstein integrals give rise to the
embeddings (\refer{e: embedding J}).

In Section \refer{s: A property of the discrete series}
we show that the discrete part $L^2_d(\spX\col \tau)$ of $L^2(\spX \col \tau)$
is finite dimensional. This fact can be derived from the description
of the discrete series by T.~Oshima and T.~Matsuki in \bib{OMds}.
We show that it can be obtained from
\bib{BSpl1} and weaker information on the discrete
series, also due to \bib{OMds}, namely the rank condition and the fact
that the $\DX$-characters of $L^2_d(\spX)$ are real and regular. The
mentioned result implies that the parameter
space $\cAtwoQ(\tau)$ of the Eisenstein integral equals $L^2_d(\spX_Q\col \tau_Q).$
Accordingly, it may be decomposed in an orthogonal finite dimensional
sum of isotypical subspaces $\cAtwoQ(\tau)_\xi,$ where $\xi \in
\spX_{Q,ds}^\wedge,$ the collection of discrete series for $\spX_Q.$

In Section \refer{s: Eisenstein integrals and induced representations}
we  explain the connection of the Eisenstein integrals with the principal series.
Let $\dK$ be the unitary dual of $K,$ i.e., the collection of equivalence
classes of irreducible unitary representations of $K.$
If $V$ is a locally convex space equipped with a continuous
representation of $K,$ then by $V_K$ we denote the subspace of $K$-finite vectors;
for $\types\subset \dK$ a finite subset we denote by
$V_\types$ the subspace of $V_K$ consisting
of vectors whose  $K$-types belong to $\types.$
Let $\types\subset \dK$ be a finite subset.
We define $\Vtypes$ to be the space of continuous
functions $K \to \C$ that are left $K$-finite with types contained in the set $\types.$
Moreover, we define $\tau_\types$ to be the restriction of the right regular representation
of $K$ to $\Vtypes.$ Let $\gd_e: \Vtypes \to \C$ be evalutation in $e.$
Then $F\mapsto \gd_e\after F$
is a natural isomorphism from $L^2(X\col \tau_\types)$ onto $ L^2(X)_\types.$
Its inverse, called
sphericalization, is denoted by $\sphiso_\types.$

For $\xi \in \spX^\wedge_{Q,ds},$  we denote by  $\barVxi$
the space of continuous linear $M_Q$-equivariant maps $\Hxi \to L^2(\spX_Q).$ This space
is a finite dimensional Hilbert space. We denote by $L^2(K\col \xi)$ the space of
the induced representation $\Ind_{K\cap M_Q}^K(\xi|_{K \cap M_Q}).$ It is well known that the
induced representation
(\refer{e: princ ser rep}) may be realized as a $\nu$-dependent representation in
$L^2(K\col \xi),$
which we shall denote by $\pi_{Q,\xi, \nu}$  as well; this is the so-called compact picture of
(\refer{e: princ ser rep}).

If $\types\subset \dK$ is a finite subset,
there is a natural isometry from $\barVxi \otimes L^2(K\col \xi)_\types$ into
$\cAtwoQ(\tau_\types),$ denoted $T \mapsto \psi_T.$
We show in Section \refer{s: Eisenstein integrals and induced representations}
that we may use the Eisenstein integrals to define a map
$J_{Q, \xi, \nu}: \barVxi \otimes L^2(K\col \xi)_K \to C^\infty(\spX)_K$ by the formula
\begin{equation}
\naam{e: J in intro}
J_{Q,\xi, \nu}(T)(x) = \gd_e [ E^\circ_\types(Q\col \nu \col x) \psi_T] .
\end{equation}
Here $\types\subset \dK$ is any finite subset such that
$T \in \barVxi \otimes \Ci(K\col \xi)_\types$
and $E^\circ_\types$ denotes the Eisenstein integral with $\tau = \tau_\types.$
The map $J_{Q,\xi,\nu}$ is a priori well-defined for $\nu$
in the complement of
the union of a certain set
$\Hyp(Q, \xi)$ of hyperplanes in $\faQqdc.$  This unio
is disjoint from $i\faQqd.$

The main result of the section is Theorem \refer{t: intertwining prop J}. It asserts that
$\Hyp(Q,\xi)$ is locally finite and that, for $\nu$ in the complement of
$\cup \cH(Q, \xi),$
the map $J_{Q, \xi, \nu}$ is
$(\fg, K)$-equivariant for the infinitesimal representations associated
with $1 \otimes \pi_{Q, \xi, -\nu}$ and $L.$ The proof of this result is
given in the next two
sections.  In the first of these we prepare for the proof by showing
that  $\pi_{Q, \xi, \nu}$ is finitely generated, with local uniformity
in the parameter $\nu,$ see Proposition \refer{p: locally uniform generators}. This result
is needed for the proof of the
local finiteness of $\Hyp(Q,\xi).$

In Section \refer{s: Differentiation of spherical functions} the $(\fg, K)$-equivar\-ian\-ce
of the map $J_{Q, \xi, \nu}$ is established. The $K$-equivar\-ian\-ce readily follows
from the definitions. For the $\fg$-equivariance it is necessary to compute
derivatives of the Eisenstein integral of the form
$L_X E^\circ_\tau(Q\col \nu)\psi,$ for $\psi \in \cAtwoQ(\tau)$ and $X \in \fg.$ The computation
is achieved by introducing a meromorphic family of spherical functions
$\derF: \faQqdc \times \spX \to \fgdc \otimes \Vtau$
by the formula
$$
\derF_\nu(x)(Z) = L_Z [E^\circ_\tau(Q\col \nu \col \dotvar)\psi](x),
$$
for $\nu \in \faQqdc,$ $x \in \spX$ and $Z \in \fg_\iC.$  The function
$\der F_\nu$ is $\dertau$-spherical, with
$\dertau := \Ad_K^\vee\otimes \tau$ and
$\Ad_K := \Ad|_K.$ It has the same annihilating ideal in $\DX$ as the Eisenstein
integral $E^\circ_\tau(Q\col \nu)\psi.$
Moreover, its asymptotic behavior on $\spX$ can be expressed in terms of
that of $E^\circ_\tau(Q\col \nu).$ By the mentioned characterization of Eisenstein
integrals this enables us to show that $\derF_\nu$ equals an Eisenstein
integral of the form $E^\circ_{\dertau}(Q\col \nu) \partial_Q(\nu)\psi,$
with $\partial_Q(\nu)$  an explicitly given differential operator
$\cAtwoQ(\tau) \to \cAtwoQ(\dertau),$ see Theorem \refer{t: tilde of nE}.
The $\fg$-equivariance of $J_{Q,\xi, \nu}$ is then obtained by computing the action
of $\partial_Q(\nu)$ on $\psi_T,$ for $T \in \barVxi \otimes \Ci(K\col \xi)_\types;$
see Lemma \refer{l: d Q xi nu and partial Q} and Proposition
\refer{p: left differentiation Eis}. At the end of the section we complete
the proof of Theorem \refer{t: intertwining prop J} by establishing the
local finiteness of $\Hyp(Q,\xi),$ combining the results
of Sections \refer{s: Generators of induced representations} and
\refer{s: Differentiation of spherical functions};
see Proposition \refer{p: uniformity singularities in types}.

In Section \refer{s: The Fourier transform} we define a Fourier transform
$f \mapsto \hat f(Q\col \xi \col \nu)$ from $\Cci(\spX)_K$ to
$\barVxi \otimes L^2(K\col \xi)_K$ by transposition of the map $J_{Q, \xi, -\bar \nu}.$
It is given by the formula
$$
\hinp{\hat f(Q\col \xi \col \nu)}{T} =
\int_X f(x) \overline{J_{Q,\xi, -\bar \nu}(T)(x)}\; dx
$$
and  intertwines the $(\fg, K)$-module of $L$ with that
of $1 \otimes \pi_{Q,\xi, - \nu}.$ In view of
(\refer{e: J in intro}),  the transform $f \mapsto \hat f$
is related to the spherical Fourier transform by the formula
\begin{equation}
\naam{e: relations fous in intro}
\hinp{\hat f(Q\col \xi\col \nu)}{T} = \hinp{\Fou_Q (\sphisotypes f) (\nu)}{\psi_T},
\end{equation}
for $f \in \Cci(X)_\types.$

The established relation (\refer{e: relations fous in intro})
combined with the spherical Plancherel formula implies that the Fourier transform
$f \mapsto \hat f(Q\col \xi\col \nu)$ defines an isometry from $L^2(\spX)$
into the direct integral
\begin{equation}
\naam{e: deco pi in intro}
\pi = \sum_{Q \in \repparabs} \sum_{\xi \in \spX^\wedge_{Q,ds}}
[W:W_Q^*] \int_{i\faQqd} 1 \otimes \pi_{Q, \xi, - \nu}\; d\nu,
\end{equation}
realized in a Hilbert space $\Ltwo.$
The continuous parts of this direct integral are studied in Section \refer{s: A direct integral}.
In Section \refer{s: Decomposition of the regular representation}
it is first shown, in Theorem \refer{t: FOU is isometry},
that the Fourier transform $f \mapsto \hat f$
extends to an isometry $\FOU$ from $L^2(\spX)$ into $\Ltwo.$ Moreover, its restriction
to $\Cci(\spX)_K$ is a $(\fg, K)$-module map into ${\Ltwo}^\infty.$ By an argument
involving continuity and density, it is then shown that $\FOU$ is $G$-equivariant, see
Theorem \refer{t: G equivariance FOU}. At this stage we have established that $\FOU$
maps the regular representation $L$ isometrically into a direct integral decomposition.
For this to give the Plancherel decomposition, we need to show that the image of $\FOU$
is a direct integral with representations that are irreducible and mutually inequivalent
outside a set of Plancherel measure zero. This is done in
Lemma \refer{l: FOU Q isometry onto LtwoQOmega}
and Proposition \refer{p: distinct reps}.
In the process we use results of F.~Bruhat and Harish-Chandra
on irreducibility and equivalence of unitarily parabolically induced representations,
see Theorem \refer{t: HC irreducibility}.
The Plancherel theorem is formulated in Theorem \refer{t: plancherel with dir int}.
Finally, in Theorem \refer{t: image FOU} a precise description of the image of $\FOU$ is given.

At this point it is still not  clear that our description
of the Plancherel decomposition uses the same parametrizations
as the one in  Delorme's paper \bib{Dpl}.
It is the object of the last section
to show that this is indeed the case.
As said, a key idea is to use
the automatic continuity theorem, Theorem \refer{t: automatic continuity},
due to Casselman and Wallach, see \bib{Cas} and \bib{Wal2}.
It implies that the map $J_{Q,\xi, \nu}$
has a continuous linear extension,
hence can be realized by taking the
matrix coefficient with an $H$-fixed distribution vector
of $\Ind_Q^G(\xi\otimes \nu \otimes 1).$ By means of
the description of such vectors
in \bib{CDdvn}, combined with an asymptotic analysis, it is shown
that our Eisenstein integral is related to Delorme's by the identity
(\refer{e: identity of Eis}),
see Corollary \refer{c: comparison nE}.

Finally, the constants $[W:W^*_Q]$ occurring in our formula
(\refer{e: deco pi in intro})
differ from those in the similar formula of Delorme. This is due to different choices of normalizations
of measures, as is explained in the final part of the paper.

\eqsection{Notation and preliminaries} \naam{s: Notation and
preliminaries} Throughout this paper, we use all notation and
preliminaries from \bib{BSpl1}, Sect.~\refplg. In particular, $G$
is a group of Harish-Chandra's class, $\gs$ an involution of $G$
and $H$ an open subgroup of $G^\gs,$ the group of fixed points for
$\gs.$ The associated reductive symmetric space is denoted by
$$
\spX = G/H.
$$
All occurring measures will be normalized according to the
conventions described in \bib{BSpl1}, end of Section \refpli.

Apart from the references just given, we shall give precise
references to \bib{BSpl1}
for additional notation, definitions and results.

\eqsection{A property of the discrete series} \naam{s: A property
of the discrete series} In this section we discuss an important
result on the discrete part of $L^2(\spX),$ which is a consequence
of the classification of the discrete series by T.~Oshima and
T.~Matsuki in \bib{OMds}. In our approach to the Plancherel
formula via the residue calculus, we obtain it as a consequence of
the rank condition and the regularity of the infinitesimal
character, also due to \bib{OMds}, see \bib{BSpl1}, Rem.~\refplj.

In the rest of this section we assume that $(\tau, \Vtau)$ is
a finite dimensional unitary representation
of $K.$ A function $f: \spX \to \Vtau$ is called $\tau$-spherical
if $f(kx) = \tau(k)f(x),$ for all $x \in \spX$ and $k \in K.$ The Hilbert space
of square integrable $\tau$-spherical functions is denoted by $L^2(\spX\col \tau).$
Its discrete part, denoted $L^2_d(\spX\col \tau)$ is defined as in
\bib{BSpl1}, \S~\refplk. The Fr\'echet space of $\tau$-spherical Schwartz functions,
denoted $\cC(\spX\col \tau),$ is defined as in \bib{BSpl1},
Eqn.~(\refpll). The subspace of $\DX$-finite functions in
$\cC(\spX\col \tau)$ is denoted by $\cA_2(\spX\col \tau).$

\begin{prop}
\naam{p: Ltwodtau equal to cAtwotau}
Let $(\tau,\Vtau)$ be a finite dimensional unitary representation
of $K.$ Then
\begin{equation}
\naam{e: Ltwodtau equal to cAtwotau}
L^2_d(\spX\col \tau) = \cA_2(\spX\col \tau).
\end{equation}
Moreover, each of the spaces above is finite dimensional.
\end{prop}

\proof By the reasoning at the end of
the proof of Lemma \refpln{} in \bib{BSpl1} it follows that the
space on the right-hand side of (\refer{e: Ltwodtau equal to
cAtwotau}) is contained in the space on the left-hand side. If the
center of $G$ is not compact modulo $H,$ then it follows
from
\bib{OMds}, see \bib{BSpl1}, Thm.~\refplm, that $\spX$ has no
discrete series. Hence, $L^2_d(\spX) = 0$ and we obtain (\refer{e:
Ltwodtau equal to cAtwotau}).

On the other hand, if $G$ has a compact center modulo $H$  the
result is part of \bib{BSpl1}, Lemma \refpln. \qed

If $(\xi, \Hxi)$ is an irreducible unitary representation of $G,$
let
$
\Hom_G(\Hxi, L^2(\spX))
$
 denote the space of
$G$-equivariant continuous linear maps
from $\Hxi$ into $L^2(\spX).$ This
space is non-trivial if and only if (the class of) $\xi$ belongs to $\discserG,$
the collection of equivalence classes of discrete series
representations of $\spX.$
If $\xi \in \discserG,$ then the mentioned space
is finite dimensional, by the finite multiplicity of the discrete series,
see \bib{Bfm}, Thm.\ 3.1.

For any irreducible unitary representation $\xi,$ the canonical
map from the tensor product $\Hom_G(\Hxi, L^2(\spX)) \otimes \Hxi$
to $L^2(\spX)$ is an embedding, which is $G$-equivariant
for the representations $1 \otimes \xi$ and $L,$ respectively.
We denote its image by $L^2(\spX)_\xi$ and equip the space $\Hom_G(\Hxi, L^2(\spX))$
with the unique inner product that turns the mentioned embedding
into an isometric $G$-equivariant isomorphism
\begin{equation}
\naam{e: canonical embedding discrete series}
m_\xi:\;\; \Hom_G(\Hxi, L^2(\spX)) \otimes \Hxi\;\;
{\buildrel\simeq \over \longrightarrow}\;\; L^2(\spX)_\xi.
\end{equation}
Obviously the space on the right-hand side of
(\refer{e: canonical embedding discrete series}) depends on $\xi$
through its class $[\xi],$ and will therefore also be indicated
with index $[\xi]$ in place of $\xi.$

With the notation just introduced, it follows that
\begin{equation}
\naam{e: Ltwod as dir sum of types}
L^2_d(\spX) = \widehat\oplus_{\omega \in \discserG} L^2(\spX)_\omega,
\end{equation}
with orthogonal summands. Here and elsewhere, the hat over the summation symbol
indicates that the closure of the algebraic direct sum is taken.

If $\omega$ is an equivalence class of an irreducible unitary representation
of $G,$ we write
$L^2(\spX \col \tau)_\omega: = L^2(\spX \col \tau)\cap [L^2(\spX)_\omega \otimes \Vtau].$
It is readily seen that
this space is non-trivial if and only if $\omega$ belongs to $\discserG$ and  has
a $K$-type in common with the contragredient of
$\tau.$ The collection of $\omega$ with this property is denoted by
$\discserGtau.$

\begin{lemma}
\naam{l: finite ds with K type}
The collection $\discserGtau$ is finite. Moreover,
\begin{equation}
\naam{e: Ltwodtau as dir sum}
L^2_d(\spX\col \tau) = \oplus_{\omega \in \discserGtau} \;\;L^2(\spX\col \tau)_\omega,
\end{equation}
where the direct sum is orthogonal and all the summands are finite dimensional.
\end{lemma}

\proof
That the direct sum decomposition is orthogonal and has closure
$L^2_d(\spX\col \tau)$ follows from the similar properties
of (\refer{e: Ltwod as dir sum of types}).
The space on the left-hand side of (\refer{e: Ltwodtau as dir sum})
is finite dimensional, by Proposition \refer{p: Ltwodtau equal to cAtwotau}.
Since all summands on the right-hand side are non-trivial,
the collection parametrizing these summands is finite.
\qed

\begin{rem}
\naam{r: types in cAtwo}
It follows from Proposition \refer{p: Ltwodtau equal to cAtwotau}
that the spaces $L^2(\spX\col \tau)_\omega,$ for $\omega \in \discserG,$
are contained in $\cA_2(\spX\col \tau);$
we therefore also denote them by $\cA_2(\spX\col \tau)_\omega.$
Note that $L^2(\spX\col \tau)_\omega = 0$ for $\omega$ an irreducible unitary
representation of $G$ that does not belong to $\discserG.$
Accordingly, we put $\cA_2(\spX\col \tau)_\omega = 0$ for such $\omega.$
In view of what has been said, the decomposition
(\refer{e: Ltwodtau as dir sum}) may be rewritten as
\begin{equation}
\naam{e: dir sum cArwotau}
\cA_2(\spX\col \tau) = \oplus_{\omega \in \discserG(\tau)} \;\; \cA_2(\spX\col \tau)_\omega.
\end{equation}
\end{rem}

Let $C(K)_K$ denote the space of right $K$-finite continous functions on $K.$
If $\types$ is a finite subset of $\dK,$ the unitary dual of $K,$
then by $C(K)_\types$ we denote the subspace
of $C(K)_K$ consisting of functions with right $K$-types contained in the set
$\types.$ If $\gd \in \dK,$ then $\gd^\vee$ denotes the contragredient representation.
Accordingly, we put $\types^\vee: = \{\gd^\vee \mid \gd \in \types\}.$
We define
\begin{equation}
\naam{e: defi Vtypes}
\Vtypes : = C(K)_{\types^\vee}
\end{equation}
and equip this space with the restriction of
the right regular representation of $K;$ this restriction is denoted by $\tau_\types.$
We endow $\Vtypes$ with the $L^2(K)$-inner product defined by means of normalized Haar measure.
By $\gd_e$ we denote the map $\Vtypes \to \C,$ $\gf \mapsto \gf(e).$

\begin{lemma}
\naam{l: sphericalization}
Let $E$ be a complete locally convex space equipped with a continuous representation of
$K.$ Then the map
$I \otimes \gd_e$ restricts to a topological linear isomorphism
 from $(E \otimes \Vtypes)^K$ onto
$E_\types.$ If $E$ is equipped with a continuous pre-Hilbert structure
for which $K$ acts unitarily, then the isomorphism is an isometry.
In particular, this yields natural isometries
$$
L^2(\spX\col \tautypes) \simeq L^2(\spX)_\types,
\quad
\Cci(\spX\col \tautypes) \simeq \Cci(\spX)_\types,
$$
where the last two  spaces are equipped with the inner products inherited from
the first two spaces.
\end{lemma}

\proof
This is well known and easy to prove.
\qed
The inverse of the isomorphism $I \otimes \gd_e$ will be denoted
by $\sphiso = \sphisotypes;$ see \bib{BSft}, text before Lemma 5, for
similar notation.
Given  a finite subset $\types \subset \dK$ we shall
write $\discserG(\types)$  for $\discserG(\tautypes),$ the set of discrete
series representations that have a $K$-type contained in $\types.$ The
following result  is now an immediate consequence of
Lemma \refer{l: finite ds with K type}.

\begin{cor}
Let $\types \subset \dK$ be a finite set of $K$-types. Then $\discserG(\types)$
is a finite set.
\end{cor}

We end this section with two simple relations between $\sphisotypes$ and $\sphisotypesp,$
for finite subsets $\types, \typesp \subset \dK$ with $\types \subset \typesp.$
Let $E$ be a complete locally convex space equipped with a continuous representation of $K.$
We denote by $\rmi_{\typesp, \types}: E_\types \to E_\typesp$ the natural inclusion map
and by $P_{\types,\typesp}: E_\typesp \to E_\types$ the $K$-equivariant projection map.
Likewise, the inclusion map $\Vtypes \to \Vtypesp$ and the $K$-equivariant
projection $\Vtypesp \to \Vtypes$ (relative to $\tautypesp, \tautypes$)
are denoted by $\rmi_{\typesp, \types}$ and $P_{\types, \typesp},$ respectively.
By $K$-equivariance, the maps $I\otimes  \rmi_{\typesp, \types}$ and
$I \otimes P_{\types, \typesp}$ induce maps
$$
I\otimes  \rmi_{\typesp, \types}:\;
(E \otimes \Vtypes)^K \to  (E \otimes \Vtypesp)^K,\qquad
I\otimes  P_{\types, \typesp}: \;
(E \otimes \Vtypesp)^K \to (E \otimes \Vtypes)^K.
$$

\begin{lemma}
\naam{l: sphiso and rmi}
Let notation be as above. Then
$$
\sphisotypesp \after (I \otimes  \rmi_{\typesp, \types}) =
 \rmi_{\typesp, \types} \after \sphisotypes,
\qquad
\sphisotypes \after (I \otimes  P_{\types, \typesp}) =
 P_{\types, \typesp} \after \sphisotypesp.
$$
\end{lemma}

\proof
The first identity is immediate from the definitions.
The second identity follows from the first by using that the  maps $P_{\types, \typesp}:
E_\typesp \to E_\types$ and $P_{\types, \typesp}:  \Vtypesp \to \Vtypes$ may both
be characterized by the identities $P_{\types, \typesp} \after \rmi_{\typesp, \types} = I$
and $P_{\types, \typesp} \after \rmi_{\typesp, \typesp\setminus \types} = 0.$
\qed

\eqsection{Eisenstein integrals and induced representations}
\naam{s: Eisenstein integrals and induced representations}
Let $Q \in \allparabs.$ We denote by
$\discserQ$ the collection of equivalence
classes of unitary irreducible representations
$\xi \in M_Q$ such that $\xi$ is a discrete series representation
of $\spXQv,$ for some  $v \in \NKaq.$

In this section we describe the relation of the normalized Eisenstein
integral $\nE(Q\col \nu)$ with the induced representations
$\Ind_{Q}^G(\xi \otimes \nu \otimes 1),$ where $\nu \in \faQqdc$ and
$\xi \in \discserQ.$
In the rest of this section we assume $\xi \in \discserQ$ to be
fixed.

Let $\QcW\subset \NKaq$ be a choice of representatives for $W_Q\bs
W / \WKH,$ see \bib{BSpl1}, text after Eqn.~(\refplo). For $v\in
\QcW,$ we equip $\spXQv$ with the left $M_Q$-invariant measure
$dx_{Q,v},$ specified at the end of \bib{BSpl1}, Section \refpli.
Moreover, we define $\bar V(Q,\xi, v) = \bar V(\xi, v)$ by
\begin{equation}
\naam{e: defi barVxiv}
\barVxiv: = \Hom_{M_Q}(\Hxi, L^2(\spXQv)).
\end{equation}
As mentioned in Section \refer{s: A property of the discrete series},
this space is finite dimensional. In accordance with the mentioned
section, we equip it with the unique inner product
that turns the natural map
\begin{equation}
\naam{e: m xi v isometry}
m_{\xi, v} : \;\;\;\barVxiv \otimes \Hxi \;\; \buildrel \simeq \over \longrightarrow \;\;
L^2(\spXQv)_\xi,
\end{equation}
into an isometric $M_Q$-equivariant isomorphism.
We define the formal direct sums
\begin{equation}\naam{e: defi bar V xi}
\barVxi: = \oplus_{v \in \QcW}\;\; \barVxiv,\qquad L^2_{Q,\xi}: =
\oplus_{v \in \QcW}\;\; L^2(\spXQv)_\xi
\end{equation}
and equip them with the direct
sum inner products. The first of these direct sums will also be denoted
by $\bar V(Q, \xi).$ The second of these direct sums is a unitary $M_Q$-module.
The direct sum of the maps $m_{\xi,v}$
as $v$ ranges over $\QcW$ is an isometric isomorphism
\begin{equation}
\naam{e: m xi isometry}
m_\xi: \;\;\barVxi \otimes \Hxi \;\; \buildrel \simeq \over \longrightarrow \;\;
L^2_{Q,\xi}
\end{equation}
that intertwines the natural $M_Q$-representations.

\begin{rem}
\naam{r: remark on discserQ}
If $Q$ is minimal, then $\discserQ$ coincides with the set
$\Mps,$ defined in \bib{Bps1}, p.\ 368. Moreover, $\QcW = \cW$ is a
choice of representatives for $W/\WKH$ in $\NKaq.$ If $v \in \cW,$
and $\eta \in \Hxi^{M \cap v H v^{-1}},$ then the map
$j_\eta: \Hxi \to L^2(M/M\cap vHv^{-1}),$
defined by $j_\eta(v)(m) = \hinp{v}{\xi(m) \eta)},$ is an $M$-equivariant map.
Moreover, $\eta \mapsto j_\eta$ defines an anti-linear map from $V(\xi, v)$
onto $\Hom_M(\Hxi, L^2(M/M\cap vHv^{-1})).$ This gives an identification
of $\overline{V(\xi, v)}$ with $\barVxiv.$ We recall from \bib{Bps1}, p.\ 378,
that we equipped
$V(\xi, v) = \Hxi^{M \cap v H v^{-1}}$ with the restriction of the
inner product from $\Hxi.$ By the Schur orthogonality relations
this implies that the inner product on $\overline{V(\xi, v)}$ coincides
with $\dim(\xi)$ times the inner product on $\barVxiv.$ Let $V(\xi)$ be defined
as in \bib{Bps1}, Eqn.~(5.1). Then $\overline{V(\xi)}\simeq \barVxi$
and the inner product
on $\overline{V(\xi)}$ coincides with  $\dim(\xi)$ times the inner product on $\barVxi.$
\end{rem}

For $\nu \in \faQqdc,$ let $L^2(Q\col \xi \col \nu)$
denote the space of measurable functions
$G \to \Hxi,$ transforming according to the rule
$$
\gf(man x) = a^{\nu + \rho_Q}\, \xi(m)\, \gf(x),\qquad  (x\in G,\,
(m,a,n) \in M_Q \times A_Q \times N_Q),
$$
and satisfying $\int_K \|\gf(k)\|_\xi^2\; dk < \infty.$ As usual we identify measurable functions
that are equal almost everywhere. The space $L^2(Q\col \xi \col \nu)$ is a Hilbert space
for the inner product given by
\begin{equation}
\naam{e: sesquilinear pairing}
\hinp{\gf}{\psi} = \int_{K} \hinp{\gf(k)}{\psi(k)}_\xi\, dk.
\end{equation}
The restriction of the right regular representation of $G$ to this
space is denoted by $\Ind_Q^G(\xi \otimes \nu \otimes 1),$ or more briefly by
$\pi_{Q, \xi, \nu} = \pi_{\xi, \nu}.$

Let $C^\infty(Q\col \xi \col \nu)$ denote the subspace of $L^2(Q\col \xi\col \nu)$ consisting
of functions that are smooth $G \to \Hxiinfty.$ This subspace is $G$-invariant;
the associated $G$-representation in it is continuous for the usual Fr\'echet topology.

\begin{rem}
\naam{r: smooth vectors}
It follows from \bib{BW}, \S~III.7, that the Fr\'echet $G$-module
$C^\infty(Q\col \xi \col \nu)$
equals the $G$-module of smooth vectors for the representation
$\pi_{Q,\xi, \nu},$ equipped with its natural Fr\'echet topology.
\end{rem}

It will be convenient to work with the compact picture of the induced
representation $\pi_{\xi, \nu}.$
Let $L^2(K\col \xi)$ denote the space of square integrable functions
$\gf: K \to \Hxi$ that transform according to the rule
\begin{equation}
\naam{e: rule for Ci K xi}
\gf(mk) = \xi(m) \gf(k),\qquad (k \in K, m\in K_Q).
\end{equation}
Multiplication induces a diffeomorphism $Q \times_{K_Q} K \simeq G.$
Hence, restriction to $K$ induces an isometry from $L^2(Q\col \xi\col \nu)$
onto $L^2(K \col\xi).$ This isometry restricts to a topological linear
isomorphism  from
$\Ci(Q\col \xi \col \nu)$ onto the subspace
$\Ci(K\col \xi)$ of functions in $L^2(K\col \xi)$ that are smooth $K \to \Hxiinfty,$
where the latter space is equipped with the usual Fr\'echet topology.
Via the isometric restriction map we
transfer
$\pi_{\xi, \nu}$ to a $G$-representation in $L^2(K\col \xi),$ also
denoted by $\pi_{Q, \xi, \nu} = \pi_{\xi, \nu}.$

Let $(\tau, \Vtau)$ be a finite dimensional unitary representation of $K.$ We
define
\begin{equation}
\naam{e: Ci K xi tau}
L^2(K\col \xi \col \tau): = [L^2(K\col \xi) \otimes \Vtau]^K.
\end{equation}
By finite dimensionality of $\tau,$ the space in (\refer{e: Ci K xi tau}) is
finite dimensional and contained in
 $C(K, \Hxi) \otimes \Vtau.$

Let $\ev_e$ denote the evalutation map $C(K, \Hxi) \to \Hxi,$ $\gf \mapsto \gf(e),$
and let $\ev_e \otimes I$ denote the induced map
$L^2(K\col \xi\col \tau) \to \Hxi \otimes \Vtau.$

\begin{lemma}
\naam{l: Frobenius reciprocity}
\vbox{\hspace{2cm}}
 \begin{enumerate}
\itema
The map $\ev_e\otimes I$ defines an isometric isomorphism
from $L^2(K\col \xi \col \tau)$ onto the space $(\Hxi \otimes \Vtau)^{K_Q}.$
\itemb
The space $L^2(K\col \xi \col \tau)$ equals its subspace
$\Ci(K\col \xi \col \tau): = [\Ci(K \col \xi) \otimes \Vtau]^K.$
\end{enumerate}
\end{lemma}
\proof
Observe that $L^2(K\col \xi)$ is the representation space for
$\Ind_{K_Q}^K (\xi|_{K_Q}).$ Hence (a) follows by Frobenius reciprocity.
It is readily checked that $\ev_e \otimes I$ maps $\Ci(K\col \xi \col \tau)$ onto
$
(\Hxiinfty \otimes \Vtau)^{K_Q}.
$
The latter space equals
$(\cH_{\xi K_Q} \otimes \Vtau)^{K_Q} = (\cH_{\xi} \otimes \Vtau)^{K_Q};$ hence (b) follows.
\qed
Given $T \in \barVxi \otimes L^2(K\col \xi \col \tau)$ we may now define the element
$\psi_T \in L^2_{Q,\xi} \otimes \Vtau $
by
$$
\psi_T = [m_\xi \otimes I] \after[I\otimes \ev_e \otimes I] (T).
$$
We agree to denote the map $\ev_e \otimes I:\;
L^2(K\col \xi \col \tau) \to (\Hxi\otimes \Vtau)^{K_Q}$ also by $\gf \mapsto \gf(e).$
With this notation,
if $T = \eta \otimes \gf,$ with $\eta \in \barVxi$ and $\gf \in L^2(K \col \xi \col \tau),$
then
\begin{equation}
\naam{e: formula for psi T v new}
\psi_{T, v} = [\eta_v  \otimes I] (\gf(e)),\qquad (v \in \QcW).
\end{equation}
We recall from Remark \refer{r: types in cAtwo}, applied to
the space $\spXQv$ in place of $\spX,$ for $v \in \QcW,$ that
$
[L^2(\spXQv)_\xi \otimes \Vtau]^{K_Q} \simeq \cA_2(\spXQv \col \tau_Q)_\xi,
$
naturally and isometrically. The space
\begin{equation}
\naam{e: defi cAtwoQ xi}
\cAtwoQ(\tau)_\xi:=
\oplus_{v \in \QcW}\;\;\; \cA_2(\spXQv \col \tau_Q)_\xi
\end{equation}
is a subspace of the space $\cAtwoQ(\tau),$ defined in
\bib{BSpl1}, Eqn.~(\refplp), as the similar direct sum without the
indices $\xi$ on the summands. It follows from the above
discussion combined with (\refer{e: defi bar V xi}) that summation
over $\QcW$ naturally induces an isometric isomorphism
\begin{equation}
\naam{e: L two Q tau new}
(L^2_{Q,\xi} \otimes \Vtau)^{K_Q} \simeq
\cAtwoQ(\tau)_\xi,
\end{equation}
via which we shall identify.

\begin{lemma}
\naam{l: T to psi T isometry}
The map $T \mapsto \psi_T$ is an isometry from $\barVxi \otimes L^2(K\col \xi \col \tau)$
onto $\cAtwoQ(\tau)_\xi.$
\end{lemma}

\proof
It follows from Lemma \refer{l: Frobenius reciprocity} that
\begin{equation}
\naam{e: I otimes eve otimes I}
I \otimes \ev_e \otimes I: \;\;\; \barVxi \otimes L^2(K\col \xi \col \tau)\;
\to\; \barVxi \otimes [\Hxi \otimes \Vtau]^{K_Q}
\end{equation}
is an isometric isomorphism.
The map $m_\xi \otimes I$ is an isometry from $\barVxi \otimes \Hxi \otimes \Vtau$ onto
$L^2_{Q,\xi} \otimes \Vtau,$ which intertwines the $K_Q$-actions
$1 \otimes \xi|_{K_Q} \otimes \tau_Q$ and $L|_{K_Q} \otimes \tau_Q.$ Therefore,
it induces an isometry between the subspaces of $K_Q$-invariants, which by
(\refer{e: L two Q tau new}) is identified with an isometry
\begin{equation}
\naam{e: mxi otimes I}
m_\xi \otimes I: \;\; \; \barVxi \otimes [\Hxi \otimes \Vtau]^{K_Q}
\buildrel \simeq \over \longrightarrow \cAtwoQ(\tau)_\xi.
\end{equation}
Since $T \mapsto \psi_T$ is the composition of (\refer{e: I otimes eve otimes I}) with
(\refer{e: mxi otimes I}), the result follows.
\qed
It follows from Lemma \refer{l: sphericalization} that
$$
L^2(K\col \xi \col \tautypes) \simeq L^2(K\col \xi)_\types,
$$
with an isometric isomorphism. The latter space is
equal to $\Ci(K\col \xi)_\types,$ in view of Lemmas \refer{l: Frobenius reciprocity} (b)
and \refer{l: sphericalization}. Accordingly, the map
$T \mapsto \psi_T,$ defined for $\tau = \tautypes,$ may naturally be viewed
as an isometric isomorphism
\begin{equation}
\naam{e: psi T on types}
T \mapsto \psi_T, \qquad \bar V(Q,\xi)  \otimes \Ci(K\col \xi)_\types \;\;
\buildrel \simeq \over \longrightarrow \;\;\cA_{2,Q}(\tautypes)_\xi.
\end{equation}
Moreover, it is given by the following formula, for
$T = \eta \otimes \gf \in \bar V (Q,\xi) \otimes
\Ci(K\col \xi)_\types;$
$$
\pr_v \psi_T  = \eta_v (\gf(e)), \qquad (v \in \QcW).
$$

We now come to the connection with the normalized Eisenstein
integral $E^\circ_\tau(Q\col \nu) = \nE(Q\col \nu),$ defined as in
\bib{BSpl1}, Def.~\refplb{}.
The Eisenstein
integral is meromorphic in the variable $\nu \in \faQqdc,$ as a
function with values in $\Ci(\spX) \otimes \Hom(\cAtwoQ, \Vtau).$
If $\psi \in \cAtwoQ,$ we agree to write $\nE(Q\col \psi \col \nu
\col \dotvar) = \nE(Q\col \nu \col \dotvar)\psi.$ Then $\nE(Q\col
\psi \col \nu) \in \Ci(\spX\col \tau),$ for generic  $\nu \in
\faQqdc.$

We need a `functorial' property of the normalized Eisenstein integral
that we shall now describe.
Let $(\taup, \Vtaup)$ be a second finite dimensional unitary representation
of $K,$ and let $S: \Vtau \to \Vtaup$ be a $K$-equivariant linear map.
Then via action on the last tensor component,
$S$ naturally induces linear maps
$\Ci(K\col \xi \col \tau) \to \Ci(K\col \xi \col \tau'),$
$\cAtwoQ(\tau)_\xi \to \cAtwoQ(\taup)_\xi$ and
$\Ci(\spX\col \tau) \to \Ci(\spX \col \taup)$ that we
all denote by $I \otimes S.$

\begin{lemma}
\naam{l: functoriality of nE}
Let $S: \Vtau \to \Vtaup$ be a $K$-equivariant map as above.
\begin{enumerate}
\itema
Let $T \in  \barVxi\otimes\Ci(K\col \xi \col \tau).$ Then
$
\psi_{[I\otimes I\otimes S]T} = [I\otimes S]\psi_T.
$
\itemb
Let $\psi \in \cAtwoQ(\tau).$ Then
$$
[I\otimes S]E^\circ_\tau(Q\col \psi \col \nu) = E^\circ_{\taup}(Q \col [I \otimes S] \psi \col \nu),
$$
as a meromorphic $\Ci(\spX\col \tau)$-valued identity in the variable $\nu \in \faQqdc.$
\end{enumerate}
\end{lemma}
\proof (a) is a straightforward consequence of the definitions.
Assertion (b) follows from the characterization of the Eisenstein
integral in \bib{BSpl1}, Def.~\refplb{}.
More
precisely, it follows from the mentioned definition and
\bib{BSpl1}, Prop.~\refplc{} (a), that the family $f =
\nE(Q\col \psi)$ belongs to $\cE_Q^\hyp(\spX\col \tau).$ See
\bib{BSpl1}, Def.~\refpls, for the definition of the latter space.
Moreover, still by \bib{BSpl1}, Prop.~\refplc{},
for $\nu$ in a non-empty open subset $\Omega$ of $\faQqdc,$ each
$v \in \QcW$ and all $X \in \faQq$ and $m \in \spXQvp,$
\begin{equation}
\naam{e: q of E types}
q_{\nu - \rho_Q}(Q,v\mid f_\nu, X, m) =  \psi_v(m).
\end{equation}
It readily follows from the definitions that $g: (\nu,x) \mapsto S(f(\nu,x))$
belongs to $\cE_Q^\hyp(\spX\col \taup);$ moreover, (\refer{e: q of E types}) implies that
$$
q_{\nu - \rho_Q}(Q,v\mid g_\nu , X, m) =
S ( \psi_v(m)) = [\pr_v [I\otimes S] \psi] (m),
$$
for all $\nu \in \Omega,$ each $v\in \QcW,$ and all $X \in \faQq$
and $m \in \spXQvp.$ In view of \bib{BSpl1},
Def.~\refplb{} and Prop.~\refplc{} (a), this
implies that $g = \nE(Q\col [I \otimes S] \psi).$ \qed

If $\types \subset \dK$ is a finite subset and $\psi \in \cAtwoQ(\tautypes),$
we denote the associated normalized Eisenstein integral
$E^\circ_{\tautypes}(Q\col \psi \col \nu)$ also by
$\nEtypes(Q\col \psi \col \nu).$ This Eisenstein integral is
a smooth $\tautypes$-spherical function, depending meromorphically
on the parameter $\nu \in \faQqdc.$

Lemma \refer{l: functoriality of nE}
implies an obvious relation between the Eisenstein integrals
$\nEtypes(Q\col \psi \col \nu)$ for different
subsets $\types.$
If $\types \subset \typesp$ are finite subsets of $\dK,$ then
$\Vtypes\subset \Vtypesp.$
The associated
inclusion map is denoted by $\itypes;$ it intertwines $\tautypes$
with $\tautypesp.$ {}From Lemmas \refer{l: sphiso and rmi} and
\refer{l: functoriality of nE} (a) it follows that
\begin{eqnarray}
\psi_{[I \otimes \itypes]T} &=&
\psi_{[I \otimes I \otimes \itypes][I \otimes \sphiso_\types]T}
\nonumber\\
\naam{e: psi T and itypes}
&=&
[I\otimes \itypes] \psi_T,\qquad (T \in \barVxi\otimes \Ci(K\col \xi)_\types ).
\end{eqnarray}
Moreover, from Lemma \refer{l: functoriality of nE} (b) it follows that
\begin{equation}
\naam{e: E psi and itypes}
\nEtypesp(Q\col [I \otimes \itypes] \psi \col \nu) =
[I \otimes \itypes]  \nEtypes(Q\col  \psi \col \nu),
\qquad (\psi \in \cAtwoQ(\tautypes)).
\end{equation}
We have similar formulas for the $K$-equivariant projection operator
$\Ptypes: \Vtypesp \to \Vtypes.$ {}From Lemmas \refer{l: sphiso and rmi}
and \refer{l: functoriality of nE} it follows that
\begin{equation}
\naam{e: psi T and Ptypes}
\psi_{[I \otimes I \otimes \Ptypes]T} =
[I\otimes \Ptypes] \psi_T,
\qquad (T \in \barVxi\otimes \Ci(K\col \xi)_\typesp),
\end{equation}
and
\begin{equation}
\naam{e: E psi and Ptypes}
\nEtypes(Q\col [I \otimes \Ptypes] \psi \col \nu) =
[I \otimes \Ptypes]  \nEtypesp(Q\col  \psi \col \nu),
\qquad (\psi \in \cAtwoQ(\tautypesp)).
\end{equation}

We recall from
\bib{BSpl1}, \S~\refplt, that a $\gSr(Q)$-hyperplanes in $\faQqdc$ is a hyperplane
of the form $(\ga^\perp)_\iC + \xi,$ with $\ga \in \gSr(Q)$ and
$\xi \in \faQqdc.$ The hyperplane is said to be real if $\xi$ may
be chosen from $\faQqd.$
If $\types \subset \dK$ is a finite subset, then by
\bib{BSpl1}, Prop.~\refplu, there exists a locally finite collection $\Hyp$ of real
$\gSr(Q)$-hyperplanes in $\faQqdc$
such that for each $T \in \barVxi \otimes\Ci(K\col \xi)_\types $
the function
$\nu \mapsto \nEtypes(Q\col \psi_T \col \nu)$ has a singular
locus contained in $\cup \Hyp.$ We denote by $\Hyp(Q,\xi,\types)$ the minimal
collection with this property.
It follows from the definition just given  that
$\types \subset \typesp \implies \Hyp(Q,\xi,\types) \subset
\Hyp(Q,\xi, \typesp).$ Let $\Hyp(Q,\xi)$ denote the union of the
collections $\Hyp(Q,\xi, \types),$ as
$\types$ ranges over the collection of
finite subsets of $\dK.$
Then
\begin{equation}
\naam{e: intersection i faQqd with Hyp}
i\faQqd \cap \;\cup \Hyp(Q,\xi) = \emptyset,
\end{equation}
by the regularity theorem for the normalized Eisenstein integral,
see \bib{BSpl1}, Thm.~\refplv.

For $\nu \in \faQqdc\setminus \cup\Hyp(Q,\xi),$ we define the linear map
$$
J_{Q,\xi,\nu} = J_{\xi, \nu}:\;\;   \barVxi \otimes \Ci(K\col \xi)_K \to \Ci(\spX)_K
$$
by
\begin{equation}
\naam{e: defi of J}
J_{\xi, \nu}(T)(x) =
\nEtypes(Q\col \psi_T\col \nu\col x)(e), \qquad (x \in \spX),
\end{equation}
for $\types \subset \dK$ a finite subset and
$T \in  \barVxi\otimes \Ci(K\col \xi)_\types.$
This definition is unambiguous in view of (\refer{e: psi T and itypes})
and (\refer{e: E psi and itypes}).

\begin{thm}
\naam{t: intertwining prop J}
Let $Q \in \allparabs$ and $\xi \in \discserQ.$
The collection $\Hyp(Q,\xi)$ consists of real $\gSr(Q)$-hyperplanes and is locally finite.
Its union is disjoint from  $i\faQqd.$
Let $\nu \in \faQqdc$ be in the complement of this union.
Then $J_{Q, \xi, \nu}$ is a $(\fg, K)$-intertwining map
from $\barVxi \otimes \Ci(K\col \xi)_K,$ equipped with the induced representation
$1 \otimes \pi_{Q, \xi, - \nu},$
to $\Ci(\spX)_K,$ equipped with the $(\fg, K)$-module structure induced
by the left regular representation of $G$ in $\Ci(\spX).$
\end{thm}

The proof of this theorem will be given in the next two
sections. In Section \refer{s: Generators of induced representations} we investigate
uniformity of generators for $\pi_{Q, \xi, \nu}$ relative to the parameter
$\nu.$ In Section \refer{s: Differentiation of spherical functions}
we shall investigate
the effect of left differentiations on
left spherical functions.

\eqsection{Generators of induced representations}
\naam{s: Generators of induced representations}
In this section we  show that
that the parabolically induced representations,
introduced in Section \refer{s: Eisenstein integrals and induced representations}, are
generated by finitely many $K$-finite vectors, with local uniformity in the continuous
induction parameter.

\begin{prop}
\naam{p: locally uniform generators}
Let $Q\in \allparabs$ and let $\xi$ be a unitary representation
of $M_Q$ of finite length.
Assume that $\Omega \subset \faQqdc$ is a bounded subset.
Then there exists a finite subset $\types \subset \dK$ such
that, for all $\nu \in \Omega,$
\begin{equation}
\naam{e: Ci K xi types generates}
\pi_{Q, \xi, \nu}(U(\fg)) C^\infty(K\col \xi)_\types =
C^\infty(K\col \xi)_K.
\end{equation}
\end{prop}

\begin{rem}
In particular, the result holds for $\gs = \Cartan;$ then $Q$ is an
arbitrary parabolic subgroup of $G$ and $\faQq$ equals its usual Langlands split component
$\fa_Q.$
\end{rem}

\proof
It suffices to prove the result for $\xi$ irreducible. We shall do this
by a method given for $\xi$ tempered in \bib{Wal1}, \S~5.5.5.
Let
$$
\omega := \{ \nu \in \faQqdc \mid \;\;
\inp{\Re \nu - \rho_Q}{\ga} > 0, \;\;\; \forall\; \ga \in \gDr(Q).\}
$$
Then for $\nu \in \omega$
we may define the standard intertwining operator $A(\nu)= A(\bar Q\col Q \col \xi \col \nu)$
from $\Ci(Q\col \xi \col \nu)$ to $\Ci(\bar Q\col \xi \col \nu),$ by
$$
A(\nu) f(x) = \int_{\bar N_Q} f(\bar n\, x)\; d\bar n,\qquad (x \in G),
$$
where $d\bar n$ denotes a choice of Haar measure on $\bar N_Q.$
The integral is absolutely convergent; this follows by an argument that
involves
estimates completely analogous to the ones
given for $Q$ minimal in \bib{Bps2}, proof of Lemma 15.6.
It also  follows from these estimates that,
for $f \in \Ci(K\col \xi),$ the function $A(\nu)f \in \Ci(K\col \xi)$ depends
holomorphically on $\nu \in \omega.$

\begin{lemma}
\naam{l: limit of matrix coefficient}
Let $f, g \in \Ci(K\col \xi),$ $\nu \in \omega$
and $X \in \faQqp.$ Then
\begin{equation}
\naam{e: limit and intertwining operator}
\lim_{t \to \infty} e^{t(-\nu + \rho_Q)(X) }\hinp{\pi_{Q,\xi, \nu}(m\exp tX)f}{g}
=
\hinp{\xi(m)[ A(\nu)f](e)}{g(e)}_\xi.
\end{equation}
\end{lemma}

\proof
See \bib{Wal2}, Lemma 10.5.1.
\qed
\medno
{\bf  Completion of the proof of Prop.~\refer{p: locally uniform generators}:\ }
{}From (\refer{e: limit and intertwining operator}) it can
be deduced, by an argument due to Langlands \bib{Lan}, Lemma 3.13, see also Milicic \bib{Mil}, Proof of Thm.~1, that
if $f \in \Ci(Q\col \xi\col \nu)_K$ and $A(\nu)f \neq 0,$ then $f$ is a cyclic
vector for $\pi_{\xi, \nu}$ in the sense that the $(\fg, K)$-module generated
by $f$ equals  $\Ci(Q\col \xi\col \nu)_K.$ See also \bib{Wal2}, Cor.~10.5.2.
We can now prove
the result in the case that the closure of $\Omega$ is contained in $\omega.$
Indeed, assume this to be the case and let $\nu_0 \in \overline \Omega.$
Since $f \mapsto A(\nu_0)f(e)$ can be expressed as a convolution
operator with non-trivial kernel, there exists a finite set $\types \subset \dK$
and a function
$f \in \Ci(K\col \xi)_\types$ such that $A(\nu_0) f(e) \neq 0;$ by continuity in
the parameter $\nu$ there exists an open neighborhood $\omega_0$ of $\nu_0$ in
$\omega$ such that $A(\nu) f(e) \neq 0$ for all $\nu \in \omega_0.$
{}From what we said above, it follows that (\refer{e: Ci K xi types generates}) holds for all
$\nu \in \omega_0.$ By compactness of the set $\overline \Omega,$ the
result now readily follows in case $\overline \Omega$ is contained in $\omega.$

We shall now use tensoring with a finite dimensional representation
to extend the result to  an arbitrary bounded subset $\Omega \subset \faQqdc.$

Let $P \in \minparabs$ be such that $P\subset Q.$
Let $\gD_Q(P): = \{\ga \in \gD(P)\mid \ga|_{\faQq} = 0 \}$ and put
$\gD(Q) = \gD(P)\setminus \gD_Q(P).$
We fix $n \in \N$ such that $\inp{\Re \nu - \rho_Q}{\ga}/\inp{\ga}{\ga} > - 8n$ for all $\nu \in \bar \Omega$
and $\ga \in \gD(Q).$
We fix $\mu \in \faqd$ with the property that $\inp{\mu}{\ga}/\inp{\ga}{\ga}$ equals
$8 n$ for all $\ga \in \gD(Q)$ and zero for all $\ga \in \gD_Q(P).$
Then $\mu + \overline \Omega \subset \omega.$ Hence there exists
a finite subset $\typesp \subset \dK$ such that
$\pi_{\xi, \nu + \mu}(U(\fg))C^\infty(K\col \xi)_\typesp  = C^\infty(K\col \xi)_K,$
for all $\nu \in \overline \Omega.$

It follows from the condition on $\mu$ that $\inp{\mu}{\ga}/2 \inp{\ga}{\ga} \in 4 \Z$
for all $\ga \in \gD(P).$ Since $\gS$ is a possibly non-reduced root
system, this implies that $\inp{\mu}{\ga}/2 \inp{\ga}{\ga} \in 2 \Z$ for all $\ga \in \gS.$
According to \bib{Bps2}, Cor.~5.7 and Prop.~5.5,
there exists a class one finite dimensional
irreducible  $G$-module $(F,\pi)$ of $\gD(P)$-highest $\faq$-weight $\mu;$
the highest weight space $F_\mu$ is one dimensional, and $M_\gs = M_{P\gs}$ acts trivially on it.
Since $M_{Q\gs}$ centralizes $\faQq,$ it normalizes $F_\mu.$ By compactness it follows
that $(K_Q)_e$ acts trivially on $F_\mu.$
Since $\mu$ vanishes on $\staQq= \faq \cap \fm_Q,$ it follows that $\stAQq$ also acts trivially
on $F_\mu.$ Finally, since $M_{Q\gs}$ is generated by $M_\gs, (K_Q)_e$ and $\stAQq,$ it follows
that $M_{Q\gs}$  acts by the identity on $F_\mu.$

Let $e_\mu \in F_\mu$ be a non-trivial highest weight vector.
Then the map $m: F^* \to \Ci(G)$ defined by $m(v)(x) = v(\pi(x^{-1})e_\mu)$
is readily seen to be an equivariant map from $F^*$ into $\Ci(Q\col 1 \col -\mu ).$
The map $M_\nu: \Ci(Q\col \xi \col \nu + \mu ) \otimes F^* \mapsto \Ci(Q\col \xi \col \nu)$
given by $M_\nu(\gf \otimes v) = m(v) \gf$ is $G$-equivariant, for every $\nu \in \faQqdc.$

Let $v_K \in F^*$ be a non-trivial $K$-fixed vector, then, since
$G = QK,$ the function $m(v_K)$ is nowhere vanishing. {}From this we
see that $M_\nu$ is surjective, for every $\nu \in \faQqdc.$ It
follows  that the $U(\fg)$-module generated by $V_\nu: =
M_\nu(\Ci(Q\col \xi \col \nu + \mu )_\typesp \otimes F^*)$ equals
$\Ci(Q\col \xi \col \nu)_K,$ for all $\nu \in \overline \Omega.$
Let $\types \subset \dK$ be the collection of all $K$-types
occurring in $\gd \otimes F^*$ for some $\gd \in \typesp.$ Then
$\types$ is a finite set and  $V_\nu \subset \Ci(Q\col \xi \col
\nu)_\types,$ for all $\nu \in \faQqdc.$ Hence, (\refer{e: Ci K xi
types generates}) follows for all $\nu \in \Omega.$ \qed
\eqsection{Differentiation of spherical functions} \naam{s:
Differentiation of spherical functions} In this section we assume
$(\tau, \Vtau)$ to be a finite dimensional unitary representation
of $K.$ We shall investigate the action of $L_Z,$ for $Z \in \fg,$
on the Eisenstein integral $\nE(Q\col \nu).$ Here $L$ denotes the
infinitesimal left regular representation. As a preparation, we
shall first investigate the action of $L_Z$ on functions from the
space $\Ci(\spX_+\col \tau),$ defined in \bib{BSpl1}, \S~\refplw.
Secondly, we shall investigate the action of $L_Z$ on families
from $\cE^\hyp_Q(\spX\col \tau),$ defined in \bib{BSpl1},
Def.~\refpls.

Given a function $F \in \Ci(\spXp\col \tau),$
we define the function
$\der{F}: \spXp \to \fgdc \otimes \Vtau\simeq \Hom(\fg_\iC, \Vtau)$
by
$$
\der{F}(x)(Z) =  L_Z F(x),
\qquad (x \in \spXp,\; Z \in \fg_\iC).
$$
One readily checks that
$$
\der{F}(kx)(Z) = \tau(k) \der{F}(x)(\Ad(k^{-1}) Z),
\qquad (x \in \spXp,\, k\in K,\, Z \in \fg_\iC).
$$
Hence, $\der{F}$ is a spherical function of its own right.
In fact, let $\Ad^\vee_K$ denote the restriction to $K$
of the coadjoint representation of $G$
in $\fgdc$ and put $\dertau : = \Ad_K^\vee \otimes \tau.$ Then
$$
\der{F} \in \Ci(\spXp \col \dertau).
$$
Our first objective is to show that if $F$ has a certain
converging expansion towards infinity along $(Q,v),$ for $Q \in
\allparabs$ and $v \in \NKaq,$ then $\derF$ has a similar
expansion, which can be computed in terms of that of $F.$ As a
preparation, we study sets consting of points of the form $m av,$
where $v \in \NKaq,$ $m \in M_{Q\gs}$ and  $a \to \infty$ in
$\AQqp.$ They describe regions of convergence for the expansions
involved, in the spirit of \bib{BSanfam}, \S~3. We will also
describe decompositions of elements of $\fg$ along such sets, in a
fashion similar to \bib{BSanfam}, \S~4. These will be needed to
compute the expansion of $\derF.$

Let $Q \in \allparabs.$ We define the function $R_{Q,v}: \MoneQ \to\, ]\,0,\infty\,[$
as in \bib{BSanfam}, Section 3. Recall that $R_{Q,v}$
is left $K_Q$- and right $\MoneQ \cap v H v^{-1}$-invariant; thus, it may be viewed
as a function on $\spXoneQv.$ If $Q = G,$ then $R_{Q,v}$ equals the constant
function $1$  and if $Q \neq G,$ then according to \bib{BSanfam}, Lemma 3.2,
it is given by
$$
R_{Q,v}(au) = \max_{\ga \in \gS(Q)} a^{-\ga},
$$
for $a \in \Aq$ and $u \in \NKQaq.$
The inclusion map $M_Q \to \MoneQ$ induces an embedding via which we
may identify $\spXQv$ with a sub $M_Q$-manifold of $\spXoneQv.$
{}From
\bib{BSanfam}, Lemma 3.2, we recall that $R_{Q,v} \geq 1$ on $\spXQv.$
\begin{lemma}
\naam{l: R Q v and conjugation}
Let $v \in \NKaq$ and put $Q' = v^{-1} Q v.$ Then
$$
R_{Q,v}(m) = R_{ Q' , 1}(v^{-1} m v), \qquad (m \in \MoneQ).
$$
\end{lemma}

\proof
This follows immediately from the characterization of $R_{Q,v}$ given above.
\qed

In accordance with  \bib{BSanfam}, Eqn.~(3.7),
we define, for $v \in \NKaq$ and $R > 0,$
$$
M_{1Q,v} [R] := \{ m \in \MoneQ \mid R_{Q,v}(m)  < R \},
$$
and
$M_{Q\gs, v}[R]: = M_{Q\gs} \cap M_{1Q, v}[R].$
Note that $M_{1Q,1}[R]$ and $M_{Q\gs, 1}[R]$ equal the sets $M_{1Q}[R]$ and
$M_{Q\gs}[R],$ defined in \bib{BSanfam},
text preceding Lemma 4.7, respectively.
Finally, for  $R > 0$ we define
\begin{equation}
\naam{e: defi AQqp R}
\AQqp(R):= \{ a \in \AQq \mid a^{-\ga} < R \text{for all} \ga \in \gDr(Q)\}.
\end{equation}

\begin{lemma}
\naam{l: Q v sets and conjugation}
Let $v \in \NKaq$ and put $Q' = v^{-1} Q v.$ Let $R > 0.$

\begin{enumerate}
\itema
$
M_{1Q,v}[R] = v M_{1Q'}[R] v^{-1},\quad
M_{Q\gs,v}[R] = v M_{Q'\gs}[R] v^{-1}.
$
\itemb
$\AQqp(R) = v A_{Q'\iq}^+(R) v^{-1}.$
\end{enumerate}
\end{lemma}

\proof
Assertion (a) follows readily from combining Lemma
\refer{l: R Q v and conjugation} with the definitions of the sets involved.
Assertion (b) is clear from (\refer{e: defi AQqp R}).
\qed
\medbreak

We define the open dense subset $\MoneQp$ of $\MoneQ$
as in \bib{BSanfam}, Eqn.\ (4.3). Write $\fg^\pm:= \ker (-I \pm\Cartan\gs)$
and put $\fg_\ga^\pm:=
\fg_\ga \cap \fg^\pm,$ for $\ga \in \gS.$ Write $\HoneQ:= \MoneQ \cap H.$
Then by \bib{BSanfam}, Cor.~4.2,
\begin{eqnarray}
M_{1Q}' & = & K_Q \;[M_{1Q}'\cap \Aq]\; \HoneQ,\nonumber\\
\naam{e: descr MoneQp new}
M_{1Q}' \cap \Aq &=&
\{a \in \Aq\mid
 a^{\ga} \neq 1\;
\mbox{\rm for all}\, \ga \in \Sigma(Q)\, \mbox{\rm with}\; \fg_\ga^+ \neq 0 \}.
\end{eqnarray}
In particular, $M_{1Q}'$ is a left $K_Q$- and right $H_{1Q}$-invariant
open dense subset of $\MoneQ.$
If $v \in \NKaq,$ then by $\MoneQpv$ we denote the analogue of the set $\MoneQp$ for
the pair $(G, vH v^{-1}).$

\begin{lemma}
\naam{l: defi MoneQpv}
Let $v\in \NKaq$ and put $Q' = v^{-1} Q v.$ Then
$\MoneQpv: = v M_{1 Q'}' v^{-1}.$
\end{lemma}

\proof
This readily follows from the definition.
\qed

\begin{lemma}
\naam{l: multiplication of convergence sets}
Let $v \in \NKaq.$
\begin{enumerate}
\itema
$M_{1Q,v}[1] \subset M_{1Q,v}'.$
\itemb
Let $R_1, R_2 > 0.$ Then $M_{Q\gs,v}[R_1] \AQqp(R_2) \subset M_{1Q,v}[R_1R_2].$
\end{enumerate}
\end{lemma}

\proof
For $v =1,$ the results are given in \bib{BSanfam}, Lemma 4.7. Let now $v$ be arbitrary
and put $Q' = v^{-1}Q v.$
Using Lemma \refer{l: Q v sets and conjugation}
(a) with $R=1$ and Lemma  \refer{l: defi MoneQpv} we obtain (a)  from the
similar statement with $Q' , 1$ in place of $Q,v.$ Likewise,
assertion (b) follows
by application of Lemma \refer{l: Q v sets and conjugation}.
\qed
\medbreak
We now come to the investigation of decompositions in $\fg,$ needed
for the study of the asymptotic behavior of $\derF.$
Write $\fkQ := \fk \cap (\fn_Q + \barfn_Q).$
Then $I+\Cartan: X \mapsto X + \Cartan X$
is a linear isomorphism from $\barfn_Q$ onto $\fkQ.$
For $\ga \in \gS$ we put $\fk_\ga^\pm: = (I + \Cartan)(\fg_{-\ga}^\pm).$
Then $\fkQ$ is the direct sum of the spaces $\fk_\ga^{\pm},$
for $\ga \in \Sigma(Q).$

\begin{lemma}
\naam{l: fn Q subset of sum}
Let $v \in \NKaq.$ If $m \in \MoneQpv,$ then $\fn_Q \subset \fkQ \oplus \Ad(mv) \fh.$
\end{lemma}
\proof
For $v=1$ this follows from \bib{BSanfam}, Lemma 4.3 (b), with $\bar Q$ in place of $Q.$
If $v$ is arbitrary, put $Q' = v^{-1}Q v.$ Then for $m \in \MoneQpv$ we have
$v^{-1}mv \in M_{1 Q'}',$
hence $\Ad(v^{-1})\fn_Q =  \fn_{Q'} \subset \fkQp \oplus \Ad(v^{-1} m v)\fh,$
and the result follows by application of $\Ad(v).$
\qed
\medbreak
By the above lemma, for $m\in \MoneQpv$ we may define a linear map
$\Phi(m)= \Phi_{Q,v}(m) \in \Hom(\fn_Q, \fkQ)$ by
\begin{equation}
\naam{e: defi Phi}
X \in \Phi(m)X + \Ad(mv)\fh,\qquad (X \in \fn_Q).
\end{equation}
It is readily seen that $\Phi_{Q,v}$ is an analytic $\Hom(\fn_Q, \fkQ)$-valued
function on $\MoneQpv.$

\begin{lemma}
If  $m\in \MoneQpv,$ $k\in K_Q$ and $h \in \MoneQ \cap vH v^{-1},$ then
$$
\Phi(km h) = \Ad(k) \after \Phi(m) \after \Ad(k)^{-1}.
$$
\end{lemma}

\proof
Since $\MoneQ$ normalizes $\fn_Q$ and  $K_Q$ normalizes $\fkQ$ the result is an immediate
consequence of the definition in equation (\refer{e: defi Phi}).
\qed

\begin{lemma}
\naam{l: Phi Q v and conjugation}
Let $v \in \NKaq$ and put $Q' = v^{-1}Q v.$ Then, for all $m \in \MoneQpv,$
$$
\Phi_{Q,v}(m) = \Ad(v) \after \Phi_{Q', 1}(v^{-1} m v) \after \Ad(v)^{-1},
$$
\end{lemma}

\proof
This follows from (\refer{e: defi Phi}), by the same
reasoning as in the proof of Lemma
\refer{l: fn Q subset of sum}.
\qed
\medbreak
Let $\Psi = \Psi_Q: \MoneQp \to \Hom(\barfn_Q, \fkQ)$
be defined as in \bib{BSanfam}, Eqn.\ (4.4).
Then, for $X \in \barfn_Q$
and $m \in \MoneQp,$
\begin{equation}
\naam{e: defi Psi}
X \in  \Ad(m)^{-1} \Psi(m) X + \fh.
\end{equation}

\begin{lemma}
\naam{l: Phi in terms of Psi}
Let $m \in \MoneQp.$ Then
\begin{equation}
\naam{e: Phi in terms of Psi}
\Phi_{Q,1}(m) = - \Psi(m) \after \gs \after \Ad(m^{-1}).
\end{equation}
\end{lemma}

\proof
If $X \in \fn_Q$ and $m \in \MoneQp,$ then
$\gs \Ad(m^{-1})X \in \barfn_Q,$ so that
$
\gs \Ad(m^{-1})X
$
belongs to $\Ad(m^{-1}) \Psi(m) \gs \Ad(m^{-1}) X + \fh.$
Since
$\Ad(m^{-1})X \in  -\gs \Ad(m^{-1})X +\fh,$ this implies that
$$
\Ad(m^{-1})X \in - \Ad(m^{-1}) \Psi(m) \gs \Ad(m^{-1}) X + \fh.
$$
Comparing with the definition of $\Phi_{Q,1}(m)$ given in (\refer{e: defi Phi})
with $v=1,$
we obtain the desired identity.
\qed
\medbreak

In the formulation of the next result we use the terminology of neat
convergence of exponential polynomial series, introduced
in \bib{BSanfam}, \S~1.

\begin{prop}
\naam{p: series for Phi}
Let $v \in \NKaq.$
There exist unique  real analytic $\Hom(\fn_Q, \fkQ)$-valued functions
$\Phi_\mu = \Phi_{Q,v,\mu}$ on $ M_{Q\gs},$
for $\mu \in \N\DrQ,$
such that, for every $m \in M_{Q\gs}$ and all
$a \in \AQqp(R_{Q,v}(m)^{-1}),$
\begin{equation}
\naam{e: series for Phi}
\Phi_{Q,v}(ma)  =  \sum_{\mu \in \N \DrQ} a^{-\mu} \Phi_\mu(m),
\end{equation}
with absolutely convergent series. Moreover, $\Phi_0 = 0.$ Finally, for every $R>1$
the series in (\refer{e: series for Phi}) converges neatly on $\AQqp(R^{-1})$
as  a $\DrQ$-power series with coefficients in $\Ci(M_{Q\gs, v}[R])\otimes \Hom(\fn_Q,\fkQ). $
\end{prop}

\proof
We first assume that $v= 1.$
Let $\Psi_\mu: \MQgs \to \End(\barfn_Q)$ be as in \bib{BSanfam}, Prop.\ 4.8.
Then it follows from combining the mentioned
proposition with (\refer{e: Phi in terms of Psi})
 that, for $m \in M_{Q\gs}$ and $a \in \AQqp(R_{Q,1}(m)^{-1}),$
$$
\Phi(ma) = - (I + \Cartan) \after
\sum_{\mu \in \N \DrQ} a^{-\mu}  \Psi_\mu(m) \after \gs \after \Ad(ma)^{-1},
$$
with absolutely convergent series.
We now see that the restriction of $\Phi(ma)$ to $\fg_\ga,$ for $\ga \in \gS(Q),$
equals
$$
- (I + \Cartan) \after
\sum_{\mu \in \N \DrQ} a^{-\mu - \ga}
\Psi_\mu(m) \after \gs \after \Ad(m)|_{\fg_{\ga}}.
$$
Put $\Phi_0 = 0$ and, for $\nu \in \N\DrQ\setminus \{0\},$
define $\Phi_\nu(m)\in \Hom(\fn_Q,\fkQ)$ by
$$\Phi_\nu(m)|_{\fg_\ga}: =
-(I + \Cartan)\after  \Psi_{\nu - \ga}(m) \after \gs \after \Ad(m)|_{\fg_{\ga}}
$$
if $\nu -\ga \in \N\DrQ,$ and by $\Phi_\nu(m)|_{\fg_\ga} = 0$ otherwise.
Then (\refer{e: series for Phi}) follows with absolute convergence.
All remaining assertions
about convergence follow from the analogous assertions in \bib{BSanfam},
Prop.\ 4.8.

We now turn to the case that $v$ is general. Let $Q' = v^{-1} Q v,$
and define
$$
\Phi_{Q,v,\mu}(m) = \Ad(v) \after \Phi_{Q', 1, \Ad(v)^{-1}\mu}(v^{-1} m v)\after \Ad(v)^{-1},
$$
for  $\mu \in \N \gDr(Q)$ and $m \in M_{Q\gs}.$ Then all assertions follow
from the similar assertions with $Q', 1$ in place of $Q,v,$
by application of Lemmas \refer{l: Phi Q v and conjugation} and
\refer{l: Q v sets and conjugation}.
\qed

We now come to the behavior of $L_Z,$ for $Z \in \fg_\iC,$
at points of the form $mav,$
with $v \in \NKaq,$ $m \in M_{Q\gs}$ and $a \to \infty$ in $\AQqp.$
We start by observing that
\begin{equation}
\naam{e: fg as direct sum Lie Q and fk}
\fg = \fn_Q \oplus \faQq \oplus (\fmQgs \cap \fp) \oplus \fk,
\end{equation}
as a direct sum of linear spaces. Accordingly, we write, for $Z \in \fg_\iC,$
\begin{equation}
\naam{e: deco Z in 4 terms}
Z = Z_\inn + Z_\inda + Z_\indm + Z_\ik,
\end{equation}
with terms in the
complexifications of the summands in (\refer{e: fg as direct sum Lie Q and fk}),
respectively. If $\fl$ is a real Lie algebra, then by $U(\fl)$ we denote
the universal enveloping algebra of its complexification, and by
$U_k(\fl),$ for $k \in \N,$ the subspace of elements of order at most $k.$
 For $Z \in \fg_\iC$
we define the element
$D_0(Z) = D_{Q,v,0}(Z)$ of $U_1(\fmQgs) \otimes U_1(\faQq) \otimes \End(\Vtau)$ by
\begin{equation}
\naam{e: D zero of Z}
D_0(Z): =  Z_\indm \otimes I \otimes I + I \otimes Z_\inda \otimes
I + I\otimes I \otimes
 \tau(\check Z_\ik),
\end{equation}
where $X \mapsto \check X$ denotes the canonical anti-automorphism of $U(\fg).$
If, moreover, $m \in M_{Q\gs},$
we define, for $\mu \in \N \DrQ\setminus \{0\},$ the element
$D_\mu(Z, m)= D_{Q,v,\mu}(Z,m)$ of $U_1(\fmQgs) \otimes U_1(\faQq) \otimes \End(\Vtau)$ by
$$
D_\mu(Z, m):= I \otimes I \otimes \tau(\Phi_{Q,v,\mu}(m) \check Z_\inn).
$$
Finally, if $m \in \MQgs$ and $a \in \AQqp(R_{Q,v}(m)^{-1}),$
we define the element
$D_{Q,v}(Z,a, m) \in U_1(\fmQgs) \otimes U_1(\faQq) \otimes \End(\Vtau)$
by
\begin{equation}
\naam{e: expansion for D}
D_{Q,v}(Z, a, m) =  \sum_{\mu \in \N\DrQ} a^{-\mu} D_\mu(Z, m),
\end{equation}
where we have put $D_0(Z,m) = D_0(Z).$
We also agree to write
$$
D_{Q,v}^+(Z, a, m) := D_{Q,v}(Z, a, m) - D_0(Z).
$$
It follows from Prop.\ \refer{p: series for Phi} that, for each
$R > 1,$  the series (\refer{e: expansion for D})
is neatly convergent on $\AQqp(R^{-1})$ as a $\DrQ$-exponential
series with values in
$C^\infty(M_{Q\gs,v}[R]) \otimes U_1(\fmQgs) \otimes
U_1(\faQq) \otimes \End(\Vtau).$  Moreover,
\begin{equation}
\naam{e: D plus and Phi}
D_{Q,v}^+(Z,a,m) =
I \otimes I \otimes \tau(\Phi_{Q,v}(ma)\check Z_{\rm n}).
\end{equation}
In the formulation of the following
result we use the notation of the paper \bib{BSanfam}, Sections 1-3.
Via the left regular representation, we view
$U(\fm_{Q\gs}) \otimes U(\faQq) \otimes \End(\Vtau)$
as the algebra of right-invariant differential operators
on $M_{1Q} \simeq M_{Q\gs}\times \AQq,$
with coefficients in $\End(\Vtau).$

\begin{prop}
\naam{p: tilde F in terms of F}
Let $F \in \Cep(\spXp\col \tau).$ Then
$\der{F} \in \Cep(\spXp \col \dertau).$ Moreover, if $Q \in \allparabs$ and
$v \in \NKaq,$ then $\Exp(Q,v\asmid \der{F}) \subset \Exp(Q,v\asmid F) - \N \DrQ.$
Finally, for every $Z \in \fg_\iC,$ the $\DrQ$-exponential expansion
\begin{equation}
\naam{e: expansion tilde F along Q v}
\der{F}(mav)(Z) = \sum_{\xi} a^{\xi}\,q_\xi(Q,v \asmid \der{F}, \log a, m)(Z)
\end{equation}
along $(Q,v)$ arises from the similar expansion
\begin{equation}
\naam{e: expansion F along Q v}
F(mav)  = \sum_{\xi} a^{\xi}\,q_\xi(Q,v \asmid F, \log a, m)
\end{equation}
by the formal application of the expansion (\refer{e: expansion for D}).
In particular, if $\xi$ is a leading exponent of $F$ along $(Q,v),$
then, for every $Z \in \fg_\iC,$
\begin{equation}
\naam{e: equation for leading coeff F}
q_\xi(Q,v\asmid \der{F}, \log( \dotvar ) , \dotvar)(Z)
= [ D_{Q,v,0}(Z) - \xi( Z_\inda) ]\, q_\xi(Q,v\asmid F, \log(\dotvar), \dotvar).
\end{equation}
\end{prop}
\proof
It is obvious that $\tF \in C^\infty(\spXp\col \tau).$ We shall investigate
its expansion along $(Q,v),$ for $Q \in \allparabs$ and $v \in \NKaq.$
We start by observing that, for  $R > 1,$  the expansion (\refer{e: expansion F along Q v})
converges neatly on $\Aqp(R^{-1})$ as a $\DrQ$-exponential polynomial expansion
in the variable $a,$ with coefficients in the space $\Ci(\spXQvp[R]\col \tau_Q),$
see \bib{BSanfam}, Thm.~3.4.

If $\gf$ is a smooth function on a Lie group $L,$ with values in a complete
locally convex space, then for $X \in \fl$ and $x \in L$ we put
$\gf(X;x):= d/dt\,\gf(\exp tX x)|_{t = 0}.$
Accordingly, it follows from (\refer{e: deco Z in 4 terms})
that for $Z \in \fg_\iC,$ and $m \in \MQgs$ and $a \in \AQq$ with $mav \in
\spXp,$ we have  \begin{eqnarray}
\tilde F (mav )(Z) & = & F(\check Z ; ma v)\nonumber\\
&=&
\tau(\check Z_\ik)F(ma v) + F(\check Z_\indm; mav)\nonumber \\
\naam{e: four terms of F}
&& \;\;\;\;\;\;\;\;\;\;\;\;+ F( \check Z_\inda ; mav ) +
F(\check Z_{\rm n} ; ma v).
\end{eqnarray}
The sum of the first three terms allows an expansion that
is obtained by the termwise formal application of $D_{Q,v,0}(Z)$
to the expansion (\refer{e: expansion F along Q v}), by \bib{BSanfam},
Lemmas 1.9 and 1.10. Moreover,
the resulting expansion converges on $\Aqp(R^{-1})$ as a $\DrQ$-exponential polynomial expansion
in the variable $a,$ with coefficients in the space $\Ci(\spXQvp[R], \Vtau).$
Thus, it remains
to discuss the last term in (\refer{e: four terms of F}).
Since $F$ is right $H$-invariant and left $\tau$-spherical,
we see by application of (\refer{e: defi Phi}) and (\refer{e: D plus and Phi})
that the mentioned
term may be rewritten as
\begin{eqnarray*}
F(\check Z_{\rm n} ; mav) &=&
F(\Phi_{Q,v}(ma) \check Z_{\rm n}; mav)\\
&=&
\tau( \Phi_{Q,v}(ma)\check Z_{\rm n}) F( mav )\\
&=&
D_{Q,v}^+(Z,a, m) F(\dotvar v)(ma).
\end{eqnarray*}
It follows from Proposition \refer{p: series for Phi} that the series for $D_{Q,v}^+(Z)$
converges neatly on $\Aqp(R^{-1})$ as a $\DrQ$-exponential polynomial expansion
in the variable $a,$ with coefficients in the space
$\Ci(M_{Q, \gs, v}[R]) \otimes \End(\Vtau).$
{}From \bib{BSanfam}, Lemma 1.10, it now follows that
$F(\check Z_{\rm n} ; mav)$ admits a $\DrQ$-exponential polynomial
expansion that is obtained by the obvious formal application of the
series for $D_{Q,v}^+(Z, a, m)$ to the series for $F(mav).$  The resulting series
converges neatly on $\AQqp(R^{-1})$  as a $\DrQ$-exponential polynomial
expansion in the variable $a$ with coefficients in $\Ci(M_{Q\gs, v}[R], \Vtau).$
It follows that $\tilde F(mav)(Z)$ has an expansion of the
type asserted along $(Q,v),$ with exponents as indicated.

In particular, if $Q$ is minimal,  it follows that $\tilde F(ma v)$
allows a neatly converging $\gD(Q)$-exponential polynomial expansion in the variable
$a \in \Aqp(Q),$ with coefficients in $\Ci(\spXQv) \otimes \fg_\iC^* \otimes
\Vtau.$  This implies that
$\tF$ belongs to the space $C^\ep(\spXp\col \tilde\tau),$ defined in
\bib{BSanfam}, Def.~2.1.

It remains to prove the assertion about the leading exponent $\xi$
for $F$ along $(Q,v).$ {}From the above
discussion we readily see that the term in the expansion
(\refer{e: expansion tilde F along Q v}) with exponent
$\xi$ is obtained from the application of the constant term $D_{Q,v, 0}(Z)$ of
$D_{Q,v}(Z, a, m)$ to the term in the expansion
(\refer{e: expansion F along Q v}) with exponent $\xi.$
This yields
$$
a^\xi q_\xi(Q,v\asmid \der{F}, \log a , m)(Z)
= D_{Q,v,0}(Z) [ ( m, a) \mapsto  a^\xi q_\xi(Q,v\asmid F, \log a, m)].
$$
Now use that $a^{-\xi}\after  D_{Q,v,0}(Z) \after a^\xi =
D_{Q,v,0}(Z) + \xi(\check Z_\inda)$
to obtain (\refer{e: equation for leading coeff F}).
\qed

We can now describe the action of $L_Z,$ for $Z \in \fg_\iC,$ on
families from the space $\cE^\hyp_Q(\spX\col \tau),$ defined in
\bib{BSpl1}, Def.~\refpls.

\begin{thm}
\naam{t: tilde of family F}
Let $F \in \cE_Q^\hyp(\spX \col \tau).$ Then the family
$\der F: \faQqdc \times \spX \to \fgdc \otimes \Vtau,$ defined by $(\der F)_\nu =
(F_\nu)^\sim$ belongs to $\cE_Q^\hyp(\spX \col \dertau).$
Moreover, for every $Z \in \fgc$ and all $\nu$ in an open dense subset of $\faQqdc,$
\begin{eqnarray}
\lefteqn{\!\!\!\!\!\!\!\!\!\!\!\!\!\!\!\!\!\!
q_{\nu - \rho_Q}(Q,v\asmid \der{F}_\nu \col \log( \dotvar) \col \dotvar)(Z)}
\qquad
\nonumber\\
\naam{e: equation for q of family tF}
\qquad &=&
[D_{Q,v,0}(Z) - (\nu -\rho_Q)( Z_\inda)]\,
 q_{\nu - \rho_Q}(Q,v\asmid F_\nu \col \log(\dotvar) \col \dotvar).
\end{eqnarray}
\end{thm}

\proof There exist $\gd \in \DQmaps$ and a finite subset $Y
\subset \staQqdc$ such that $F \in \cE_{Q,Y}^\hyp(\spX\col
\tau\col \gd).$ Let $\Hyp = \Hyp_F,$ $d = d_F$ and $k = \dega F$
be defined as in the text following \bib{BSpl1}, Def.~\refplx.
Then $F$ satisfies all conditions of the mentioned definition. It
follows from the characterization of the expansions for $\der F$
in Proposition \refer{p: tilde F in terms of F} that $\der F$
satisfies the hypotheses of \bib{BSpl1}, Def.~\refplx. with
$\dertau$ in place of $\tau,$ with the same $Y, \Hyp, d, k.$ In
particular, $\der F$ belongs to $C^{\ep, \hyp}_{Q,Y}(\spXp\col
\dertau).$

Since $F_\nu$ is annihilated by the ideal $I_{\gd, \nu}$ for
generic $\nu \in \faQqdc,$ the same holds for $\tF_\nu,$ and we
see that $\tF \in \cE^\hyp_{Q,Y}(\spXp\col \dertau\col \gd),$ see
\bib{BSpl1}, Def.~\refply.

Let now $s \in W,$ $P \in \parone$ such that $s(\faQq) \not\subset
\faPq$ and $v \in \NKaq.$ Then there exists an open dense subset
$\Omega \subset \faQqdc$ such that $F$ satisfies the condition
stated in \bib{BSpl1}, Def.~\refplz. It follows from Propositon
\refer{p: tilde F in terms of F} and the fact that the functions
$m \mapsto D_{P,v,\mu}(Z,m)$ are smooth on all of $M_{P\gs},$ for
$Z \in \fgc, \mu \in \N\DrP,$ that $\tilde F$ also satisfies the
condition of \bib{BSpl1}, Def.~\refplz, with the same set
$\Omega.$ We conclude that $\tF \in \cE_{Q,Y}^\hyp(\spXp\col
\dertau\col \gd)_\glob.$ In view of \bib{BSpl1}, Lemma \refbpla,
$\nu \mapsto F_\nu$ is a meromorphic $\Ci(\spX\col \tau)$-valued
function on $\faQqdc.$ Hence, $\nu \mapsto \tF_\nu$ is a
meromorphic $\Ci(\spX\col \dertau)$-valued function on $\faQqdc.$
In view of
\bib{BSpl1}, Def.~\refpls, we now infer that $\tF \in \cE_Q^\hyp(\spX
\col \dertau).$

Finally, for  $\nu$ in an open dense subset of $\faQqdc,$ the element $\nu - \rho_Q$
is a leading exponent for $F$ along $(Q,v).$ Thus,
(\refer{e: equation for q of family tF})
follows from (\refer{e: equation for leading coeff F}).
\qed

Next, we apply the above result to the normalized Eisenstein
integral $\nE(Q\col \psi \col \nu),$ defined for $\psi \in \cAtwoQ.$
Let $v\in \QcW.$ Given a function $\psi_v\in \cAtwo(\spXQv\col \tau_Q)$
and an element
$\nu \in \faQqdc$
we define the function
\begin{equation}
\naam{e: partial psi v}
\partial_{Q,v}(\nu)\psi_v: \;\;\spXQv \to \fgdc \otimes \Vtau
\end{equation}
by
$$
\partial_{Q,v}(\nu)\psi_v(x)(Z) =
[L_{Z_\indm} - (\nu - \rho_Q)(Z_\inda)]\psi_v (x)
- \tau(Z_\ik)[\psi_v(x)],
$$
for $x \in \spXQv$ and $Z\in \fgc.$
Clearly, the function $\partial_{Q,v}(\nu) \psi_v$ is a $\D(\spXQv)$-finite
Schwartz function with values in $\fgdc \otimes \Vtau.$ Since $K_Q$ normalizes
the decomposition
(\refer{e: fg as direct sum Lie Q and fk})
and centralizes $\faQq,$ one readily checks
that the function is $\dertau_Q$-spherical, with $\dertau_Q:= \dertau|_{K_Q}.$
Hence,
$$
\partial_{Q,v}(\nu) \psi_v \in \cAtwo(\spXQv \col \dertau_Q).
$$
We define  the map $\partial_Q(\nu): \cAtwoQ(\tau) \to \cAtwoQ(\dertau)$ as
the direct sum, for $v \in \QcW,$ of the maps $\partial_{Q,v}(\nu):
\cAtwo(\spXQv\col \tau_Q) \to \cAtwo(\spXQv \col \dertau_Q).$

\begin{thm}
\naam{t: tilde of nE}
Let $\psi \in \cAtwoQ(\tau)$ and let the family
$F: \faQqdc \times \spX \to \Vtau$ be defined by
$$
F(\nu, x) = E^\circ_\tau(Q\col \psi \col \nu\col x).
$$
Then the family $\der F: \faQqdc \times \spX \to \fgdc \otimes \Vtau,$ defined
by $(\der F)_\nu = (F_\nu)^\sim,$ is given by
$$
\der F (\nu, x) = E^{\circ}_{\dertau}(Q \col \partial_Q(\nu) \psi\col \nu  \col
x). $$
\end{thm}

\proof It follows from \bib{BSpl1}, Def.~\refplb{}
and Prop.~\refplc{}, that the family $F$
belongs to $\cE^\hyp_Q(\spX\col \tau)$ and that the family $G: =
\nE(Q\col \partial_Q(\nu)\psi)$ belongs to $\cE^\hyp_Q(\spX\col
\dertau).$ Let $v \in \QcW.$ Then it follows  from the mentioned
proposition, combined with
\bib{BSanfam}, Thm.~7.7, Eqn.~(7.14), that, for
$\nu $ in an open dense subset of $\faQqdc$ and all
$X \in \faQq$ and $m \in \spXQvp,$
\begin{eqnarray}
\naam{e: q of F is psi}
q_{\nu - \rho_Q}(Q,v\mid F_\nu, X, m) &=&  \psi_v(m),\\
\naam{e: q of G is partial psi}
q_{\nu - \rho_Q}(Q,v\mid G_\nu, X, m) & = &  \partial_{Q,v}(\nu)\psi_v(m).
\end{eqnarray}
{}From Theorem \refer{t: tilde of family F} we see
that $\tF \in \cE^\hyp_Q(\spX\col \dertau).$
Moreover, combining (\refer{e: q of F is psi}) and
(\refer{e: equation for q of family tF})
we infer
that, for $Z \in \fgc,$ $\nu$ in an open dense subset of $\faQqdc$
and all $X \in \faQq$
and $m \in \spXQvp,$
$$
q_{\nu - \rho_Q}(Q,v\mid \tF_\nu, X, m)(Z)
=
[D_{Q,v,0}(Z) - (\nu -\rho_Q)(Z_\inda)] [(m,a) \mapsto \psi_v(m) ].
$$
{}From (\refer{e: D zero of Z}) we see that the expression on the right-hand side of this
equation equals $[\partial_{Q,v}(\nu) \psi_v(m) ](Z);$ hence
\begin{equation}
\naam{e: q of tF and partial psi}
q_{\nu - \rho_Q}(Q,v\mid \tF_\nu, X, m) = \partial_{Q,v}(\nu) \psi_v(m).
\end{equation}
Comparing (\refer{e: q of tF and partial psi}) with (\refer{e: q of G is partial psi})
we deduce that the family $\tF - G \in \cE^\hyp_Q(\spX\col\dertau)$
satisfies the hypothesis of the vanishing theorem,
\bib{BSpl1}, Thm.~\refbplb.
Hence, $\tF = G.$
\qed
\medbreak

Given $\nu \in \faQqdc$ and $\gf \in \Ci(K\col \xi \col \tau)$ we define the function
$$
d(Q, \xi, \nu)\,\gf \in C^\infty(K\col \xi) \otimes \fgdc \otimes \Vtau
$$
by
\begin{equation}
\naam{e: defi d Q xi nu}
[d(Q,\xi, \nu)\, \gf] (k, Z) = ([ \pi_{\xi, -\nu}(Z) \otimes I ]\,\gf)(k),
\end{equation}
for $k \in K$ and $Z \in \fg_\iC.$
One readily verifies that $d(Q,\xi, \nu) \gf \in \Ci(K\col \xi \col \dertau).$

\begin{lemma}
\naam{l: d Q xi nu and partial Q}
Let $T \in \barVxi \otimes \Ci(K\col \xi \col \tau).$
Then, for all $\nu \in \faQqdc,$
$$
\psi_{[I \otimes d(Q, \xi,\nu)]T} = \partial_Q(\nu) \psi_T.
$$
\end{lemma}

\proof
By linearity it suffices to prove this for $T = \eta \otimes \gf,$
with $\eta \in \barVxi$ and
$\gf \in \Ci(K\col \xi \col \tau).$ Let $v \in \QcW$ and $Z \in \fgc.$ Then
combining (\refer{e: defi d Q xi nu}) with the decomposition
(\refer{e: deco Z in 4 terms}) we infer that
$$
[d(Q, \xi, \nu) \gf](e)(Z) = [\xi(Z_\indm) \otimes I - I \otimes \tau(Z_\ik)] \gf(e)
-(\nu - \rho_Q)(Z_\inda) \gf(e).
$$
By equivariance, $\eta_v$ maps $\Hxiinfty$ into $L^2(\spX_{Q,v})_\xi^\infty \subset \Ci(\spXQv),$
intertwining the $(\fm_Q, K_Q)$-actions.
Using formula (\refer{e: formula for psi T v new}) we now obtain that
\begin{eqnarray*}
\lefteqn{ \psi_{[I \otimes d(Q,\xi, \nu)]T, v} (\dotvar)(Z) =} \\
&=&
[\eta_v \otimes I ] ( [\xi (Z_\indm)\otimes I]\gf(e)  )
- [(\nu - \rho_Q)(Z_\inda) I\otimes I + I \otimes \tau(Z_\ik)] (\eta_v \otimes I)(\gf(e))
\\
&=&
[L_{Z_\indm} - (\nu - \rho_Q)(Z_\inda)] (\psi_{T,v})(\dotvar)
- \tau(Z_\ik)[ \psi_{T, v}(\dotvar)]
\\
&=&
(\partial(Q,\nu)\psi_T)_v(\dotvar)(Z).
\end{eqnarray*}
\qed

\begin{cor}
\naam{c: relation with derF for Eis}
Let $T \in \barVxi\otimes \Ci(K\col \xi \col \tau)$
and let the family $F: \faQqdc \times \spX \to \Vtau$
be defined by
$$
F_\nu = \nE(Q\col \psi_T \col \nu).
$$
Then the family $\der F: \nu \mapsto (F_\nu)^\sim$ is given by
\begin{equation}
\naam{e: formula der F}
\der F_\nu = \nE(Q \col  \psi_{[I \otimes d(Q,\xi, \nu)]T} \col \nu).
\end{equation}
\end{cor}

\proof
This follows from Theorem \refer{t: tilde of nE} and Lemma \refer{l: d Q xi nu and partial Q}.
\qed

As a consequence of the above, we can now express derivatives of the normalized
Eisenstein integral in a form needed for the proof of Theorem \refer{t: intertwining prop J}.

\begin{prop}
\naam{p: left differentiation Eis}
Let $\types \subset \dK$ be a finite subset, and let
$\typesp \subset \dK$ be the union of the collections
of $K$-types occurring in $\Ad_K \otimes \gd,$ as $\gd \in \types.$
Let $T \in  \barVxi\otimes \Ci(K\col \xi)_\types .$ Then
$(I \otimes \pi_{\xi, -\nu}(Z) )T \in \barVxi \otimes \Ci(K \col
\xi)_\typesp,$  for all $Z \in \fgc$ and $\nu \in \faQqdc.$ Moreover, for all
$Z \in \fgc,$ $x \in \spX$  and $k \in K,$
\begin{equation}
\naam{e: intertwining prop Eis}
L_{\Ad(k)^{-1} Z} E^\circ_{\types} (Q\col \psi_{T} \col \nu) (x)(k) =
E^\circ_\typesp(Q\col \psi_{[I \otimes \pi_{\xi, -\nu}(Z)]T}\col \nu)(x)(k),
\end{equation}
as a meromorphic identity in $\nu \in \faQqdc.$
\end{prop}

\proof
Let $\tau = \tautypes.$ We shall use the natural identification
$\Ci(K\col \xi)_\types \simeq \Ci(K \col \xi\col \tau)$ of
Lemma \refer{l: sphericalization},
so that
$\psi_T$ may be viewed as an element  of
$\barVxi \otimes \Ci(K \col \xi\col \tau).$

Define the family $F$ as
in Corollary \refer{c: relation with derF for Eis}. We shall derive the identity
(\refer{e: intertwining prop Eis}) from (\refer{e: formula der F}) by using the
functorial properties of Lemma \refer{l: functoriality of nE}.

Fix $Z \in \fgc.$
We define the matrix coefficient map
$m_Z: \fgdc \to \Ci(K)$ by
$$
m_Z(\zeta)(k) = \zeta(\Ad(k^{-1})Z), \qquad (\zeta \in \fgdc,\; k \in K).
$$
The map $m_Z$ intertwines the representation $\Ad_K^\vee$ of $K$ in
$\fgdc$ with the right regular
representation of $K$ in $C(K).$ In particular, it
maps into the finite dimensional space $C(K)_{\typesovee},$ with
$\typeso\subset \dK$ the set of $K$-types in $\Ad_K.$
We define the equivariant map
$$
S_1:= m_Z \otimes I:\;\; \fgdc \otimes \Vtypes \to C(K)_{\typesovee} \otimes
\Vtypes. $$
On the other hand, we define the map  $S_2: C(K)_{\typesovee} \otimes \Vtypes \to C(K)$
by $\phi \otimes \psi \mapsto \phi \psi.$
This map intertwines $\tau_{\typeso} \otimes \tautypes$ with the right regular
representation of $K$ in $C(K),$ hence maps into $C(K)_\typespvee.$
The space $C(K)_{\typesovee} \otimes \Vtypes$ may be naturally identified with
a finite dimensional $K$-submodule of $C(K\times K),$ the latter being
equipped with the diagonal $K$-action from the right. Under this identification
the map $S_2$ corresponds with
the restriction of the map $\Delta^*: C(K\times K) \to C(K)$ given by
$\Delta^*\gf(k) = \gf(k,k).$

The map $S = S_2 \after S_1: \;\fgdc \otimes \Vtau \to \Vtypesp$
is $K$-equivariant. We shall apply $I\otimes S$ to both sides of
the identity (\refer{e: formula der F}). Application of $I\otimes S_1$
to the left-hand side yields
$(I \otimes S_1)[\der{F}_\nu(\dotvar)](k) = \der{F}_\nu(\dotvar, \Ad(k^{-1})Z),$ which in
turn equals $L_{\Ad(k^{-1}) Z} F_\nu.$
By application of $I \otimes S_2$ to the latter function we find
\begin{eqnarray}
(I\otimes S )[\der{F}_\nu(\dotvar)](k) &=&
L_{\Ad(k^{-1})Z} F_\nu (\dotvar)(k)\nonumber\\
\naam{e: I otimes S der F as L Z Eis}
&=&
L_{\Ad(k^{-1})Z} \nE(Q\col \psi_{T} \col \nu) (\dotvar) (k).
\end{eqnarray}
On the other hand, from Lemma \refer{l: functoriality of nE} we see
that application of $I\otimes S$ to the expressions on both sides  of
(\refer{e: formula der F}) yields
\begin{equation}
\naam{e: I otimes S on right}
(I\otimes S )\der{F}_\nu=
\nE(Q \col  \psi_{[I \otimes I \otimes S][I \otimes d(Q,\xi, \nu)]T} \col \nu).
\end{equation}
We observe that $(I \otimes S)\after d(Q,\xi, \nu)$ is a linear map from
$\Ci(K\col \xi\col \tautypes)$ to $\Ci(K \col \xi \col  \tautypesp)$
and claim that the following diagram commutes, for every $\nu \in \faQqdc,$
\begin{equation}
\begin{array}{ccc}
\Ci(K\col \xi \col \tautypes) &
           {\buildrel { \scriptstyle (I\otimes S)\after d(Q,\xi, \nu)}
\over{\hbox to 80pt{\rightarrowfill}}}  &
                \Ci(K\col \xi \col \tautypesp)\\
\downarrow & & \downarrow \\
\Ci(K\col \xi)_\types &
           {\buildrel {\scriptstyle\pi_{\xi, -\nu}(Z)}\over{\hbox to 80pt{\rightarrowfill}}} &
                \Ci(K\col \xi)_\typesp
\end{array}
\end{equation}
Here the vertical arrows represent the natural isomorphisms of
Lemma \refer{l: sphericalization}.
We denote
both of these isomorphisms by $\gf \mapsto \gf'.$
It suffices to prove the claim, since its validity implies
that $ \pi_{\xi, -\nu}(Z)$ maps $\Ci(K\col \xi)_\types$ into $\Ci(K\col \xi)_\typesp$
and that the expression on the right-hand side of (\refer{e: I otimes S on right}) equals
the one on the right-hand side of (\refer{e: intertwining prop Eis}).
Combining this with (\refer{e: I otimes S der F as L Z Eis})
we obtain (\refer{e: intertwining prop Eis}).

To see that the claim holds,
let $\gf \in \Ci(K\col \xi \col \tautypes) = (\Ci(K\col \xi) \otimes \Vtypes)^K.$
The associated element $\gf' \in \Ci(K \col \xi)_\types$ is given by
$$
\gf'(k) = \gf(k)(e),\qquad(k\in K).
$$
The element
$(I \otimes S_1)d(Q,\xi, \nu)\gf$ of
$ [\Ci(K\col \xi) \otimes C(K)_{\typesovee} \otimes \Vtypes]^K$
is given by
\begin{eqnarray*}
[(I \otimes S_1)d(Q,\xi, \nu)\gf](k)(k_1)
&=&  [d(Q, \xi, \nu)\gf(k)](\Ad(k_1^{-1})Z)\\
&=&
[I\otimes \pi_{\xi, -\nu}(\Ad(k_1^{-1}) Z) ] \gf (k);
\end{eqnarray*}
see (\refer{e: defi d Q xi nu}).
Hence, the element $(I\otimes S)d(Q,\xi,\nu)\gf \in [\Ci(K\col \xi) \otimes
\Vtypesp]^K$  is given by
$$
(I\otimes S)d(Q,\xi,\nu)\gf(k)(k_1) =
(I \otimes \pi_{\xi, -\nu}(\Ad(k_1^{-1}) Z))\gf (k)(k_1).
$$
The natural isomorphism from
$[\Ci(K\col \xi) \otimes \Vtypesp]^K$ onto $\Ci(K\col \xi)_\typesp$
is induced
by the map $I \otimes \gd_e,$ where $\gd_e: \Vtypesp \to \C$ denotes
evaluation at $e$ (see (\refer{e: defi Vtypes})).
Hence,
$$
((I\otimes S)d(Q,\xi,\nu)\gf)'(k) = [(I \otimes \pi_{\xi, -\nu}(Z)) \gf] (k)(e) =
[\pi_{\xi, -\nu}(Z) \gf'] (k).
$$
This establishes the claim.
\qed

We shall apply the above result in combination with Proposition
\refer{p: locally uniform generators} to obtain the assertion of
Theorem \refer{t: intertwining prop J} about finiteness. If $H
\subset \faQqdc$ is a $\gSr(Q)$-hyperplane, we denote by $\ga_H$
the shortest root of $\gSr(Q)$ such that $H$ is  a translate of
$(\ga_H^\perp)_\iC.$ Thus, $\inp{\ga_H}{\dotvar}$ equals a
constant $c$ on $H;$ we denote by $l_H$ the linear polynomial
function $\inp{\ga_H}{\dotvar} - c.$ In accordance with
\bib{BSpl1}, Eqn.~(\refbplc), given a locally finite collection $\Hyp$
of $\gSr(Q)$-hyperplanes in $\faQqdc$ and a map $d: \Hyp \to \N,$
we define, for every bounded subset $\omega$ of $\faQqdc,$ the
polynomial $\pi_{\omega, d}$ by
\begin{equation}
\naam{e: defi pi omega d}
\pi_{\omega, d}(\nu) =
\prod_{H \in \Hyp\atop H \cap \omega \neq \emptyset} l_H^{d(H)}
\end{equation}

\begin{prop}
\naam{p: uniformity singularities in types}
Let $Q \in \allparabs,$ $\xi \in \discserQ.$ Then $\Hyp(Q,\xi)$
is a locally finite collection of real $\gSr(Q)$-hyperplanes.
Moreover, there exists a map $d: \Hyp(Q,\xi) \to \N$ such
that, for every finite dimensional unitary representation $\tau$ of $K,$
every $T \in \barVxi \otimes \Ci(K\col \xi \col \tau) $ and
every bounded open subset $\omega \subset \faQqdc,$ the $\Ci(\spX\col \tau)$-valued
function
\begin{equation}
\naam{e: pi omega d times nE Q}
\nu \mapsto \pi_{\omega, d}(\nu) \nE(Q\col \psi_T \col \nu)
\end{equation}
is holomorphic on $\omega.$
Here $\pi_{\omega, d}$ is defined by (\refer{e: defi pi omega d}) with
$\cH = \cH(Q,\xi).$
\end{prop}
\proof Select any bounded open subset $\omega \subset \faQqdc.$
Let $\types \subset \dK$ be a finite set associated with $Q, \xi,
- \omega$ as in Proposition \refer{p: locally uniform generators}.
According to \bib{BSpl1}, Prop.~\refplu, there exists a map $d :
\Hyp(Q,\xi, \types) \to \N$ with the property that, for every $T
\in \barVxi \otimes \Ci(K\col \xi)_\types,$ the map $\nu \mapsto
\nEtypes(Q\col \psi_T \col \nu)$ belongs to $\Mer(\faQqd,
\Hyp(Q,\xi, \types), d, \Ci(\spX\col \tautypes)).$ See
\bib{BSpl1}, \S~\refplt, for the definition of the latter space.

Let $\typesp \subset \dK$ be an arbitrary finite subset.
Fix $\nu_0 \in \omega.$
Then
by Proposition \refer{p: locally uniform generators}
there exists $k \in \N$ such that the map
$$
M_\nu:\;\; U_k(\fg) \otimes C(K\col \xi)_\types \to C(K\col \xi)_K,
\;u\otimes \gf \mapsto \pi_{\xi, -\nu}(u) \gf
$$
has image containing $C(K\col \xi)_\typesp$ for $\nu = \nu_0.$
On the other hand, let $\typespp\subset \dK$ be the finite collection of $K$-types occurring in
$\gd_1 \otimes \gd_2,$ with $\gd_1\in \dK$ a $K$-type occurring in the adjoint representation
of $K$ in $U_k(\fg)$ and
with $\gd_2 \in \types.$
Then the image of $M_\nu$ is contained in $\Ci(K\col \xi)_\typespp$ for all $\nu \in \faQqdc.$
Let $\Ptypespp$ denote the $K$-equivariant projection from $C(K\col \xi)_{\typespp}$ onto
$C(K\col \xi)_\typesp.$ Then $\Ptypespp \after M_{\nu_0}$ is surjective.
Hence there exists a finite dimensional subspace $E \subset U_k(\fg) \otimes
C(K\col \xi)_\types$ such that
$R_\nu:= \Ptypespp \after M_{\nu}|_E: E \to C(K\col \xi)_\typesp$
is a bijection for $\nu = \nu_0.$
By continuity and finite dimensionality, there exists an open
neighborhood $\omega_0$ of $\nu_0$ in $\omega$ such that $R_\nu$ is a bijection
for all $\nu \in \omega_0.$ By Cramer's rule,
the inverse $S_\nu:= R_\nu^{-1} \in \Hom(C(K\col \xi)_\typesp, E)$ depends
holomorphically on $\nu \in \omega_0.$ Let $(u_i\mid 1 \leq i \leq I)$
be a basis of $U_k(\fg)$ and $(\gf_j \mid 1 \leq j \leq J)$ a basis
of $C(K\col \xi)_\types.$ Then there exist holomorphic
$[\Ci(K\col \xi)_\typesp]^*$-valued functions $s_{ij}$ on $\omega_0,$
for $1 \leq i \leq I,\; 1 \leq j \leq J,$ such that
$$
S_\nu(\gf) =
\sum_{1 \leq i \leq I\atop 1 \leq j \leq J} s_{ij}(\nu, \gf) u_i \otimes \gf_j,
\qquad  (\nu \in \omega_0),
$$
for $\gf \in \Ci(K \col \xi)_\typesp.$
Let $\gf \in C(K\col \xi)_\typesp.$
Then $\gf  = \Ptypespp \after M_\nu \after S_\nu (\gf),$
hence
$$
\gf = \sum_{i,j} s_{ij}(\nu, \gf) \, \Ptypespp \pi_{\xi, -\nu}(u_i)\gf_j.
$$
Let $\eta \in \barVxi.$ Then it follows from the above by application of
(\refer{e: psi T and Ptypes}) and (\refer{e: E psi and Ptypes}) that
\begin{eqnarray*}
\pi_{\omega, d}(\nu)
\nEtypesp(Q\col \psi_{\eta \otimes \gf} \col \nu)
&=&
\sum_{i,j} s_{ij}(\nu, \gf) \pi_{\omega,d}(\nu)
\nEtypesp(Q\col \psi_{\eta \otimes \Ptypespp \pi_{\xi, -\nu}(u_i)\gf_j}\col \nu)\\
&=&
\sum_{i,j} s_{ij}(\nu,\gf) (I\otimes \Ptypespp) [ \pi_{\omega,d}(\nu) \nEtypespp(Q\col \psi_{\eta \otimes
\pi_{\xi, -\nu}(u_i)\gf_j} \col \nu)].
\end{eqnarray*}
Applying $I \otimes \gd_e = \sphiso_\typesp^{-1}$ and using Lemma \refer{l: sphiso and rmi} and
Proposition \refer{p: left differentiation Eis}
we infer that
\begin{eqnarray}
\lefteqn{
\pi_{\omega, d}(\nu)
\nEtypesp(Q\col \psi_{\eta \otimes \gf} \col \nu \col \dotvar)(e)
}
\nonumber\\
&=&
\sum_{i,j} s_{ij}(\nu,\gf) \Ptypespp [ \pi_{\omega,d}(\nu) \nEtypespp(Q\col \psi_{\eta \otimes
\pi_{\xi, -\nu}(u_i)\gf_j} \col \nu\col \dotvar)(e)]
\nonumber\\
&=&
\sum_{i,j} s_{ij}(\nu, \gf) \Ptypespp L_{u_i}
[\pi_{\omega, d}(\nu) \nEtypes(Q\col \psi_{\eta \otimes \gf_j} \col \nu\col  \dotvar)(e)]
\naam{e: nEtypesp in nEtypes new}
\end{eqnarray}
{}From this we conclude that the expression on the left-hand side of the above
equation depends holomorphically on $\nu \in \omega_0,$ as a function
with values in $\Ci(\spX).$ Since $\nu_0$ was arbitrary,
it follows that the expression on the left-hand side of
(\refer{e: nEtypesp in nEtypes new}) in fact
depends holomorphically on $\nu \in \omega.$
Hence, every $H \in \Hyp(Q,\xi, \typesp)$ with
$ H \cap \omega \neq \emptyset$ must be contained in $\Hyp(Q,\xi, \types).$
This shows that the collection $\Hyp(Q,\xi)$ is locally finite.
The argument also shows that there exists a map $d: \Hyp(Q,\xi) \to \N$
such that the assertion of the proposition holds for every $\tau$
of the form $\tau = \tautypesp,$ with $\typesp \subset \dK$  a finite subset.
The general result now follows by application of the functorial property
of Lemma \refer{l: functoriality of nE}.
\qed
\begin{cor}
\naam{c: holomorphy J}
Let $d: \Hyp(Q,\xi) \to \N$ be as in
Proposition \refer{p: uniformity singularities in types}.
Then,  for every $T \in \barVxi \otimes \Ci(K\col \xi)_K$ and  every
bounded open subset $\omega \subset \faQqdc,$
the function
$$
\nu \mapsto \pi_{\omega,d}(\nu)J_{Q,\xi,\nu}(T)
$$
extends to a holomorphic
$\Ci(\spX)$-valued
function on $\omega.$
\end{cor}
\proof
This follows from Proposition \refer{p: uniformity singularities in types} and
the definition of $J_{Q,\xi, \nu},$ see (\refer{e: defi of J}).
\qed
We can now finally give the promised proof.
\medno
{\bf Proof of Theorem \refer{t: intertwining prop J}:\ }
The properties of $\Hyp(Q,\xi)$ have been established in
Proposition \refer{p: uniformity singularities in types} and
(\refer{e: intersection i faQqd with Hyp}). Let  $\nu \in \faQqdc.$
That $J_{Q,\xi,\nu}$ is a $\fg$-equivariant map follows from formula
(\refer{e: defi of J}) combined
with formula (\refer{e: intertwining prop Eis}) with $k=e.$
It remains to establish the $K$-equivariance of $J_{Q,\xi, \nu}.$
Let $\types \subset \dK$
be a finite subset and let $T = \eta \otimes \gf
\in \barVxi \otimes \Ci(K\col \xi)_\types.$
We denote the natural isomorphism $\Ci(K\col\xi)_\types \to \Ci(K\col \xi \col \tautypes)$
of
Lemma \refer{l: sphericalization}
by $\sphiso = \sphisotypes.$
Let $k_1 \in K.$ Then one readily checks that
$\sphiso \after \pi_{\xi, -\nu}(k_1) = (I \otimes S)\after \sphiso,$
with $S$ the endomorphism of $\Vtypes = C(K)_\typesvee$
given by restriction of the left translation $L_{k_1}.$ In particular,
$S$ intertwines $\tautypes = R|_{\Vtypes}$ with itself.
By the identification discussed in the text before
(\refer{e: psi T on types}) we have \begin{eqnarray*}
\psi_{[I \otimes \pi_{\xi, -\nu}(k_1)] T}
&=&
\psi_{[I \otimes \sphiso \pi_{\xi, -\nu}(k_1)]T}
\\
&=&
\psi_{(I\otimes I \otimes S)(I \otimes \sphiso)T}.
\end{eqnarray*}
By Lemma \refer{l: functoriality of nE} (a), combined with the identification
mentioned above,  the latter expression equals
$$
(I \otimes S)\psi_{(I \otimes \sphiso)T} = (I \otimes S) \psi_T.
$$
Applying Lemma \refer{l: functoriality of nE} (b) we now
find that
\begin{eqnarray*}
J_{Q,\xi,\nu}([I \otimes \pi_{\xi, -\nu}(k_1)] T) &=&
E^\circ_{\types}(Q\col
\psi_{[I \otimes  \pi_{\xi, - \nu}(k_1)]T} \col \nu)(\dotvar)(e)\\
&=&
[[I \otimes S]E^\circ_{\types}(Q\col \psi_{T} \col \nu)(\dotvar)](e) \\
&=&
\nEtypes(Q\col \psi_T \col \nu)(\dotvar)(k_1^{-1}) \\
&=&
L_{k_1} J_{Q,\xi,\nu}(T).
\end{eqnarray*}
\qed

\eqsection{The Fourier transform}
\naam{s: The Fourier transform}
Let $Q \in \allparabs$ and $\xi \in \discserQ.$ We will use the map $J_{Q,\xi},$
introduced in (\refer{e: defi of J}), to define a $(\fg, K)$-equivariant Fourier transform
for functions from $\Cci(\spX)_K.$

We define the collection $\Hypvee(Q,\xi)$ of hyperplanes in $\faQqdc$ by
$$
\Hypvee(Q,\xi): = \{ - H \mid H \in \Hyp(Q,\xi) \}.
$$
Since $\Hyp(Q,\xi)$ is a locally finite collection of real $\gSr(Q)$-hyperplanes,
see Theorem \refer{t: intertwining prop J},
the same holds for $\Hypvee(Q,\xi).$ It also follows from the mentioned theorem
that $\cup \Hypvee(Q,\xi)$ is disjoint from $i\faQqd.$

Since $\Hyp(Q,\xi)$ consists of real $\gS_r(Q)$-hyperplanes, every hyperplane
of $\Hyp(Q,\xi)$ is invariant under the complex conjugation $\gl \mapsto \bar \gl$
in $\faQqdc,$ defined with respect to the real form $\faQqd.$ Hence,
$\Hypvee(Q,\xi) = \{- \bar H \mid H \in \Hyp(Q,\xi)\}.$ This justifies
the following definition.

\begin{defi}
Let $f \in \Cci(\spX)_K.$ For
$\nu \in \faQqdc\setminus\cup\Hypvee(Q,\xi),$
the Fourier transform $\hat f (Q\col\xi\col\nu)$ is defined to be the element
of $\barVxi \otimes \Ci(K\col \xi)_K$ determined by
\begin{equation}
\naam{e: defi hat f}
\hinp{\hat f (Q\col \xi \col \nu)}{T} =
\int_\spX f(x) \overline{J_{Q,\xi, - \bar \nu}(T)(x)}\; dx,
\end{equation}
for all $T \in \barVxi \otimes \Ci(K\col \xi)_K.$
\end{defi}

\begin{lemma}
\naam{l: equivariance hat f}
Let $\nu \in \faQqdc\setminus \cup\Hypvee(Q,\xi).$
Then the map $f \mapsto \hat f(Q\col\xi\col\nu)$ from $\Cci(\spX)_K$ to
$\barVxi \otimes \Ci(K\col \xi)_K$ intertwines the $(\fg, K)$-module
structure on $\Ci(\spX)_K$ coming from the left regular representation
with the $(\fg, K)$-module structure on $\barVxi \otimes \Ci(K\col \xi)_K$
coming from $1 \otimes \pi_{Q, \xi, -\nu}.$
\end{lemma}

\proof
The spaces $\Cci(\spX)_K$ and $\barVxi \otimes \Ci(K\col \xi)$
are equipped with the natural $L^2$-type inner products. The first of these inner products
is equivariant for the $(\fg, K)$-module structure coming from the left regular representation.
The second induces a sesquilinear pairing of $\barVxi \otimes \Ci(K\col \xi)$ with itself
that is equivariant for the $(\fg, K)$-module structures coming from
$1 \otimes \pi_{Q,\xi, -\nu}$ and $1 \otimes \pi_{Q,\xi, \bar \nu},$ respectively.
On the other hand, it follows from (\refer{e: defi hat f})
that the map $f \mapsto \hat f (Q,\xi, \nu)$ is adjoint to the map
$J_{Q,\xi , - \bar \nu},$ with respect to the given inner products.
Therefore, the $(\fg, K)$-intertwining property of $f \mapsto \hat f (Q, \xi, \nu)$
follows by transposition
from the similar property for $J_{Q,\xi, -\nu},$ asserted in
Theorem \refer{t: intertwining prop J}.
\qed

If $d: \Hyp(Q,\xi) \to \N$ is a map, we define the map $\dvee: \Hypvee(Q,\xi) \to \N$ by
$\dvee(H) = d(-H),$ for $H \in \Hypvee(Q,\xi).$

\begin{lemma}
\naam{l: conjugate of pi omega d}
Let $d,\dvee$ be as above and let $\omega \subset \faQqdc$ be a bounded subset.
Then
\begin{equation}
\naam{e: conjugate of pi omega d}
\overline{\pi_{-\bar \omega, d}(-\bar \nu)} = (-1)^N \pi_{\omega, \dvee}(\nu),
\qquad (\nu \in \faQqdc),
\end{equation}
with $N = \sum_H d(H),$ for $H \in \Hyp(Q,v),$  $H \cap \omega \neq \emptyset.$
\end{lemma}
\medno
\proof
In the notation of the text preceding (\refer{e: defi pi omega d}),
we have, for $H \in \Hypvee(Q,\xi),$
\begin{equation}
\naam{e: overline l H}
\overline{l_H(-\bar \nu) } = - \inp{ \ga_H }{ -\nu } - c_H=
\inp{\ga_{-H}}{\nu} + c_{ -H } = - l_{ -H }(\nu).
\end{equation}
Moreover, since $H$ is real, $-\bar \omega \cap H \neq \emptyset$ is equivalent to
$\omega \cap (-H) \neq \emptyset.$
In view of the definition in (\refer{e: defi pi omega d}),
the  identity (\refer{e: conjugate of pi omega d}) follows from
(\refer{e: overline l H}) raised to the power $\dvee(H) = d(-H),$
by taking the product over all $H \in \Hypvee(Q,\xi)$ with $\omega\cap (-H) \neq \emptyset.$
\qed
\medno
\begin{lemma}
Let $d: \Hyp(Q,\xi) \to \N$ be as in Proposition
\refer{p: uniformity singularities in types}
and define $\dvee: \Hypvee(Q,\xi) \to \N$ by $\dvee(H) = d(-H).$
Then, for every bounded open subset $\omega \subset \faQqdc,$
every finite subset $\types \subset \dK$ and every $f \in \Ci(\spX)_\types,$
the function
\begin{equation}
\naam{e: pi times hat f}
\nu \mapsto \pi_{\omega,\dvee}(\nu) \hat f(Q, \xi, \nu),
\end{equation}
originally defined on $\omega\setminus \cup \Hypvee(Q,\xi),$ extends
to a holomorphic $\barVxi \otimes \Ci(K\col \xi)_\types$-valued function
on $\omega.$
\end{lemma}
\medno
\proof
Let $f$ be fixed as above.
It follows from Lemma \refer{l: equivariance hat f}
that the function (\refer{e: pi times hat f})
has values in the finite dimensional space $\barVxi \otimes \Ci(K\col \xi)_\types.$
Hence, it suffices to establish the holomorphic continuation
of the function that results
from (\refer{e: pi times hat f}) by taking the inner product
with a fixed element $T$ from $\barVxi \otimes \Ci(K\col \xi)_\types.$
In view of (\refer{e: defi hat f}) the resulting function equals
$$
\nu \mapsto \pi_{\omega, \dvee}(\nu) \hinp{f}{J_{Q,\xi, -\bar \nu}(T)}
=
\hinp{f}{\pi_{-\bar \omega, d}(-\bar \nu)J_{Q,\xi, -\bar \nu}(T)},
$$
see Lemma \refer{l: conjugate of pi omega d}. We may now apply
Corollary \refer{c: holomorphy J}, with $-\bar \omega$ in place of
$\omega,$ to finish the proof. \qed \medbreak The following result
gives the connection between the present Fourier transform and the
spherical Fourier transform, defined in \bib{BSpl1}, \S~\refple.

\begin{lemma}
\naam{l: first relation hat f and Fou F}
Let $\types\subset \dK$ be a finite subset,
let $f \in \Ci(\spX)_\types$ and let $F = \sphiso_\types(f) \in \Ci(\spX \col \tau_\types)$
be the corresponding spherical function, see Lemma \refer{l: sphericalization}.
Let $\cF_Q F$ be the
$\tautypes$-spherical Fourier transform of $F.$
Then, for every $T \in \barVxi\otimes  \Ci(K\col \xi)_\types ,$
$$
\hinp{\hat f(Q\col \xi \col \nu)}{T} = \hinp{\cF_Q F( \nu)}{\psi_T},
\qquad (\nu \in \faQqdc\setminus \cup \Hypvee(Q,\xi)).
$$
\end{lemma}

\proof
It follows from (\refer{e: defi hat f}) and (\refer{e: defi of J})
that
$$
\hinp{\hat f_\types (Q\col \xi \col \nu)}{T}
=\int_\spX f(x)
\overline{\nEtypes(Q\col \psi_T \col -\bar \nu)(x)(e)}\; dx.
$$
One may now proceed as in \bib{BSft}, p.~539, proof of Prop.~3, displayed
equations, but in reversed order.
\qed

\eqsection{A direct integral}
\naam{s: A direct integral}
In this section we assume that $Q \in \allparabs$ and $\xi \in \discserQ$
are fixed. We will define and study a direct integral representation $\pi_{Q,\xi}$
that will later appear as a summand in the Plancherel decomposition.

We equip
\begin{equation}
\naam{e: defi Hilb Q xi}
\Hilb(Q,\xi):=\barVxi \otimes  L^2(K \col \xi)
\end{equation}
with the tensor product Hilbert structure and the natural structure
of $K$-module. Moreover, we define
$$
\Ltwo(Q,\xi): = L^2(i\faQqd,\, \Hilb(Q,\xi),\, [W\col W_Q^*]\, d\mu_Q),
$$
the space of square integrable functions $i\faQqd \to
\Hilb(Q,\xi),$ equipped with the natural $L^2$-Hilbert structure
associated with the indicated measure. Recall that $d\mu_Q$ is
Lebesgue measure on $i\faQqd,$ normalized as in \bib{BSpl1}, end
of \S~\refpli.

By unitarity of the representations
$\pi_{Q, \xi, \nu},$ for $\nu \in i \faQqd,$ the prescription
$$
(\pi_{Q, \xi}(x) \gf)(\nu): =
[I \otimes \pi_{Q, \xi, -\nu}(x)] \gf(\nu),
\qquad (\gf \in \Ltwo(Q,\xi), \;x \in G),
$$
defines a homomorphism $\pi_{Q,\xi}$
from $G$ into the group $\rmU(\Ltwo(Q,\xi))$
of unitary transformations of $\Ltwo(Q,\xi).$
\begin{lemma}
\naam{l: pi Q xi continuous rep}
The homomorphism
$\pi_{Q,\xi}: G \to \rmU(\Ltwo(Q,\xi))$
is a unitary representation of $G$ in $\Ltwo(Q,\xi).$
\end{lemma}

\begin{rem}
It follows from the result above that $\pi_{Q,\xi}$ provides a
realization of
the following direct integral
$$
\pi_{Q,\xi} \simeq \int^\oplus_{i\faQqd}
  I_{\barVxi}\otimes \pi_{Q,\xi, -\nu}  \; [W:W_Q^*]\,d\mu_Q(\nu).
$$
\end{rem}

For the proof of Lemma \refer{l: pi Q xi continuous rep}
it is convenient to define
a dense subspace of $\Ltwo(Q,\xi)$ by
\begin{equation}
\naam{e: defi Ltwoc Q Omega}
\Ltwosub(Q,\xi): =  C_c(i\faQqd, \barVxi \otimes \Ci(K\col \xi)).
\end{equation}
This space
is equipped with a locally convex topology
in the usual way.
Thus, if $\cA \subset i\faQqd$ is compact,
let
$
\Ltwo_{0\cA}(Q,\xi)
$
denote the space $C_\cA(i\faQqd, \barVxi \otimes \Ci(K\col \xi))$ of continuous functions
from (\refer{e: defi Ltwoc Q Omega})
with support contained in $\cA.$ This space is equipped with the
Fr\'echet topology determined by the seminorms
$$
\gf \mapsto \sup_{\nu \in \cA} s(\gf(\nu)),
$$
where $s$ ranges over the continuous seminorms on $\barVxi \otimes \Ci(K\col \xi).$
Moreover,
$$
\Ltwosub(Q,\xi) = \cup_{\cA}\;\; \Ltwo_{0\cA}(Q,\xi)
$$
is equipped with the direct limit locally convex topology.
Thus topologized,  $\Ltwosub(Q,\xi)$ is a complete locally convex space.
\begin{lemma}
\naam{l: smooth rep on Ltwosub}{\ }
\begin{enumerate}
\itema
The space $\Ltwosub(Q,\xi)$ is $G$-invariant.
\itemb
The restriction of $\pi_{Q,\xi}$ to  $\Ltwosub(Q,\xi)$ is
a smooth representation of $G.$
Moreover, if $\gf \in \Ltwosub(Q,\xi)$ and $u \in U(\fg),$ then
$\pi_{Q,\xi}(u)\gf$ is given by
\begin{equation}
\naam{e: pi of u on gf}
[\pi_{Q,\xi}(u) \gf]( \nu) =
[I \otimes \pi_{Q,\xi, -\nu}(u)] \; \gf(\nu).
\end{equation}
\end{enumerate}
\end{lemma}

\proof
Let $\cA \subset i \faQqd$ be compact. Then it
 is a straightforward consequence of the definition that the
space $\LtwosubA(Q,\xi)$ is $G$-invariant. In particular, this
implies (a).

For (b) it suffices to prove that the restriction of
$\pi_{Q,\xi}$ to $\LtwosubA(Q,\xi)$ is smooth
and that (\refer{e: pi of u on gf}) holds for $\gf \in \LtwosubA(Q,\xi).$

Let $\gf \in \LtwosubA(Q,\xi).$ We consider the function
$\Phi: \cA \times G \times K \to \barVxi \otimes \Hxii$ defined by
$$
\Phi(\nu, x, k) : = [I \otimes \ev_k]([\pi_{Q,\xi}(x) \gf](\nu)) =
[I \otimes \ev_k\pi_{Q,\xi, -\nu}(x)] \gf(\nu),
$$
where $\ev_k$ denotes the map $\Ci(K\col \xi) \to \cH_\xi^\infty$ induced
by evaluation in $k.$ We recall
that the multiplication map $N_Q \times A_Q \times \exp (\fm_Q\cap \fp) \times K \to G$
is a diffeomorphism. Accordingly, we write
\begin{equation}
\naam{e: Q K deco of x}
x = n_Q(x) a_Q(x) m_Q(x) k_Q(x), \qquad (x \in G),
\end{equation}
where $\nu_Q,$ $a_Q,$ $m_Q$ and $k_Q$ are smooth maps from $G$
into $N_Q,$ $A_Q,$ $\exp(\fm_Q \cap \fp)$ and $K_Q,$
respectively. Applying (\refer{e: Q K deco of x}) with $kx$
in place of $x,$ we find that
$$
\Phi(\nu, x, k) = a_Q(kx)^{-\nu + \rho_Q} \xi(m_Q(kx)) \gf(\nu, k_Q(kx)).
$$
{}From this expression we see that the
function $\Phi$ is continuous, and smooth in the variables $(x,k).$
Moreover, for every $C^\infty$ differential operator
$D$ on $G \times K,$ the $ \barVxi\otimes \Hxii $-valued function
$$
(\nu , x, k) \mapsto D[\Phi(\nu)](x,k)
$$
is continuous. This implies that the
$C_\cA(i\faQqd,  \barVxi\otimes \Ci(K\col \xi) )$-valued
function $x\mapsto \Phi(\cdot, x, \cdot)$ is smooth on $G;$
hence $\gf$ is a smooth vector for $\pi_{Q,\xi}.$
Let $D$ be any $C^\infty$ differential operator
on $G$ and let $\nu \in \cA.$ Then evaluation in $\nu$ induces a continuous
linear operator $\Ltwo_{0\cA} \to \barVxi \otimes \Ci(K\col\xi).$ Hence,
for all $x \in G,$
$$
D[\pi_{Q,\xi}(\dotvar) \gf](x)(\nu) = D[ \ev_\nu(\pi_{Q,\xi}(\dotvar)\gf)](x)
= D(\pi_{Q,\xi, -\nu}(\dotvar)\gf(\nu))(x).
$$
Applying this with $D = R_u$ and $x =e$ we obtain (\refer{e: pi of u on gf}).
\qed
\medno
{\bf Proof of Lemma \refer{l: pi Q xi continuous rep}:\ }
Put $\pi := \pi_{Q,\xi}.$
It suffices to show that the map $G \times \Ltwo(Q,\xi) \to \Ltwo(Q,\xi),$
$(x,\gf) \mapsto \pi(x) \gf$ is continuous.
Since $\pi$ is a homomorphism from the group $G$ into $\rmU(\Ltwo(Q,\xi)),$
it suffices to show that for
any fixed  $\gf \in \Ltwo(Q,\xi)$ we have
\begin{equation}
\naam{e: lim pi x gf new}
\lim_{x \to e}\pi(x)\gf = \gf \text{in} \Ltwo(Q,\xi).
\end{equation}
Moreover, again because $\pi$ maps into $\rmU(\Ltwo(Q,\xi)),$ it suffices
to prove (\refer{e: lim pi x gf new}) for $\gf$ in a
dense subspace of $\Ltwo(Q,\xi).$
Let $\gf \in \Ltwosub(Q,\xi).$ Then $\pi(x) \gf \to \gf$ in $\Ltwosub(Q,\xi),$
as $x \to e,$
by Lemma \refer{l: smooth rep on Ltwosub}. By continuity
of the inclusion map $\Ltwosub(Q,\xi) \embeds \Ltwo(Q,\xi),$
this implies (\refer{e: lim pi x gf new}).
\qed

We end this section by establishing some other useful properties
of the invariant subspace $\Ltwosub(Q,\xi).$
If $\types\subset \dK,$ then one readily verifies that
\begin{equation}
\naam{e: Ltwosub types}
\Ltwosub(Q,\xi)_\types  = C_c(i\faQqd, \barVxi \otimes \Ci(K\col \xi)_\types).
\end{equation}
The space of $K$-finite vectors in $\Ltwosub(Q,\xi)$ equals the union
of the spaces (\refer{e: Ltwosub types}) as $\types$ ranges over the collection of
finite subsets
of $\dK.$ The natural $U(\fg)$-module structure of $\Ltwosub(Q,\xi)_K$ is given
by formula (\refer{e: pi of u on gf}).

\begin{lemma}
\naam{l: density in smooth vectors two}
Let $(\rho,W)$ be a continuous representation
of $G$ in a complete locally convex space, and let $U$ be a dense
$G$-invariant subspace of $W.$ If $U$ is contained in $W^\infty$ then it
is dense  in $W^\infty$
for the $C^\infty$-topology.
\end{lemma}

\begin{rem}
For $W$ a Banach space, this result is Thm.~1.3 of \bib{Poul}. The following proof
is an adaptation of the proof given in the mentioned paper.
\end{rem}

\proof
Replacing $U$ by its closure in $W^\infty$ if necessary, we may
as well assume that $U$ is closed in $W^\infty.$
Fix a choice $dg$ of Haar measure on $G.$
If $\gf \in \Cci(G),$ then the map
$$
\rho(\gf) := \int_G \gf(g)\rho(g)\; dg
$$
is continuous linear from $W$ to $W^\infty,$ as can be seen from a straightforward
estimation. Moreover, since $U$ is closed in $W^\infty,$ the map $\rho(\gf)$ maps
$U$ into itself.
Let $W_1$ be the collection of vectors in $W^\infty$ of the
form $\rho(\gf) w_0,$ with $\gf \in \Cci(G)$ and $w_0 \in W^\infty.$
Then $W_1$ is dense in $W^\infty.$ Hence, it suffices to show
that $W_1$ is contained in $U.$
Select $w_1 \in W_1$ and let $N_1$ be an open neighborhood of $0$ in $W^\infty.$
Then it suffices to show that $U \cap (w_1 + N_1) \neq \emptyset.$

Write $w_1 = \rho(\gf) w_0,$ with $w_0 \in W^\infty$ and $\gf \in \Cci(G).$
By the mentioned continuity of $\rho(\gf),$ there exists an open
neighborhood $N_0$ of $0$ in $W$ such that $\rho(\gf)N_0 \subset N_1.$
By density of $U$ in $W$, the intersection $U \cap (w_0 + N_0)$ is non-empty.
Hence,
\begin{eqnarray*}
\emptyset &\subsetneq& \rho(\gf)[U \cap (w_0 + N_0)]\\
& \subset &
U \cap \rho(\gf)(w_0 + N_0)
\\
&\subset& U \cap (w_1 + N_1).
\end{eqnarray*}
\qed

\begin{lemma}
\naam{l: density in Ltwoinfty}
The space $\Ltwosub(Q,\xi)_K$ is dense in $\Ltwo(Q,\xi)^\infty$
with respect to the
natural Fr\'echet topology of the latter space.
\end{lemma}

\proof
The inclusion map $j: \Ltwosub(Q,\xi) \to \Ltwo(Q, \xi)$
is continuous, intertwines the $G$-actions and has dense image.
{}From Lemma \refer{l: smooth rep on Ltwosub} it follows
that $\Ltwosub(Q,\xi)^\infty = \Ltwosub(Q,\xi).$
By equivariance of $j$ it follows that $\Ltwosub(Q,\xi)$ is contained in
$\Ltwo(Q,\xi)^\infty.$ By application of Lemma
\refer{l: density in smooth vectors two} we
see that $\Ltwosub(Q, \xi)$ is dense in $\Ltwo(Q,\xi)^\infty.$
The conclusion now follows since $\Ltwosub(Q, \xi)_K$ is dense
in $\Ltwosub(Q, \xi).$
\qed

\eqsection{Decomposition of the regular representation}
\naam{s: Decomposition of the regular representation}
Up till now, for $Q \in \allparabs,$ the expression $\xi \in \discserQ$
meant, by abuse of language, that $\xi$ is an irreducible unitary
representation of $M_Q$ with equivalence class $[\xi]$ contained in $\discserQ.$
{}From now on it will be convenient to distinguish between representations
and their classes.
For every $Q \in \allparabs$ and all $\omega \in \discserQ$ we select
an irreducible unitary representation $\xi = \xi_\omega$ in a Hilbert space
$\cH_\omega:= \Hxi,$ with $[\xi] = \omega.$ Moreover, we put
$\fH(Q,\omega) :=
\fH(Q, \xi_\omega),$
see (\refer{e: defi Hilb Q xi}).

For $\types \subset \dK$ a finite subset, let $\discserQ(\types)$
denote the collection of $\omega \in \discserQ$ that
have a $K_Q$-type in common with $\tau_\types.$

\begin{lemma}
\naam{l: finiteness discserQtypes}
Let $Q \in \allparabs$ and let
 $\types \subset \dK$ be a finite subset. Then $\discserQ(\types)$
is finite. Moreover,
$$
\cA_{2,Q}(\tautypes) = \oplus_{\omega \in \discserQ(\types)}\;\; \cA_{2,Q}(\tautypes)_\omega,
$$
where the direct sum is finite and orthogonal.
\end{lemma}
\proof
The collection $\discserQ(\types)$ is finite by
Lemma \refer{l: finite ds with K type},
applied to the spaces $\spXQv,$
for $v \in \QcW.$ Put $\tau = \tautypes$ and fix $v \in \QcW.$
We note that $\cA_2(Q,v\col \tau_Q)_\omega = 0$
for $\omega \in \discserQ(\types)\setminus M^\wedge_{v H v^{-1}, ds}(\tau_Q),$
by Remark  \refer{r: types in cAtwo}, with $\spXQv,\tau_Q$ in place of $\spX, \tau.$
Moreover, by the same remark,
$$
\cA_2(\spX_{Q,v}\col \tau_Q) = \oplus_{\omega \in \discserQ(\types)} \;\;\;
\cA_2(\spXQv\col \tau_Q)_\omega,
$$
with orthogonal summands. The result now follows by summation over
$v \in \QcW,$ in view of (\refer{e: defi cAtwoQ xi}), and
\bib{BSpl1}, Eqn.~(\refplp). \qed

\begin{lemma}
\naam{l: fH Q omega types non trivial}
Let $Q \in \allparabs,$ $\types\subset \dK$  a finite subset and  $\omega \in \discserQ.$
Then $\fH(Q,\omega)_\types \neq 0$ if and only if $\omega \in \discserQ(\types).$
\end{lemma}
\proof
We have that $\fH(Q,\omega)_\types =
\bar V(Q,\xi_\omega) \otimes L^2(K\col \xi_\omega)_\types,$
with non-trivial first component in the tensor product. Hence, $\fH(Q,\omega)_\types$
is non-trivial if and only if $L^2(K\col \xi_\omega)_\types$ is.
Since $L^2(K\col \xi_\omega)$ is the representation space
for $\Ind_{K_Q}^K(\xi_\omega),$ the assertion
follows by Frobenius reciprocity.
\qed
If $Q \in \allparabs,$ we define the Hilbert space
\begin{equation}
\naam{e: fH Q as dir sum}
\fH(Q): = \widehat\oplus_{\omega \in \discserQ} \;\;\fH(Q, \omega),
\end{equation}
where the hat over the direct sum symbol indicates that the natural
Hilbert space completion of the algebraic direct sum
is taken.
Let $\types \subset \dK$ be a finite subset.
In view of Lemmas \refer{l: finiteness discserQtypes} and
\refer{l: fH Q omega types non trivial} it follows
that we may define a map
$$
\Psi_\types= \Psi_{Q,\types}:\;\; \fH(Q)_\types \to
\cAtwoQ(\tautypes)
$$
as the direct sum of the isometries
$T \mapsto \psi_T: \; \fH(Q, \omega)_\types \to  \cAtwoQ(\tautypes)_\omega,$
for $\omega \in \discserQ(\types),$ see
(\refer{e: psi T on types}). The following result is immediate.
\begin{lemma}
\naam{l: Psi types isometry}
$\Psi_\types$ is an isometric isomorphism from $\fH(Q)_\types$ onto
$\cAtwoQ(\tautypes).$
\end{lemma}
\medbreak
If $f \in \Cci(\spX)_K$ then for $Q \in \allparabs, $ $ \omega \in \discserQ$ and
$\nu \in i\faQqd$ we define
the Fourier transform $\hat f(Q\col \omega\col \nu) \in
\bar V(\xi_\omega)\otimes
\Ci(K\col\xi_\omega)_K $
by
$$
\hat f(Q\col \omega\col \nu):= \hat f(Q\col \xi_\omega \col \nu).
$$
This definition is justified since it follows from
(\refer{e: intersection i faQqd with Hyp})
that $\nu \notin \cup \Hypvee(Q,\xi).$

\begin{prop}
\naam{p: L two norms of Fou and fou}
Let $\types \subset \dK$ be a finite set of $K$-types, let
$f \in \Cci(\spX)_\types$ and let $F = \sphiso_\types(f) \in \Cci(\spX\col \tautypes)$ be
the associated spherical function, see Lemma \refer{l: sphericalization}.
 Then for each $Q \in \allparabs$ and
every $\nu \in i \faQqd,$
\begin{equation}
\naam{e: Fou Q as sum of psi fou}
\cF_Q F(\nu) = \sum_{\omega \in \discserQ} \psi_{\hat f(Q\col \omega\col \nu)}\quad
\text{in} \quad \cAtwoQ(\tautypes).
\end{equation}
If $\hat f(Q\col \omega\col \nu)$ is non-zero, then $\omega \in \discserQ(\types);$ in particular,
the above sum is finite. Finally,
$$
\|\cF_Q F(\nu )\|^2 =
\sum_{\omega \in \discserQ}
 \| \hat f(Q\col \omega \col \nu)\|^2_{\barVxiomega \otimes L^2(K\col \xi_\omega)}.
$$
\end{prop}

\proof
It follows from Lemma \refer{l: equivariance hat f} that
$\hat f(Q\col\omega\col \nu)$ is an element of
$\fH(Q, \omega)_\types,$ for every $\omega \in \discserQ.$ Hence, if this
element
is non-zero, then $\omega$ belongs to the finite set $\discserQ(\types).$

The identity (\refer{e: Fou Q as sum of psi fou}) follows from
Lemma \refer{l: first relation hat f and Fou F},
since $\Psi_\types$
is a surjective isometry. The final assertion follows by once more using that
$\Psi_\types$ is an isometry.
\qed

The following result will turn out to be the Plancherel identity
for $K$-finite functions. We recall from \bib{BSpl1},
Def.~\refplq, that two parabolic subgroups $P,Q \in \allparabs$
are said to be associated if their $\gs$-split components $\faPq$
and $\faQq$ are conjugate under $W.$ The equivalence relation of
associatedness on $\allparabs$ is denoted by $\sim.$ Let
$\repparabs\subset \allparabs$ be a choice of representatives for
$\allparabs/\sim.$ If $Q \in \allparabs,$ then $W^*_Q$ denotes the
normalizer of $\faQq$ in $W.$

\begin{thm}
\naam{t: FOU is isometry}
Let $f \in \Cci(\spX)_K.$ Then
$$
\|f\|_2^2 =
\sum_{Q \in \repparabs}\sum_{\omega \in \discserQ} [W:W_Q^*]\;
\int_{i \faQqd} \|\hat f (Q\col \omega \col \nu)\|^2\; d\mu_Q(\nu).
$$
\end{thm}

\proof
This follows from
\bib{BSpl1},
Thm.~\refplr, combined with Proposition \refer{p: L two norms of
Fou and fou}. \qed

Our next goal is to show that the above indeed corresponds with
a direct integral decomposition for the left regular representation
$L$ of $G$ in $L^2(\spX).$

Let $Q \in \repparabs.$ For $\omega \in \discserQ,$ the direct integral
representation
$\pi_{Q,\omega} : = \pi_{Q, \xi_\omega}$ of $G$ in
$\Ltwo(Q,\omega): = \Ltwo(Q, \xi_\omega)$ is unitary, see
Lemma \refer{l: pi Q xi continuous rep}.
We define
\begin{equation}
\naam{e: defi Ltwo}
\Ltwo(Q) := \widehat\oplus_{\omega \in \discserQ} \;\Ltwo(Q,\omega),
\end{equation}
where the hat over the direct sum sign has the same meaning
as in (\refer{e: fH Q as dir sum}).
We note that $\Ltwo(Q)$ is naturally  isometrically
isomorphic with the Hilbert space of $\fH(Q)$-valued $L^2$-functions on $i\faQqd,$
relative to the measure $[W\col W_Q^*]\, d \mu_Q.$
Let $\pi_Q$ be the associated direct sum of
the representations $\pi_{Q,\omega}.$ Then $\pi_Q$ is a unitary representation
of $G$ in $\Ltwo(Q).$

Finally, we define
\begin{equation}
\naam{e: defi Ltwo as dir sum}
\Ltwo: = \oplus_{Q \in \repparabs} \;\;  \Ltwo(Q)
\end{equation}
and equip it with the direct sum inner product. The direct sum being finite,
$\Ltwo$ becomes a Hilbert space in this way.
The associated direct sum $\pi = \oplus_{Q \in \repparabs} \pi_Q$ is
a unitary representation of $G$ in $\Ltwo.$

For $Q \in \repparabs$ and $\omega \in \discserQ$
we denote the natural inclusion map $\Ltwo(Q,\omega) \to \Ltwo$ by $\rmi_{Q,\omega};$
its adjoint $\pr_{Q,\omega}: \Ltwo \to \Ltwo(Q,\omega)$ is the natural projection map.
If $\gf \in \Ltwo,$ we denote its component $\pr_{Q,\omega} \gf \in \Ltwo(Q\col \omega)$
by
$\gf(Q\col \omega \col \dotvar).$
Thus,
$$
\|\gf\|^2_{\Ltwo} = \sum_{Q \in \repparabs} \sum_{\omega \in \discserQ}
[W:W_Q^*]\;
\int_{i\faQqd}
\|\gf(Q\col \omega\col \nu)\|^2 \; d  \mu_Q.
$$
It follows from Theorem \refer{t: FOU is isometry}
that the Fourier transform $\hat f$ of
an element $f \in \Cci(\spX)_K$ belongs to $\Ltwo.$ Moreover,
$$
\|f\|_{L^2(\spX)} = \|\hat f\|_{\Ltwo}.
$$

\begin{thm}
\naam{t: G equivariance FOU}
The map $f \mapsto \hat f$ has a unique extension
to a continuous linear map $\FOU: L^2(\spX) \to \Ltwo.$
The map $\FOU$ is isometric and intertwines the $G$-representations
$(L, L^2(X))$ and
$(\pi, \Ltwo).$
\end{thm}

\proof
The first two assertions are obvious from the discussion preceding
the theorem. It remains to prove the intertwining
property. Fix $Q \in \repparabs$ and $\omega \in \discserQ;$ then it suffices
to prove that $\FOU_{Q,\omega}:= \pr_{Q,\omega} \after \FOU$ intertwines
$L$ with $\pi_{Q,\omega}.$ We will do this in a number of steps.
For convenience, we write $\xi = \xi_\omega.$

\begin{lemma}
\naam{l: infi intertw prop FOU Q omega}
Let $f \in \Cci(\spX).$
\begin{enumerate}
\itema
If $k \in K,$ then $\FOU_{Q,\omega}(L_k f) = \pi_{Q,\omega}(k) \FOU_{Q,\omega}f.$
\itemb
If $u \in U(\fg),$ then
\begin{equation}
\naam{e: Ug intertwining FOU}
\hinp{\FOU_{Q,\omega}(L_u f)}{\gf} = \hinp{\FOU_{Q,\omega}(f)}{\pi(\check u )\gf}
\end{equation}
for all $\gf \in \Ltwo(Q,\omega)^\infty.$
\end{enumerate}
\end{lemma}

\proof
We first assume that $f$ is $K$-finite.
Assertion (a) is an immediate consequence of the $K$-equivariance asserted in
Lemma \refer{l: equivariance hat f}.

In view of Lemma \refer{l: density in Ltwoinfty}
it suffices to prove assertion (b) for $\gf \in \Ltwosub(Q,\omega).$
In view of Lemma \refer{l: smooth rep on Ltwosub} we may as well assume that
$u = X \in \fg.$
Then the expression on the left-hand side of (\refer{e: Ug intertwining FOU})
equals
$$
\int_{i\faQqd} \hinp{(L_X f)\hat{\,}(Q\col \xi\col \nu)}{\gf(\nu)}_{\fH(Q,\omega)}
 \; d\mu_Q(\nu).
$$
The integrand is continuous and compactly supported
as a function of $\nu.$ By the $\fg$-equivariance asserted in
Lemma \refer{l: equivariance hat f} and
unitarity of the representations $\pi_{Q,\xi, -\nu}$ for all
$\nu \in i\faQqd,$ we see that
the integral equals
\begin{eqnarray*}
\lefteqn{\!\!\!\!\!\!\!\!\!\!\!\!\!\!\!\int_{i\faQqd}
\hinp{\hat f(Q\col \xi \col \nu)}{[I \otimes \pi_{Q,\xi, -\nu}(\check X) ]\gf(\nu)}
\; d\mu_Q(\nu)\;\;\;\;\;\;\;\;\;\;\;\;\;\;\;}\\
&\qquad\qquad\qquad\qquad=&
\int_{i\faQqd}
\hinp{\hat f(Q\col \xi \col \nu)}{[\pi_{Q,\omega}(\check X)\gf](\nu)}\; d\mu_Q(\nu),
\end{eqnarray*}
see also Lemma \refer{l: smooth rep on Ltwosub}.
The  expression on the right-hand side of the latter equality
equals the one on the right-hand side of (\refer{e: Ug intertwining FOU}).
This establishes the result for $f$ in the dense subspace $\Cci(\spX)_K$ of $\Cci(\spX).$
The idea is to extend the result by an argument involving continuity.

For assertion (a) we proceed as follows. Fix $k \in K.$ Then the map
$f \mapsto L_k f$ is continuous from $\Cci(\spX)$ to $L^2(\spX).$ Since $\Fou_{Q,\omega}$
is continuous $L^2(\spX) \to \Ltwo(Q,\omega),$ by the first part of the proof
of Theorem \refer{t: G equivariance FOU}, whereas $\pi_{Q,\omega}$ is a unitary
representation, it follows that both $f \mapsto \Fou_{Q,\omega} L_k f$
and $f \mapsto \pi_{Q,\omega}(k) \Fou_{Q,\omega} f$ are continuous maps
from $\Cci(\spX)$ to $\Ltwo(Q,\omega).$ Hence, (a) follows by continuity and
density.

Finally, for the proof of (b) we fix $u \in U(\fg)$ and $\gf \in \Ltwo(Q,\omega)^\infty.$
Then the map $f \mapsto L_u f$ is continuous from $\Cci(\spX)$ to $L^2(\spX),$
whereas $\Fou_{Q,\omega}$ is continuous from $L^2(\spX)$ to $\Ltwo(Q,\omega)$
as said above.
It follows that the inner product on the left-hand side of
(\refer{e: Ug intertwining FOU}) depends
continuous linearly on $f \in \Cci(\spX).$ Since $\pi_{Q,\omega} (\check u) \gf$
is a fixed element of $\Ltwo(Q,\omega),$ the same holds for the inner product
on the right-hand side of (\refer{e: Ug intertwining FOU}).
Thus, (b) follows by continuity and
density from the similar statement for $K$-finite functions.
\qed

\begin{lemma}
\naam{l: equivariance FOU Q omega}
Let $f \in \Cci(\spX).$ Then
$$
\FOU_{Q,\omega}( L_x f) = \pi_{Q, \omega}(x) \FOU_{Q,\omega}( f )
$$
for all $x \in G.$
\end{lemma}

\proof
By Lemma \refer{l: infi intertw prop FOU Q omega} (a) it suffices to prove
the identity for $x$ in the connected component of $G$ containing $e.$
Hence it suffices to establish the identity for $x \in \exp(\fg).$
Write $\pi = \pi_{Q,\omega}$ and fix
$X \in \fg.$ Then it suffices to show
that $\pi(\exp t X)^{-1} \FOU_{Q,\omega} (L_{\exp tX} f)$ is a constant function
of $t \in \R$ with value $\FOU_{Q,\omega}(f).$ For this it suffices to show that,
for every $\gf \in \Ltwo(Q,\omega)^\infty,$ the function
$$
\Phi: t \mapsto \hinp{\pi(\exp t X)^{-1} \FOU_{Q,\omega} (L_{\exp tX} f)}{\gf}
$$
is differentiable with derivative zero.
We observe
that
$$
\Phi(t) = \hinp{\FOU_{Q,\omega}(L_{\exp tX} f)}{\pi(\exp tX) \gf}.
$$
The $L^2(\spX)$-valued function $t \mapsto L_{\exp tX} f$ is $C^1$ on $\R,$
with derivative $t \mapsto L_X L_{\exp tX}f.$
Moreover, since $\gf$ is a smooth vector,
the $\Ltwo(Q,\omega)$-valued function $t \mapsto \pi(\exp tX) \gf$
is also $C^1$ on $\R,$ with derivative
$t \mapsto \pi(X)\pi(\exp tX).$
By continuity of $\FOU_{Q,\omega}$ and the inner product
on $\Ltwo(Q,\omega),$ it follows that $\Phi(t)$ is $C^1$ with derivative given by
$$
\Phi'(t) = \hinp{\FOU_{Q,\omega}(L_X L_{\exp tX} f)}{\pi(\exp tX) \gf} +
\hinp{\FOU_{Q,\omega}(L_{\exp tX} f)}{\pi(X) \pi(\exp tX) \gf}.
$$
The latter expression equals zero by Lemma
\refer{l: infi intertw prop FOU Q omega} (b), applied with
$L_{\exp tX} f$ and $\pi(\exp tX) \gf$ in place of $f$ and $\gf,$
respectively.
\qed

{\bf End of proof of Theorem \refer{t: G equivariance FOU}:\ }
In the beginning of the proof of the theorem, we established that $\Fou_{Q,\omega}$
is a continuous linear map from $L^2(\spX)$ to $\Ltwo(Q,\omega).$
By density of $\Cci(\spX)_K$ in $L^2(\spX)$ and continuity  of
the representations $L$ and $\pi_{Q,\omega},$ the $G$-equivariance of
$\FOU_{Q,\omega}$ follows from Lemma \refer{l: equivariance FOU Q omega}.
\qed

Let $Q \in \repparabs$ and let $\types \subset \dK$ be a finite
subset. We recall from \bib{BSpl1}, Thm.~\refplr, that the
spherical Fourier transform $\Fou_Q$ associated with $\tau =
\tau_\types,$ originally defined as a continuous linear map
$\cC(\spX\col \tautypes) \to \cS(i\faQqd) \otimes
\cAtwoQ(\tautypes),$ has a unique extension to a continuous linear
map $L^2(\spX\col \tautypes) \to L^2(i\faQqd) \otimes
\cAtwoQ(\tautypes),$ denoted by the same symbol. For application
in the next section, we state the relation between the extended
spherical Fourier transform $\Fou_Q$ and $\FOU$ in a lemma. Let
$\pr_Q$ denote the projection operator  $\Ltwo \to \Ltwo(Q)$
associated with the  decomposition (\refer{e: defi Ltwo as dir
sum}).

\begin{lemma}
\naam{l: Fou Q and FOU}
Let $Q \in \repparabs$ and let $\types\subset \dK$ be a finite subset.
Let $f \in L^2(\spX)_\types$ and let $F = \sphiso_\types(f) \in L^2(\spX\col \tautypes)$
be
the associated spherical function, see Lemma \refer{l: sphericalization}.
Then
\begin{equation}
\naam{e: Fou Q and FOU}
\Fou_Q F (\nu) = \Psi_{Q,\types}(\pr_Q \after \FOU f(\nu)),
\end{equation}
for almost all $\nu \in i\faQqd.$
Here $\Psi_{Q,\types}$ is the isometry of  Lemma \refer{l: Psi types isometry}.
\end{lemma}

\proof For $f \in \Cci(\spX)_\types$ we have $\pr_Q\after \FOU f
(\nu) = \hat f(Q\col \dotvar\col \nu),$ so that (\refer{e: Fou Q
and FOU}) follows from (\refer{e: Fou Q as sum of psi fou}) and
the definition of $\Psi_{Q,\types}$ before Lemma \refer{l: Psi
types isometry}. The general result follows from this by density
of $\Cci(\spX)_\types$ in $L^2(\spX)_\types$ and continuity of the
maps $\Fou_Q\after \sphiso_\types$ and $(I \otimes
\Psi_{Q,\types})\after \pr_Q \after \FOU$ from $L^2(\spX)_\types$
to $L^2(i\faQqd) \otimes \cAtwoQ(\tau_\types),$ see \bib{BSpl1},
Thm.~\refplr, and Lemmas \refer{l: sphericalization}, \refer{l:
Psi types isometry} and Theorem \refer{t: G equivariance FOU}.
\qed

\eqsection{The Plancherel decomposition}
\naam{s: The Plancherel decomposition}
Our goal in this section is to establish the Plancherel decomposition.
For this we need to characterize
the image of the transform $\FOU,$ defined in the previous section.
To achieve this we shall first decompose
$\FOU$ into pieces corresponding with the parabolic subgroups from $\repparabs.$

Let $Q \in \repparabs.$ For $\types \subset \dK$ a finite subset,
we defined in \bib{BSpl1}, text before Thm.~\refpla{}  a subspace
$\cC_Q(\spX\col \tau_\types)$ of $\cC(\spX\col \tau_\types),$ as
the image of the wave packet transform $\Wave_Q.$ In \bib{BSpl1},
text before Cor.~\refbpld, we defined $L^2_Q(\spX\col \tautypes)$
as the closure of $\cC_Q(\spX\col \tau_\types)$ in $L^2(\spX\col
\tautypes).$ Accordingly, we denote by $L^2_Q(\spX)_\types$ the
canonical image of $L^2_Q(\spX\col \tautypes)$ in
$L^2(\spX)_\types,$ cf.\ Lemma \refer{l: sphericalization}.
Finally, we denote by $L^2_Q(\spX)$ the $L^2$-closure of the union
of such spaces for all $\types.$ Then it follows from \bib{BSpl1},
Cor.~\refbpld, that
\begin{equation}
\naam{e: deco LtwospX by repparabs}
L^2(\spX) = \oplus_{Q \in \repparabs} \;\;L^2_Q(\spX),
\end{equation}
with orthogonal $K$-invariant direct summands.

\begin{lemma}
\naam{l: FOU on LtwoQX}
Let $Q \in \repparabs.$ The space $L^2_Q(\spX)$ is $G$-invariant.
Moreover,  $\FOU$ maps $L^2_Q(\spX)$ into $\Ltwo(Q).$
\end{lemma}

\proof We shall first prove the second assertion. Fix $P \in
\repparabs,$ $P \neq Q$ and assume that $\types \subset \dK$ is a
finite subset. Then it follows from \bib{BSpl1},
Cor.~\refbpld{} and
Thm.~\refplr{} (c), that $\cF_P = 0$ on $L^2_Q(\spX\col
\tautypes).$ In view of Lemma \refer{l: Fou Q and FOU} this
implies that $\pr_P\after \FOU (f) = 0$ for every $f \in
L^2_Q(\spX)_\types;$ here $\pr_P$ denotes the orthogonal
projection $\Ltwo \to \Ltwo(P)$. By density of $L^2_Q(\spX)_K$ in
$L^2_Q(\spX)$ and continuity of the map $\pr_P \after \FOU:
L^2(\spX) \to \Ltwo(P),$ see Theorem \refer{t: G equivariance
FOU}, it follows that $\pr_P \after \FOU$ vanishes on
$L^2_Q(\spX),$ for every $P \in \repparabs \setminus \{Q\}.$ The
second assertion now follows by orthogonality of the decomposition
(\refer{e: defi Ltwo as dir sum}).

Since $\FOU$ is an isometry, its adjoint $\FOU^*$ is surjective
from $\Ltwo$ onto $L^2(\spX).$ Moreover, since $\FOU$ is compatible
with the decompositions (\refer{e: deco LtwospX by repparabs})
and (\refer{e: defi Ltwo as dir sum}), it follows by orthogonality
of the mentioned decompositions that
$$
L^2_Q(\spX) = \FOU^* (\Ltwo(Q)).
$$
By $G$-equivariance of $\FOU$
and unitarity of the representations $L$ and $\pi,$ the map $\FOU^*$ is $G$-equivariant.
If follows that $L^2_Q(\spX)$ is $G$-invariant.
\qed

We denote by $\FOU_Q$ the restriction of $\FOU$ to $L^2_Q(\spX),$ viewed
as a map into $\Ltwo(Q).$

\begin{cor}
\naam{c: FOU dir sum over Q}
The map $\FOU$ is the direct sum of the maps $\FOU_Q,$ for $Q \in \repparabs.$
\end{cor}

If $\cH$ is a Hilbert space, we denote by $\End(\cH)$ the space
of continuous linear endomorphisms of $\cH,$
equipped with the operator norm. By $\rmU(\cH)$ we denote the subspace of
unitary endomorphisms. If $P \in \allparabs,$ we define $W(\faPq) = W(\faPq\mid \faPq)$
as in \bib{BSpl1}, \S~\refbplf. Then  by \bib{BSpl1}, Cor.~\refbplg,
\begin{equation}
\naam{e: WfaPq as quotient}
W(\faPq) \simeq W^*_P/W_P.
\end{equation}

\begin{prop}
\naam{p: existence fC}
For each $s \in W(\faQq)$ there exists a measurable map $\fC_{Q,s}:
i\faQqd \to \End(\fH(Q)),$ which is almost everywhere uniquely
determined, such that $\nu \mapsto \|\fC_{Q,s}(\nu)\|$ is essentially
bounded,  and such that for every $f \in \Cci(\spX),$
\begin{equation}
\naam{e: transformation prop FOU}
\FOU_Q f (s\nu) = \fC_{Q,s}(\nu) \,\FOU_Q f (\nu),
\end{equation}
for almost all $\nu \in i \faQqd.$
For almost every $\nu \in i\faQqd$
the map $\fC_{Q,s}(\nu)$ is
unitary. Moreover, for all $s, t \in W(\faQq),$
\begin{equation}
\naam{e: transformation props fC}
\fC_{Q, st}(\nu) = \fC_{Q, s} (t\nu)\after \fC_{Q, t}(\nu).
\end{equation}
In particular, $\fC_{Q,1}(\nu) = I$ and $\fC_{Q, s}(\nu)^{-1} = \fC_{Q,s^{-1}}(s\nu),$
for all $s \in W(\faQq).$
\end{prop}

For $Q$ minimal this result is part of Prop.~18.6 in \bib{BSmc}.
In the present more general setting, we initially reason in a similar way.
For $\Omega \subset i \faQqd$ a measurable subset, we denote
by $\Ltwo_\Omega(Q)$ the closed $G$-invariant subspace of $\Ltwo(Q)$ consisting
of square integrable functions $i\faQqd \to \fH(Q)$
that vanish almost everywhere
outside $\Omega.$ The orthogonal projection
onto this subspace is denoted by $\gf \mapsto \gf_\Omega.$

The uniqueness statement of Proposition \refer{p: existence fC}
follows from the following lemma, which generalizes \bib{BSmc},
Lemma 18.7. We denote by $\fa_{Q\iq}^{*\reg}$ the collection of
elements $H\in \faQqd$ whose parabolic equivalence class relative
to $(\faqd, \gS)$ is open in $\faQqd.$ The set
$\fa_{Q\iq}^{*\reg}$ consists of finitely many connected
components, called chambers. The group $W(\faQq)$ acts freely, but
in general not transitively, on the collection of chambers;
therefore, there exists an open and closed fundamental domain for
$W(\faQq)$ in $\fa_{Q\iq}^{*\reg}.$

\begin{lemma}
\naam{l: surjectivity of fou relative to Omega}
Let $\Omega \subset i \fa_{Q\iq}^{*\reg}$ be an
open and closed fundamental domain for
$W(\faQq).$
Then $f \mapsto (\pr_Q  \FOU f)_\Omega$ maps $\Cci(\spX)$ onto a dense subspace
of $\Ltwo_\Omega(Q),$ and $\Cci(\spX)_\types$ onto a dense subspace of
$\Ltwo_\Omega(Q)_\types,$
for every finite set $\types \subset \dK.$
\end{lemma}

\proof The proof is similar to the proof of Lemma 18.7 of
\bib{BSmc}. Fix a finite subset $\types\subset K;$ then it
suffices  to prove the statement about $\Cci(\spX)_\types,$ by
density of the $K$-finite vectors. Let $T \in
\Ltwo_\Omega(Q)_{\types},$ and suppose that $\hinp{\pr_Q\FOU f}{T}
= 0$ for all $f \in \Cci(\spX)_\types.$ Then it suffices to show
that $T = 0.$ Put $T(\nu) = \sum_{\omega} T(\omega\col \nu)$ with
$T(\omega \col \nu) \in \fH(Q, \omega)_\types.$ Note that this sum
is finite by Lemmas \refer{l: finiteness discserQtypes} and
\refer{l: fH Q omega types non trivial}. Let $\tau = \tautypes,$
then $\psi_{T(\omega\col \nu)} \in \cAtwoQ(\tau).$ We put
$$
\Psi(\nu) := \Psi_\types(T( \nu)) =
\sum_{\omega \in \discserQ}
\psi_{T(\omega\col \nu)} \in \cAtwoQ(\tau).
$$
Note that for $Q$ minimal, the constants $d_\omega$ that occur in \bib{BSmc}
are absent here, see Remark \refer{r: remark on discserQ}.
Let $F \in \Cci(\spX\col \tau),$ and let $f = F(\dotvar)(e),$ then
$f \in \Cci(\spX)_\types$ and $F = \sphisotypes(f),$ see Lemma \refer{l: sphericalization}.
Moreover, as in \bib{BSmc}, proof of Lemma 18.7,
\begin{equation}
\naam{e: hinp Fou Q f Psi}
\hinp{\Fou_Q F}{  \Psi} = \hinp{\hat f(Q)}{T} = \hinp{\pr_Q \FOU f}{T} = 0.
\end{equation}
Let the space $[L^2(i\faQqd) \otimes  \cAtwoQ(\tau)]^{W(\faQq)}$ be defined
as in \bib{BSpl1}, text before Cor.~\refbplh. It follows from the definition of this space
that the restriction map $\gf \mapsto \gf|_\Omega$ is a bijection from it
onto $L^2(\Omega) \otimes \cAtwoQ(\tau).$

The image of $\Cci(\spX\col \tau)$ under $\Fou_Q$ is dense in the
space $[L^2(i\faQqd) \otimes  \cAtwoQ(\tau)]^{W(\faQq)},$ by
\bib{BSpl1}, Thm.~\refplr{}~(c). Combining this  with (\refer{e: hinp
Fou Q f Psi}) we see that $\Psi$ is perpendicular to the mentioned
space. Since $\Psi = 0$ outside $\Omega,$ we infer that
$\Psi|_\Omega$ is perpendicular to $L^2(\Omega) \otimes
\cAtwoQ(\tau).$ We conclude that $\Psi,$ hence $T,$ is zero. \qed
\medno {\bf Proof of Proposition \refer{p: existence fC}:\ } We
fix a finite subset $\types \subset \dK$ and put $\tau =
\tau_\types.$ We will first prove that there exists a measurable
map $\fC_{Q,s, \types}: i\faQqd \to \End(\fH(Q)_\types),$ such
that (\refer{e: transformation prop FOU}) is valid with $\fC_{Q,s,
\types}$ in place of $\fC_{Q,s},$ for every $f \in
\Cci(\spX)_\types.$ We define
\begin{equation}
\naam{e: fC and C}
\fC_{Q,s, \types}(\nu) = \Psi_\types^{-1} \after \nC_{Q\mid Q}(s \col \nu) \after  \Psi_\types,
\end{equation}
where $\Psi_\types$ is the isometry from $\fH(Q)_\types$ onto
$\cAtwoQ(\tau)$ defined in the text preceding Lemma \refer{l: Psi types isometry}
and where the $C$-function is defined as in \bib{BSpl1}, Def.~\refbpli,
with $\tau = \tautypes.$
We note that the $\End(\fH(Q)_\types)$-valued function
$\fC_{Q,s, \types}$  is analytic on $i\faQqd,$ by \bib{BSpl1}, Cor.~\refbplj.
It follows from Lemma \refer{l: Psi types isometry} combined with the
Maass--Selberg relations for the $C$-function, see \bib{BSpl1}, Thm.~\refbplk,
that $\fC_{Q,s, \types}(\nu)$
maps $\fH(Q)_\types$ unitarily into itself, for $\nu \in i\faQqd.$
{}From (\refer{e: fC and C}) and \bib{BSpl1},
Lemma \refbpll{}
with $P =R = Q,$ we deduce
that (\refer{e: transformation props fC})
is valid with everywhere the index $\types$ added.

In view of Lemma \refer{l: surjectivity of fou relative to Omega},
the function $\fC_{Q,s, \types}$ is uniquely determined by
the requirements in the beginning of this proof.
If $\typesp\subset \dK$ is a second finite subset with $\types \subset \typesp,$
let $\rmi_{\typesp, \types}$ denote the inclusion map $\fH(Q)_\types \to \fH(Q)_\typesp,$
and let $\pr_{\types,\typesp}$ denote the orthogonal projection
$\fH(Q)_\typesp \to \fH(Q)_\types.$
Then it follows from the uniqueness that
$$
\pr_{\types,\typesp} \after \fC_{Q,s, \typesp}(\nu)\after \rmi_{\typesp, \types} =
\fC_{Q,s, \types}(\nu),
$$
for every $\nu \in i\faQqd.$ By unitarity of the
endomorphisms $\fC_{Q,s, \types}(\nu)$ and $\fC_{Q,s, \typesp}(\nu)$
this implies that $\fC_{Q,s, \typesp}(\nu)$ leaves the subspace
$\fH(Q)_\types$ of $\fH(Q)_\typesp$ invariant, and equals $\fC_{Q,s,\types}(\nu)$
on it. Thus, we may define the endomorphism $\fC_{Q,s}(\nu)$
of $\fh(Q)$ by requiring it to be equal to $\fC_{Q,s, \types}(\nu)$ on $\fh(Q)_\types,$
for every finite subset $\types \subset \dK.$ The endomorphism defined
depends measurably on $\nu,$
has essentially bounded norm and satisfies
(\refer{e: transformation prop FOU}).
We asserted already that it is uniquely determined by these properties, in view of Lemma
\refer{l: surjectivity of fou relative to Omega}. The remaining asserted properties
of $\fC_{Q,s}(\nu)$ follow from the discussion above.
\qed

\begin{lemma}
\naam{l: FOU Q isometry onto LtwoQOmega}
Let $Q \in \repparabs$  and let $\Omega$ be
an open and closed fundamental domain for the $W(\faQq)$-action in $i \fa_{Q\iq}^{*\reg}.$
Then the map $f \mapsto |W(\faQq)|^{1/2}(\FOU_Q f)_\Omega$ defines an isometric isomorphism
from $L^2_Q(\spX)$ onto $\fL^2(Q)_\Omega,$ intertwining the restriction of $L$ to $L^2_Q(\spX)$
with the direct integral representation
\begin{equation}
\naam{e: direct integral Omega Q}
\widehat\oplus_{\omega \in \discserQ}\;\;
\int_\Omega^\oplus  1_{\bar V(Q,\omega)} \otimes
\pi_{Q, \xi, -\nu} \;\; [W: W^*_Q]\, d\mu_Q(\nu).
\end{equation}
of $G$ in $\Ltwo(Q)_\Omega.$
\end{lemma}

\proof
In view of Theorem \refer{t: G equivariance FOU} and
Lemma \refer{l: FOU on LtwoQX}, the map
$\FOU_Q$ is an isometry from $L_Q^2(\spX)$ into $\fL^2(Q),$ intertwining
the restriction of $L$ to $L_Q^2(\spX)$ with $\pi_Q: = \pi|_{\fL^2(Q)}.$
The map $\gf \mapsto \gf_\Omega$
from $\Ltwo(Q)$ to $\fL^2(Q)_\Omega$ intertwines $\pi_Q$ with the direct integral
(\refer{e: direct integral Omega Q}). Thus, it remains to show that the map
$T: f \mapsto |W(\faQq)|^{1/2} (\FOU_Q f)_\Omega$ from $L^2_Q(\spX)$ to $\fL^2(Q)_\Omega$
is isometric and onto.

To establish the first property, we note that, for $f \in L^2_Q(\spX)$ and
$s \in W(\faQq),$
$$
\|\FOU_Q f(s\nu)\| = \| \fC_{Q,s}(\nu) \FOU_Q f(\nu)\| = \|\FOU_Q f(\nu)\|,
$$
for almost every $\nu \in i \faQqd,$ by Proposition \refer{p: existence fC}.
Hence,
\begin{eqnarray*}
\|f\|_2^2 & = & \|\FOU f\|^2_{\fL^2} = \|\FOU_Q f\|^2_{\fL^2(Q)} \\
&=&
\int_{i\faQqd} \|\FOU_Q f (\nu)\|^2 \;[W:W^*_Q]\, d\mu_Q(\nu)\\
&=&
\sum_{s \in W(\faQq)} \int_{\Omega} \|\FOU_Q f (s\nu)\|^2 \;[W:W^*_Q]\, d\mu_Q(\nu)\\
&=&
|W(\faQq)| \, \int_{\Omega} \|\FOU_Q f (\nu)\|^2 \;[W:W^*_Q] \,d\mu_Q(\nu).
\end{eqnarray*}
It follows that $T$ is an isometry.
On the other hand, $T$ has dense image
in view of Corollary \refer{c: FOU dir sum over Q}
and Lemma \refer{l: surjectivity of fou relative to Omega}.
We conclude that $T$ is surjective.
\qed

We shall now investigate irreducibility and equivalence of the occurring representations
$\pi_{Q, \omega, \nu}.$

\begin{lemma}
\naam{l: comparison infinitesimal characters}
Let $\pi \in \discserG.$ Then $\pi$ has a real infinitesimal $Z(\fg)$-character
in the following sense. Let $\fj$ be a Cartan subalgebra of $\fg,$ $W(\fj)$
the Weyl group of the root system $\gS(\fj)$ of $\fj_\iC$ in $\fg_\iC.$
Let  $[\gL] \in \fj_\iC^*/W(\fj)$  be the infinitesimal character of $\pi.$
Then $\inp{\gL}{\ga}$ is real for every $\ga \in \gS(\fj).$
\end{lemma}

\proof
Let $\fb$ a $\Cartan$-stable Cartan subspace of $\fq,$ let $\gS(\fb)$ be
the root system of $\fb$ in $\fg_\iC,$ $W(\fb)$ the associated
Weyl group and $I(\fb)$ the algebra of $W(\fb)$-invariants in $S(\fb),$
the symmetric algebra of $\fb_\iC.$ Let $\gg_\fb: \DX \to I(\fb)$
be the associated Harish-Chandra isomorphism. Let $L^2(\spX)_\pi$
be defined as in (\refer{e: canonical embedding discrete series})
with $\pi$ in place of $\xi.$ We may
fix a non-zero simultanous eigenfunction $f$ for $\DX$ in $[L^2(\spX)_\pi]^\infty.$
Let $\gl \in \fbdc$ be such that $Df = \gg_\fb(D\col \gl) f$
for $D \in \DX.$ Then in particular, for each element $Z$ of
$\fZ,$ the center of $U(\fg),$
\begin{equation}
\naam{e: R Z on f}
R_Z f  = \gg_\fb(Z\col \gl)f.
\end{equation}
On the other hand, $R_Z f = L_{\check Z} f$ may be expressed in
terms of the infinitesimal character of $\pi$ as follows.
Let $\fj$ be a $\Cartan$-stable Cartan subalgebra of $\fg,$ containing
$\fb,$  let $\gS(\fj)$ and $W(\fj)$ be as in the lemma,
and let $I(\fj)$ denote the algebra
of $W(\fj)$-invariants in $S(\fj).$ We denote the canonical isomorphism
${\fZ} \to I(\fj)$ by $\gg_\fj.$ Let $\gL \in \fj_\iC$
be as in the lemma, then
\begin{equation}
\naam{e: R Z on f two}
R_Z f = L_{\check Z} f = \gg_\fj(\check Z\col \gL) f = \gg_\fj(Z \col - \gL) f.
\end{equation}
{}From (\refer{e: R Z on f}) and (\refer{e: R Z on f two}) we obtain that
\begin{equation}
\naam{e: relation gg fb and gg fj}
\gg_\fj(Z\col - \gL) = \gg_\fb( Z \col\gl),\qquad (Z \in \fZ).
\end{equation}
Let $\fl$ be the centralizer of $\fb$ in $\fg,$ let $\gS(\fl_\iC, \fj)$
be the root system of $\fj$ in $\fl_\iC,$
$\gS^+(\fl_\iC, \fj)$
a choice of positive roots and $\gd_\fl \in \fj_\iC^*$
the associated half sum of the positive roots.
By a standard computation in the universal enveloping algebra,
involving the definitions of $\gg_\fj$ and $\gg_\fb,$ it follows that
$\gg_\fb(Z\col \gl) = \gg_{\fj}(Z\col \gl - \gd_\fl),$ for all $Z \in \fZ.$
Combining this with (\refer{e: relation gg fb and gg fj})
we obtain that $- \gL$ and $\gl - \gd_\fl$ are $W(\fj)$-conjugate.

Now $\inp{\gl}{\ga}$ is real, for each $\ga \in \gS(\fb),$ by
\bib{OMds}, see also
\bib{BSpl1}, Thm.~\refplm. It readily follows that $\gL$ is real.
\qed

The following result is due to F.~Bruhat \bib{Bru} for minimal parabolic
subgroups and to Harish-Chandra in general. A proof is essentially given
in \bib{KV}.

\begin{thm}
\naam{t: HC irreducibility}
For $j=1,2,$ let  $P_j = M_jA_jN_j$ be a parabolic subgroup of $G,$ with
the indicated Langlands decomposition.
Moreover, let  $\xi_j$ be an irreducible unitary representation of $M_j$ with real
infinitesimal character and let
$\nu_j \in i \fa_j^*$ be regular with respect to the roots of  $\fa_j$
in $P_j.$ Let $\pi_j$ denote the unitarily induced representation
$\Ind_{P_j}^G (\xi_j \otimes \nu_j \otimes 1).$
\begin{enumerate}
\itema
The representation $\pi_j$ is
irreducible, for $j =1,2.$
\itemb
The representations $\pi_1$ and $\pi_2$ are
equivalent if and only if the data $(\fa_j, \nu_j, [\xi_j]),$ for $j =1,2,$
are conjugate
under $K.$
The latter condition means that there exists  $k \in K$
such that $\Ad(k) \fa_1 = \fa_2,$ $\nu_1 \after \Ad(k)^{-1} = \nu_2$
and $\xi_1^k \sim \xi_2,$ where $\xi_1^k:= \xi_1(k^{-1}(\dotvar) k)|_{M_2}.$
\end{enumerate}
\end{thm}

\proof
Taking into account the actions of the centers of $M_j$ and $G,$
one readily checks that it suffices to prove this result for
$G$ connected semisimple and with finite center. Thus, let us assume
this to be the case.

Assertion (a) follows from \bib{KV}, Thm.~4.11. Thus,
it remains to prove assertion (b).
We first establish the `if' part. If in addition to the hypothesis,
$P_2 = k P_1 k^{-1},$ then the equivalence of $\pi_1, \pi_2$ is
a trivial consequence of conjugating all induction data.
Thus, by applying conjugation we may reduce to the case that
$\fa_1 = \fa_2,$ $\nu_1 = \nu_2$ and $\xi_1 \sim \xi_2.$ Then $P_1$ and $P_2$
have the same split component. It now follows from \bib{KS}, Prop.~8.5 (v),
that there exists
a unitary intertwining operator from $\pi_1$ onto $\pi_2.$
Hence $\pi_1 \sim \pi_2.$

We shall now prove the `only' if part. Assume that $\pi_1 \sim \pi_2.$
By conjugating all induction data of $\pi_1$
with an element of $K,$ we see that we may restrict ourselves to the situation
that $P_1$ and $P_2$ contain a fixed minimal parabolic subgroup $P_0$ of $G,$ with split
component $A_0.$ In particular, $\fa_j \subset \fa_0$ for $j=1,2.$ It now follows from
\bib{KV}, p.~94, text under the heading `Equivalence', that there exists a $k \in N_K(\fa_0)$
such that $\Ad(k) \fa_1 = \fa_2$ and $\nu_1 \after \Ad(k)^{-1} = \nu_2.$ Conjugating all
data of $\pi_1$ with $k$ we see that we may as well assume that $\fa_1 = \fa_2$ and $\nu_1 = \nu_2.$
Moreover, applying \bib{KS}, Prop.~8.5 (v), as in the first part of this proof we
see that in addition we may assume that $P_1 = P_2.$ We now claim that
$\xi_1 \sim \xi_2.$ This assertion is essentially proved in \bib{KV}, proof of Thm.~4.11, but not explicitly
stated as a result. We shall indicate how to modify the mentioned proof.
We use the notation of \bib{KV}. In particular, $\xi_j = {}^\circ \gs_j.$
We follow the proof of \bib{KV}, Thm.~4.11,
after the heading `equivalence', but with
$P_1 = P_2 = P$ and $\nu_1 = \nu_2 = \nu.$ {}From $i(\pi_1, \pi_2^*) > 0$
it follows, by application of \bib{KV}, Thm.~4.10, that
$$
(M_1^{(0)} \otimes (E_1^0 \bar\otimes {E_2^0}')' \otimes \C_1')^{(P \otimes P)}
$$
has positive dimension. Now $M_1^{(0)} $ equals $\Ci(P),$ equipped with the left times right
action of $P$ (see \bib{KV}, Eqn.~(2.6)). Hence the above space is naturally isomorphic
with the space of ${\rm diag}(P \times P)$-invariants in
$(E_1^0 \bar\otimes {E_2^0}')'$
which in turn is naturally isomorphic with $\Hom_P(E_2^0, E_1^0) =
\Hom_M({}^\circ \gs_2, {}^\circ \gs_1).$ It follows that the latter space is non-trivial,
hence ${}^\circ \gs_1 \sim {}^\circ \gs_2,$
since the representations involved are irreducible.
\qed

\begin{prop}
\naam{p: distinct reps}
For $j=1,2,$ let  $Q_j \in \repparabs,$  $\omega_j \in M^{\wedge}_{Q_j, H},$
$\nu_j \in i\fa_{Q_j\,\iq}^{*\reg}.$ Then the representations
$\pi_j= \pi_{Q_j,\omega_j, \nu_j}$ are
irreducible. Moreover, they are equivalent if and only if $Q_1 = Q_2$
and there exists  $s \in W(\fa_{Q_1\iq})$ such that $\nu_2 = s\nu_1$
and $\omega_2 = s\omega_1.$
\end{prop}

\proof
Put $\xi_j = \xi_{\omega_j}$ and $Q =Q_1.$
There exists $v \in \QcW,$ such that $\omega_1$ belongs to the
discrete series of $M_Q/M_Q\cap v H v^{-1}.$ It follows from Lemma
\refer{l: comparison infinitesimal characters}
that $\omega_1$ has a real infinitesimal character
for the center of $U(\fm_Q).$ A similar statement holds for $Q_2, \omega_2.$

If $\ga$ is a root of $\fa_{Q}$ in $Q,$
then its restriction $\ga_\iq$ to $\faQq$ belongs
to $\gSr(Q).$
Moreover, $\inp{\ga}{\nu_1} = \inp{\ga|_{\faQq}}{\nu_1} \neq 0.$ Thus, it follows
that $\nu_1$ is regular with respect to the $\fa_Q$-roots in $Q.$
A similar statement holds for $\nu_2.$

Thus, Theorem \refer{t: HC irreducibility} is applicable and we conclude that
$\pi_1$ and $\pi_2$ are irreducible.

Assume that $\pi_1 \sim \pi_2.$ Then by Theorem \refer{t: HC irreducibility} (b)
we conclude that there exists $k \in K$ such that $\Ad(k) \fa_1 = \fa_2,$
$\nu_1 \after \Ad(k)^{-1} =  \nu_2$ and $\xi_1^k  \sim \xi_2.$
Since $\gs \nu_j = - \nu_j,$ for $j =1,2,$ it follows by application
of $\gs$ that also $\nu_1 \after \Ad(\gs k)^{-1} = \nu_2.$ We infer that
$\Ad((\gs k)^{-1}k)^{-1*}$ centralizes $\nu_1,$ hence
belongs to the centralizer $M_{1Q}$ of $\faQq,$ by regularity of $\nu_1.$
The mentioned element therefore centralizes $\fa_Q,$ from which we see
that $\Ad (k) = \Ad(\gs k)$ on $\fa_Q.$ This implies that
$\gs \after \Ad(k) = \Ad( k)\after \gs$ on $\fa_Q,$
hence $\Ad(k)$ maps $\fa_{Q\iq}$ onto $\fa_{Q_2 \iq}.$ We conclude
that $Q_1$ and $Q_2$ are associated, hence equal. Put $Q = Q_1 = Q_2.$

It follows from the above that
$\Ad(k)$ normalizes $\faQq.$ Hence, $s:=\Ad(k)|_{\faQq}$ belongs to $W(\faQq),$
see \bib{BSpl1}, \S~\refbplf.
Finally, it follows that $s\nu_1 = \nu_2$ and
$s[\xi_1] = [\xi_1^k] = [\xi_2].$
\qed

\begin{thm}
\naam{t: plancherel with dir int}
Let, for each $Q \in \repparabs$ an open and closed fundamental domain $\Omega_Q$ for
the action of $W(\faQq)$ on $i\fa_{Q\iq}^{*\reg}$ be given.
The Fourier transform $\FOU$ induces the following
Plancherel decomposition of the regular representation $L$ of $G$ in $L^2(\spX):$
\begin{equation}
\naam{e: deco L}
L \simeq \oplus_{Q \in \repparabs}\widehat \oplus_{\omega \in \discserQ}
\int_{\Omega_Q}^\oplus 1_{\bar V( Q, \omega)} \otimes \pi_{Q, \omega, \nu}\;
|W||W_Q|^{-1} \,d\mu_Q(\nu),
\end{equation}
and
\begin{equation}
\naam{e: Plancherel isometry}
\|f\|^2 = \sum_{Q \in \repparabs} \sum_{\omega \in \discserQ} \int_{\Omega_Q}^\oplus
\|\FOU_Q f( \omega \col - \nu)\|^2 \;  |W||W_Q|^{-1} \,d\mu_Q(\nu),
\end{equation}
for every $f \in L^2(\spX).$ In particular, for each $Q \in \repparabs$ and every
$\omega \in \discserQ,$ the induced representation $\pi_{Q, \omega, \nu}$
occurs with multiplicity $m_{Q,\omega} = \dim \bar V(Q, \omega),$ for almost
every $\nu \in \Omega_Q.$
\end{thm}
\proof
The fact that $\FOU$ induces the isometric isomorphism
of $L$ with the direct integral
as expressed by (\refer{e: deco L}) and (\refer{e: Plancherel isometry})
follows from Lemma \refer{l: FOU Q isometry onto LtwoQOmega}
applied to $-\Omega_Q$ in place of $\Omega,$ combined with Corollary
\refer{c: FOU dir sum over Q} and (\refer{e: WfaPq as quotient}).
The occurring representations  $\pi_{Q,\omega, \nu}$ are all irreducible, by
Proposition \refer{p: distinct reps}. It remains to exclude double occurrences.
For $j = 1,2,$ let  $Q_j \in \repparabs,$  $\omega_j \in M^{\wedge}_{Q_j, H},$
$\nu_j \in i\fa_{Q\iq}^{*\reg},$ and assume that $\pi_{Q_1, \omega_1, \nu_1} \sim
\pi_{Q_2, \omega_2, \nu_2}.$ Then it follows from Proposition
\refer{p: distinct reps} that
$Q_1 = Q_2.$
Moreover, there exists $s \in W(\fa_{Q_1\iq})$ such that
$(s\nu_1, s\omega_1) = (\nu_2, \omega_2).$
Since $\Omega_{Q_1} = \Omega_{Q_2}$ is a fundamental
domain for the $W(\fa_{Q_1\iq})$-action, it follows that $s =1.$
\qed

We finish this section with a description of the image
of the isometry $\FOU: L^2(\spX) \to \Ltwo.$

\begin{lemma}
Let $Q \in \repparabs$ and $s \in W(\faQq).$ Then for almost every $\nu \in i\faQqd,$
the unitary endomorphism $\fC_{Q,s}(\nu)$ of $\fH(Q)$ maps the subspace
$\fH(Q,\omega)$ onto $\fH(Q, s\omega),$ for $\omega \in \discserQ,$
intertwining the representations $\pi_{Q, \omega, \nu}$ and $\pi_{Q, s\omega, s\nu}.$
\end{lemma}
\proof
Let $s \in W.$
We will first show that, for almost every $\nu \in i\faQqd,$
the unitary endomorphism $\fC_{Q,s}(\nu)$
of $\fH(Q)$ intertwines the direct sum $\pi_{Q,\nu}$
of the representations
$1 \otimes \pi_{Q,\omega,-\nu},$ for $\omega \in \discserQ,$
with the similar direct sum $\pi_{Q,s\nu}$ of the representations
$1 \otimes \pi_{Q,\omega, -s\nu}.$ Let $\Omega$ be an open and closed
fundamental domain for the $W(\faQq)$-action on $i \fa_{Q\iq}^{*\reg}.$
Then the map $\FOU_{Q\Omega}: f \mapsto \FOU_Q(f)_\Omega$  is an
equivariant isometry
from $L^2_Q(\spX)$ onto $ \Ltwo(Q)_\Omega,$
by Lemma \refer{l: FOU Q isometry onto LtwoQOmega}.
Similarly, the map $\FOU_{Q\, s\Omega}$ is an  intertwining
isomorphism from $L^2_Q(\spX)$ onto $\Ltwo(Q)_{s\Omega}.$ Moreover,
by (\refer{e: transformation prop FOU}), for every $f \in L^2_Q(\spX)$ we have
$$
s^*(\FOU_{Q\, s\Omega} f)(\nu) = \FOU_{Q\, s\Omega} f(s\nu)
= \fC_{Q,s}(\nu) \Fou_{Q\Omega}f(\nu),
$$
for almost every $\nu \in \Omega.$ It follows that the map
$
\gf \mapsto s^{*-1}[\fC_{Q,s}(\dotvar) \gf]
$
is an equivariant isometry from
$\Ltwo(Q)_\Omega$ onto $\Ltwo(Q)_{s\Omega}.$ Let $x \in G.$
Then
$$
s^* \, \pi_{Q}(x) s^{*-1}[\fC_{Q,s}(\dotvar)\gf] = \fC_{Q,s}(\dotvar) \pi_Q(x)\gf,
$$
for every $\gf \in \fL^2(Q).$ It follows that
\begin{equation}
\naam{e: intertwining prop fC}
\pi_{Q,s \nu}(x) \fC_{Q,s}(\nu) = \fC_{Q,s}(\nu) \pi_{Q,\nu}(x)
\end{equation}
for almost every $\nu \in \Omega.$ Since $\Omega$ was arbitrary,
(\refer{e: intertwining prop fC}) holds for almost every $\nu \in i\faQqd.$
Select a countable dense subset $G_0$ of $G.$
Then there exists a subset $\cA \subset i \fa_{Q\iq}^{*\reg}$
with complement of measure zero in $i\faQqd,$
such that $\fC_{Q,s}(\dotvar)$ is represented by a function on $\cA$
with values in $\rmU(\fH(Q)),$
satisfying (\refer{e: intertwining prop fC}) for all $x \in G_0$ and $\nu \in \cA.$
By continuity of $\fC_{Q,s}(\nu)$ and of the representations
$\pi_{Q,\nu}$ and $\pi_{Q,s\nu},$
it follows that (\refer{e: intertwining prop fC})
holds for all $\nu \in \cA$
and all $x \in G.$ In view of Theorem \refer{t: HC irreducibility},
the representation $1 \otimes \pi_{Q, \tilde \omega , - s\nu},$ for
$\tilde \omega \in \discserQ$ and $s \in \cA,$ is not disjoint from
$1 \otimes \pi_{Q, \omega, - \nu},$ if and only if
$\tilde \omega  = s \omega.$ All assertions now follow
for all  $\nu \in \cA.$
\qed

It follows from the above result that for each $s \in W(\faQq),$ we may define
a unitary endomorphism $\Gamma_Q(s)$ of $\Ltwo(Q)$ by
$$
[\Gamma_Q(s) \gf ](\nu) =  \fC_{Q, s}( s^{-1} \nu) \gf(s^{-1}\nu),
$$
for $\gf \in \Ltwo(Q)$ and almost every $\nu \in i\faQqd.$
Moreover, the map $\Gamma_Q(s)$ intertwines $\pi_Q$ with itself. It follows from
Proposition \refer{p: existence fC}
that $s \mapsto \Gamma_Q(s)$ defines a unitary representation of $W(\faQq)$ in $\Ltwo(Q),$
commuting with the action of $G.$ Accordingly, the associated space $\Ltwo(Q)^{W(\faQq)}$
of invariants is a closed $G$-invariant subspace of $\Ltwo(Q).$

\begin{thm}\vbox{\hspace{1mm}}
\naam{t: image FOU}
\begin{enumerate}
\itema
For each $Q \in \repparabs,$ the image of $\FOU_Q$ equals $\Ltwo(Q)^{W(\faQq)}.$
\itemb
The image of the Fourier transform $\FOU$ is given by the following
orthogonal direct sum
$$
\image(\FOU) = \oplus_{Q\in \repparabs}\;\; \Ltwo(Q)^{W(\faQq)}.
$$
\end{enumerate}
\end{thm}

\proof
{}From Proposition \refer{p: existence fC} it follows that $\FOU_Q$ maps into
$\Ltwo(Q)^{W(\faQq)}.$
Thus, for (a) it remains to prove the surjectivity. Let $\Omega$ be
an open and closed fundamental domain in $i \fa_{Q\iq}^{*\reg}$
for the action of $W(\faQq).$ Then the map $\gf \mapsto \gf|_{\Omega}$
is a bijection from $\Ltwo(Q)^{W(\faQq)}$ onto $\Ltwo(Q)_\Omega.$
The surjectivity now follows by application of Lemma
\refer{l: FOU Q isometry onto LtwoQOmega}.

Finally, assertion (b) follows from (a) combined with Corollary
\refer{c: FOU dir sum over Q}.
\qed

\eqsection{$H$-fixed generalized vectors, final remarks}
\naam{s: H-fixed generalized vectors, final remarks}
In this section we compare our results with those
obtained by P.\ Delorme in \bib{Dpl}. This comparison relies heavily on the automatic continuity
theorem, due to W.\ Casselman and N.R.\ Wallach, see \bib{Cas}, Cor.~10.5 and
\bib{Wal2},
Thm. 11.6.7. We shall therefore
first recall this result. The group decomposes as  $G \simeq {}^\circ G \times  \exp \fc_\ip,$
where, as usual,  ${}^\circ G $ denotes the intersection of all subgroups $\ker |\chi|,$
with $\chi$ a continuous homomorphism $G \to \R^*,$ and where
$\fc_\ip = {\rm center}(\fg) \cap \fp.$
Accordingly, we define
the function $\|\dotvar\|: G \to \;]\, 0, \infty \,[$ by
\begin{equation}
\naam{e: defi norm}
\|x\exp H \| = \|\Ad(x)\|_{\rm op} \,e^{|H|},
\end{equation}
for $x \in {}^\circ G$ and $H \in \fc_\ip;$
here $\|\dotvar\|_{\rm op}$ denotes the operator norm on $\End(\fg).$
Let $\fa_\ip$ be a maximal abelian subspace of $\fp$ containing $\faq,$
and let $\gS(\fa_\fp)$
be the root system of $\fa_\fp$ in $\fg.$ Then one readily checks that
\begin{equation}
\naam{e: norm and KAK}
\|k_1 a  k_2 \| =  \max_{\ga \in \gS(\fa_\fp)} a^\ga,
\end{equation}
for $k_1,k_2 \in K$ and $a \in A_\fp \cap {}^\circ G.$ In particular, it follows
that $\|\dotvar\| \geq 1$ on $G.$ Note that it follows from (\refer{e: defi norm}) and
(\refer{e: norm and KAK}) that
\begin{equation}
\naam{e: equality of norm inverse}
\|x\| = \|x^{-1}\|\qquad (x \in G).
\end{equation}
We recall from \bib{Wal2}, 11.5.1,   that a smooth representation $\pi$
of $G$ in a Fr\'echet space $V$ is said to be of moderate growth if for
each continuous seminorm $s$ on $V$ there exists a continuous seminorm $p_s$ on $V$
and a constant $d_s \in \R$ such that
$$
s(\pi(x)v) \leq \|x\|^{d_s} p_s(v),
$$
for all $v \in V$ and $x \in G.$

\begin{thm}{\rm (The automatic continuity theorem)\ }
\naam{t: automatic continuity}
Let $(\pi_j, V_j),$ for $j=1,2,$ be smooth Fr\'echet representations of $G$
of moderate growth,
such that the associated $(\fg, K)$-modules $(V_j)_K$ are finitely generated.
Then every $(\fg, K)$-equivariant linear map $(V_1)_K \to (V_2)_K$ extends
(uniquely) to a continuous linear $G$-equivariant map $V_1 \to V_2.$ Moreover,
the image of the extension is a closed topological summand of $V_2.$
\end{thm}

\begin{rem}
A proof of this theorem, due to W.\ Casselman and N.R.\ Wallach, is given
in \bib{Cas}, Cor.~10.5 and in \bib{Wal2}, Thm.~11.6.7,
but for a somewhat different class of real reductive groups.
In \bib{D1n}, \S 1, it
is shown that the result is valid for groups of Harish-Chandra's class.
\end{rem}

For any function $f: G \to \C$ and
any non-negative real number $r \geq 0$ we define
\begin{equation}
\naam{e: defi r norm}
\|f\|_r:= \sup_{x \in G}\; \|x\|^{-r}\,|f(x)|.
\end{equation}
Moreover, we define $C_r(G)$ to be the space of continuous functions
$f: G \to \C$ with $\|f\|_r < \infty.$ Then $C_r(G),$ equipped with the norm
$\|\dotvar\|_r,$ is a Banach space.

\begin{lemma}
\naam{l: continuity L g on CrG}
For every $g \in G,$ both the left regular action $L_g$ and the right regular action $R_g$
leave the space $C_r(G)$ invariant; their restrictions to the mentioned
Banach space have operator norm at most $\|g\|^r.$
\end{lemma}

\proof
For any function $f: G \to \C,$ we define $f^\vee: G \to \C$ by $f^\vee(x) = f(x^{-1}).$
It follows from (\refer{e: equality of norm inverse}) and (\refer{e: defi r norm})
that $f \mapsto f^\vee$
is an isometry from the Banach space $C_r(G)$ onto itself, intertwining
$R_g$ with $L_g.$ Therefore, it suffices to prove the assertions for the left action.

In view of (\refer{e: defi norm}) it follows from the multiplicative
property of the operator norm that
$\|g x\| \geq \|g^{-1}\|^{-1} \| x\|$ for all $x, g \in G.$
Applying this inequality to the definition
of $\|L_g f\|_r,$ for $f \in C_r(G)$ and $g \in G,$ we see that $L_g$ acts on
$C_r(G)$ with operator norm at most $\|g^{-1}\|^r.$ The lemma now follows by application
of (\refer{e: equality of norm inverse}).
\qed

We note that the left regular representation $L$
of $G$ in $C_r(G)$ is not continuous if $G$ is not compact.
In fact, in that case there exists a function $f \in C_r(G)$ such
that $L_g f$ has no limit for $g \to e.$ However, $L$ does induce
a continuous representation in a subspace that we shall now introduce.

We define $C_r^\infty(G)$ to be the space
of smooth functions $f:G\to \C$ with $L_u f \in C_r(G)$ for all $u \in U(\fg).$
If $F \subset U(\fg)$ is a finite subset, we define the seminorm $\nu_{F,r}$ on
$\CirG$ by
$$
\nu_{F,r}(f): = \max_{u \in F}\; \|L_u f\|_r.
$$
We equip $\CirG$ with
the locally convex topology induced by the
collection of seminorms $\nu_{F,R},$ for $F \subset U(\fg)$ finite.
It is readily seen that the space $\CirG,$ thus topologized, is a Fr\'echet space.

\begin{rem}
The space $C_r^\infty(G)$ has been introduced  in \bib{Cas}, p.\ 424, where it was denoted by
$A_{{\rm umg},r}(G).$ In the mentioned paper it is asserted that
this space is a continuous $G$-module of moderate growth
for the left action. The following result expresses that in fact this $G$-module is smooth.
\end{rem}

\begin{prop}
\naam{p: CirG smooth moderate growth}
Let $r \geq 0.$ The space $\CirG$ is left $G$-invariant. Moreover, the
left regular representation $L$ of $G$ in $C_r^\infty(G)$ is a smooth
Fr\'echet representation of moderate growth.
\end{prop}

\proof
In \bib{BD}, \S~1, the space $\cA_N(G/H),$ for $N \in \N,$ is defined as the analogue
of $C_r^\infty(G),$ with respect to the norm function
$g \mapsto \|g\gs(g)^{-1}\|$ on $G/H,$ in place of the norm function $\|\cdot\|$ on $G.$
Proposition \refer{p: CirG smooth moderate growth} is the analogue
of \bib{BD}, Lemma 1.
The proof of the mentioned lemma may be transferred to
the present situation with obvious modifications.
\qed

\begin{cor}
\naam{c: CirG moderate growth}
Let $r\geq 0.$
Every closed $G$-submodule of $C_r^\infty(G)$ is a smooth Fr\'echet module of
moderate growth.
\end{cor}

\proof
Immediate. See also \bib{Wal2}, Lemma 11.5.2.
\qed

\begin{prop}
\naam{p: automatic continuity and CirG}
Let $(\pi, V)$ be a smooth Fr\'echet representation of $G$ of moderate growth,
such that $V_K$ is finitely generated. Let $r \geq 0$ and let
$T: V_K \to C^\infty_r(G)$ be a $(\fg, K)$-equivariant linear map.
The map $T$ has a unique extension to a continuous linear $G$-equivariant map
$ V \to C^\infty_r(G).$ The image of this extension is closed.
\end{prop}

\proof
Let $W$ be the closure of the image of $T$ in $C^\infty_r(G).$
Then $W$ is a closed $G$-submodule of $C_r^\infty(G),$
hence a smooth Fr\'echet module of moderate growth. Moreover,
$W_K = T(V_K)$ is finitely generated. By Theorem \refer{t: automatic continuity},
$T$ has a unique extension
to a continuous linear $G$-equivariant
map $\tilde T: V \to W.$ The image of this extension is closed and contains a dense subspace
of $W,$ hence equals $W.$
\qed

\begin{lemma}
\naam{l: continuity right regular action}
Let $r \geq 0.$ Then the space $C^\infty_r(G)$ is right $G$-invariant.
Moreover, if $y \in G,$ then the right regular action $R_y$ restricts
to a continuous linear operator of $\CirG.$
\end{lemma}

\proof
It follows from Lemma \refer{l: continuity L g on CrG}
that $R_y$ is a continuous
linear endomorphism of $C_r(G)$ with operator norm at most $\|y\|^r.$
Since the action of $R_y$ on $C^\infty(G)$ commutes with that of $L_u,$ for
every $u \in U(\fg),$ it readily follows that $R_y$ leaves the space $\CirG$ invariant
and restricts to a continuous linear endomorphism of it.
\qed

Let $r \geq 0.$ We define
$$
C_r^\infty(\spX):= C_r^\infty(G) \cap C(G/H),
$$
the space of right $H$-invariant functions in $C_r^\infty(G).$

\begin{lemma}
\naam{l: CirX is closed}
Let $r \in \R.$ The space $C_r^\infty(\spX)$
is a closed $G$-submodule of $C_r^\infty(G).$ In particular, it is a smooth Fr\'echet
$G$-module of moderate growth.
\end{lemma}

\proof
For every $h \in H,$ the map $R_h - I$ restricts to a continuous linear
operator of $\CirG,$ by Lemma \refer{l: continuity right regular action}.
The space $C_r^\infty(\spX)$ equals the intersection in $C_r^\infty(G)$
of the kernels of the operators $R_h - I,$ for $h \in H.$ Therefore, $C^\infty(\spX)$
is closed. The  remaining assertion follows by application of
Corollary \refer{c: CirG moderate growth}.
\qed

\begin{rem}
It follows from \bib{Bps2}, Cor.~12.2, that the space $C_r^\infty(\spX)$ equals the space
$\cA_{N}(G/H),$ with $N = 2r,$ defined in \bib{BD}, \S~1;
the definition in the last mentioned
paper is given for $N \in \N,$ but makes sense for  arbitrary real
 $N\geq 0.$
Accordingly, Lemma \refer{l: CirX is closed} is due to \bib{BD}; see
loc.~cit.\
Lemma 1.
\end{rem}

If $V$ is a locally convex space, we denote its continuous linear dual by $V'.$
Unless otherwise specified, we equip it with the strong dual topology.

\begin{cor}
\naam{c: of automatic continuity}
Let $(\pi, V)$ be a smooth Fr\'echet representation of $G$ of moderate growth,
such that $V_K$ is finitely generated. Let $r \in \R$ and let
$T: V_K \to C^\infty_r(\spX)$ be a $(\fg, K)$-equivariant linear map.
\begin{enumerate}
\itema
The map $T$
has a unique extension to a continuous linear $G$-equivariant map
$\tilde T: V \to C^\infty_r(\spX).$
\itemb
The linear functional
$\ev_e \after T: v \mapsto Tv(e)$
has a unique extension to a continuous linear functional $\eta_T \in V'.$
\itemc
The functional $\eta_T$ is $H$-invariant and $\tilde T$ may be represented as
the generalized matrix coefficient map given by
$$
 T(v)(x) = \eta_T(\pi(x)^{-1} v),\qquad (x \in G/H).
$$
\end{enumerate}
\end{cor}

\proof
{}From Proposition \refer{p: automatic continuity and CirG} it follows that
$T$ has a unique extension to a continuous linear $G$-equivariant map
$\tilde T: V \to \CirG.$ The image of $\tilde T$ is a closed subspace $W$ of
$\CirG$ which contains the image of $T$ as a dense subspace.
In view of
Lemma \refer{l: CirX is closed} it follows that $W\subset C^\infty_r(\spX).$
The extended functional is given by
$\eta_T = \ev_e \after \tilde T: v \mapsto \tilde T v (e).$
The assertions of (c) readily follow by $G$-equivariance.
\qed

Using the above result we shall be able to express our Eisenstein integrals
as matrix coefficients of principal series representations.
As a preparation we need to
relate the function $\|\dotvar\|,$ defined in (\refer{e: defi norm}),
to the $G = K\Aq H$-decomposition.
Following \bib{BSpl1}, Eqn.~(\refbplm), we define the distance function
$\lspX: G \to [0, \infty\,[$
by
$$
\lspX(kah) = |\log a|,
$$
for $k \in K, a \in \Aq$ and $h \in H.$

\begin{lemma}
\naam{l: distance and norm}
There exists a constant $s >0$ such that
$$
e^{\lspX(x)} \leq \|x\|^s,\qquad (x \in G).
$$
\end{lemma}

\proof
One readily sees that it suffices to prove this in case $G = {}^\circ G.$
Moreover, since the functions of $x$
on both sides of the equality are left $K$-invariant, we may reduce
to the case that $G$ is connected and semisimple, with finite center.
{}From \bib{BSbv}, Lemma 14.4, we deduce, using the equality $\|x\| = \|x^{-1}\|,$
that  $\|a\| \leq \| a h\|$ for all $a \in \Aq$ and $h \in H.$
Hence, by the $G = K\Aq H$ decomposition,
it suffices to prove the estimate
\begin{equation}
\naam{e: estimate e log a}
e^{|\log a|} \leq \|a\|^s \qquad (a \in \Aq),
\end{equation}
for some $s >0$ independent of $a.$
Let $m$ be the minimal value of the continuous function
$\max \{\ga \mid \ga \in \gS(\fa_\fp)\}$ on the unit sphere in $\faq.$ Then $m>0.$
Using (\refer{e: norm and KAK}) we see that the estimate
(\refer{e: estimate e log a}) holds for $s \geq m^{-1}.$
\qed

Let $Q \in \allparabs$ and $\xi \in \discserQ$
be fixed throughout the rest of this section.

\begin{lemma}
\naam{l: Eis of moderate growth new}
Let $\types \subset \dK$ be a finite subset,
let $\psi \in \cAtwoQ(\tautypes)_\xi$ and let $\nu_0 \in \faQqdc$ be a regular point
for the Eisenstein integral $\nE(Q\col \psi\col \nu).$
There exist an open neighborhood $U$ of $\nu_0$ and
a constant $r >0$ such that $\nu \mapsto \nE(Q\col \psi \col \nu)$
is a function on $U$ with values in $C_r^\infty(\spX) \otimes \Vtypes.$ Moreover,
the mentioned function is locally bounded.
\end{lemma}

\proof It follows from \bib{BSpl1}, Prop.~\refplu, combined with
Lemma \refer{l: distance and norm}, that there exist an open
neighborhood $\Omega$ of $\nu_0$ and a polynomial function $p \in
\Pi_{\gSr(Q)}(\faQqd),$ such that the function $\nu \mapsto p(\nu)
\nE(Q\col \nu)$ is holomorphic on $\Omega$ as a function with
values in $C_r^\infty(\spX) \otimes \Hom(\cAtwoQ(\tautypes),
\Vtypes)$ and such that for every continuous seminorm $\mu$ on the
latter tensor product space, the function $\nu \mapsto \mu(p(\nu)
\nE(Q\col \nu))$ is a bounded on $\Omega.$

Select an open neighborhood $U$ of $\nu_0$ with compact closure contained
in $\Omega$ such that $\nE(Q\col \psi \col \dotvar)$ is holomorphic on an open
neighborhood of $\overline U.$ Then it follows by a straightforward application of
Cauchy's integral formula in the variable $\nu,$
see, e.g., \bib{Bps2},
proof of Lemma 6.1, that for every continuous seminorm
$\mu'$ on $C_r^\infty(\spX) \otimes \Vtypes$ the function
$\nu \mapsto \mu'(\nE(Q\col \psi\col \nu))$ is bounded on $U.$
\qed

\begin{lemma}
\naam{l: principal series of moderate growth}
Let $\nu \in \faQqdc.$
The representation $\pi_{Q, \xi, \nu}$ of $G$ in $C^\infty(K\col \xi)$ is
a smooth Fr\'echet representation of moderate growth. Moreover, the associated
$(\fg, K)$-module
$C^\infty(K\col \xi)_K$ is finitely generated.
\end{lemma}

\proof
It follows from Remark \refer{r: smooth vectors} that
$V := C^\infty(K\col \xi),$ equipped with
$\pi_{Q, \xi, \nu},$ is the space of
$C^\infty$-vectors for the Hilbert representation $\Ind_Q^G(\xi \otimes \nu \otimes 1).$
It now follows from \bib{Wal2}, Lemma 11.5.1, that $V$ is a smooth Fr\'echet $G$-module of
moderate growth. The last assertion is well known, see also
Proposition \refer{p: locally uniform generators} for a stronger assertion.
\qed

\begin{prop}
\naam{p: continuity J}
Let $\nu \in \faQqdc\setminus \cup \Hyp(Q,\xi).$
\begin{enumerate}
\itema
There exists a constant
$r \in \R$ such that
$J_{Q,\xi, \nu}$ maps
$\bar V (Q,\xi) \otimes \Ci(K\col \xi)_K$ into the space $C_r^\infty(\spX).$
\itemb
Let $r \in \R$ be a constant as in (a).
The map $J_{Q,\xi, \nu}$ has a unique extension
to a continuous linear map from $\bar V(Q,\xi)\otimes  \Ci(K\col \xi)$ into
$C^\infty_r(\spX).$ The extension intertwines
the $G$-representations $I \otimes \pi_{Q,\xi, -\nu}$
and $L.$
\end{enumerate}
\end{prop}

\proof
Fix $\nu$ as above.
By Lemma \refer{l: principal series of moderate growth}, there exists
a finite subset $\types\subset \dK$ such that $\Ci(K\col \xi)_\types$
generates $C^\infty(K\col \xi)_K$ as the $(\fg,K)$-module associated
with $\pi_{Q,\xi, - \nu}.$ Let $r \in \R$ be associated with $\types$
as in Lemma \refer{l: Eis of moderate growth new}. Then it follows from
(\refer{e: defi of J})
that $J_{Q,\xi, \nu}$
maps $\bar V(Q, \xi) \otimes \Ci(K\col \xi)_\types$ into $C^\infty_r(\spX).$
The map $J_{Q, \xi, \nu}$ is $(\fg, K)$-equivariant, by Theorem
\refer{t: intertwining prop J}. Since  $\Ci(K\col \xi)_\types$ generates
$\Ci(K\col \xi)_K,$ whereas $C^\infty_r(\spX)$ is $(\fg, K)$-invariant,
assertion (a) follows.

Assume that $r$ is a constant as in (a). Then it follows from
Theorem \refer{t: intertwining prop J} that the map $J_{Q,\xi, \nu}$
is $(\fg, K)$-equivariant. In view of Lemma
\refer{l: principal series of moderate growth}
and assertion (a), we may apply
Corollary \refer{c: of automatic continuity}
with $T = J_{Q,\xi, \nu}.$
Assertion (b) follows.
\qed

If $\nu \in \faQqdc\setminus \cup \Hyp(Q,\xi),$
we denote the continuous linear extension of $J_{Q,\xi, \nu}$
by the same symbol.
We denote the conjugate of the topological linear dual of $\Ci(K \col \xi)$ by
$C^{-\infty}(K\col \xi).$ The $G$-representation on the latter space induced
by dualization of $\pi_{\xi, \nu},$ is denoted by $\pi_{\xi, \nu}^{-\infty}.$

The sesquilinear pairing $\Ci(K\col \xi) \times \Ci(K\col \xi) \to \C,$ given
by (\refer{e: sesquilinear pairing}) induces a continuous linear
embedding $C^\infty(K\col \xi) \embeds C^{-\infty}(K\col \xi),$
intertwining the representations $\pi_{\xi, - \bar \nu}$
and $\pi_{\xi, \nu}^{-\infty}.$ The latter may therefore be viewed as
the continuous linear extension of $\pi_{\xi, - \bar \nu}.$ Accordingly, we shall sometimes
use the notation $\pi_{\xi, - \bar \nu}$ for the representation $\pi_{\xi, \nu}^{-\infty}.$

We denote by $V(Q,\xi)$ the
conjugate space of $\bar V(Q, \xi),$  and define the linear map
$\nj(Q\col\xi\col \bar \nu): V(Q,\xi) \to C^{-\infty}(K\col \xi)$
by
\begin{equation}
\naam{e: defi of nj}
\hinp{\gf}{\nj(Q\col\xi\col \bar \nu)(\eta)} = J_{Q,\xi, \nu}(\eta \otimes \gf)(e).
\end{equation}
Then by Proposition \refer{p: continuity J} and Corollary \refer{c: of automatic continuity},
the image of $\nj(Q\col\xi\col\bar \nu)$ is contained in the subspace
of $C^{-\infty}(K\col \xi)$ consisting
of $H$-invariants
for the representation $\pi_{\xi, \bar \nu};$
we agree to denote this subspace by $C^{-\infty}(Q\col \xi\col \bar \nu)^H.$

We may now represent the Eisenstein integral as a matrix coefficient.
The following formula generalizes the similar formula for $Q$ minimal,
see \bib{BSft}, Eqn.\ (53).

\begin{lemma}
\naam{l: nE as matrix coefficient}
Let $\nu \in \faQqdc\setminus\cup\Hyp(Q,\xi),$ let $\types\subset \dK$ be a finite
subset and let
$T = \eta \otimes \gf \in \bar V(Q, \xi)\otimes C^\infty(K\col \xi)_\types.$ Then
$$
E_\types(Q\col \psi_T \col \nu)(x)(k) =
\hinp{\gf}{\pi_{Q, \xi, \bar \nu}(kx)\nj(Q\col \xi\col \bar \nu)\eta}
\qquad
(x \in \spX, \; k\in K).
$$
\end{lemma}

\proof
This follows from (\refer{e: defi of J}) and (\refer{e: defi of nj}),
by application of Corollary \refer{c: of automatic continuity} (c).
\qed

To identify our
Eisenstein integral with the one introduced
by P.~Delorme in \bib{D1n}, we recall some results from \bib{CDdvn}, \S 2.4.

For each $v \in \QcW,$ we denote by $\cV(Q,\xi, v)$ the space of
$M_Q \cap vHv^{-1}$-fixed elements in $\cH_\xi^{-\infty},$ the conjugate
of the topological linear dual of $\cH_\xi^\infty.$ The space
$\cV(Q,\xi,v)$ is finite dimensional
by \bib{Bfm}, Lemma 3.3.
We introduce the formal direct sum
$$
\cV(Q,\xi):= \oplus_{v \in \QcW}\;\; \cV(Q, \xi, v).
$$
If $u \in C^{-\infty}(Q\col \xi \col\nu)^H,$
then on an open neighborhood of any $v \in \QcW$ in $K,$ the functional
$u$ may be represented by a unique continuous function
with values in $\cH_\xi^{-\infty},$
via the sesquilinear pairing (\refer{e: sesquilinear pairing}).
Its value $\ev_v (u)$ in $v$ is therefore a well defined element
of $\cV(Q, \xi, v).$ See \bib{CDdvn}, \S~3.3, for details.
The direct sum of the maps $\ev_v,$ for $v \in \QcW,$
is denoted by
$$
\ev = \oplus_{v \in \QcW} \;\;\ev_v:\;\;\;  C^{-\infty}(Q\col \xi \col\nu)^H \to \cV(Q,\xi).
$$
We have the following result, due to \bib{Bps1} for minimal $Q$
and to \bib{CDdvn} in general.

\begin{thm}
\naam{t: CD theorem on j}
There exists a unique meromorphic function $j(Q,\xi, \dotvar)$ on $\faQqdc$
with values in $\Hom(\cV(Q,\xi), C^{-\infty}(K\col \xi))$
such that the following two conditions are fulfilled.
\begin{enumerate}
\itema
For regular values of  $\nu,$ the image of $j(Q\col \xi\col  \nu)$ is contained
in  $C^{-\infty}(Q\col \xi\col \nu)^H.$%
\itemb
For regular values of $\nu,$ we have
$
\ev \after j(Q\col \xi \col \nu) = I_{\cV(Q,\xi)}.
$
\end{enumerate}
There exists a locally finite collection $\Hyp = \Hyp(j,Q,\xi)$ of
hyperplanes in $\faQqdc$ such that each
$\nu \in \faQqdc\setminus \cup \Hyp$ is a regular value for $j(Q\col \xi \col\dotvar)$ and
the associated map $j(Q\col\xi\col \nu)$
is surjective from $\cV(Q,\xi)$ onto $C^{-\infty}(Q\col \xi\col \nu)^H.$

Finally, each $\nu \in \faQqdc$ with $\Re \nu + \rho_Q$ stricly $\gSr(Q)$-anti-dominant
is a regular value for $j(Q\col \xi\col \dotvar).$ Moreover, for such $\nu$
and every $\eta \in \cV(Q,\xi),$ the
element $j(Q\col\xi\col \nu)\eta \in C^{-\infty}(K\col \xi)$
is representable by a continuous function $u: K \to \cH_\xi^{-\infty},$
in the sense that
$$
\hinp{\gf}{j(Q\col\xi\col\nu)\eta} =
\int_{K} \hinp{\gf(k)}{u(k)}\; dk,\qquad (\gf \in \Ci(K\col \xi)).
$$
\end{thm}

\proof
This follows from \bib{CDdvn}, Prop.\ 2, Thm.\ 1 and Thm.\ 3.
\qed

The Eisenstein integrals of Delorme are built in terms of matrix
coefficients coming from a subspace $\cV_{ds}(Q,\xi)$ of $\cV(Q,\xi),$
which is defined as follows, see \bib{D1n}, \S 8.3.
Let $v \in \QcW.$ An element
$\eta \in \cV(Q,\xi, v)$ naturally determines the $M_Q$-equivariant
embedding $\iota_\eta: \cH_\xi^\infty \to \Ci(M_Q/M_Q \cap vHv^{-1}),$
given by
$$
\iota_\eta(v)(m) = \hinp{\xi(m)^{-1} v}{\eta} \qquad (m \in M_Q).
$$
We denote by $\cV_{ds}(Q,\xi, v)$ the subspace of $\eta \in \cV(Q,\xi)$
with the property that $\iota_\eta$ maps into $L^2(M_Q/M_Q \cap vHv^{-1})^\infty.$
Note that for such $\eta$ the map $\iota_\eta$ extends to a continuous linear
map $\cH_\xi \to L^2(\spXQv);$ see \bib{CDn}, Lemma 1.
Moreover, the map
$\eta \mapsto \iota_\eta$ defines a linear isomorphism from $\cV_{ds}(Q, \xi, v)$
onto $V(Q, \xi, v),$  via which we shall identify.

We define the subspace $\cV_{ds}(Q,\xi)$ of $\cV(Q,\xi)$ as the direct
sum of the spaces $\cV_{ds}(Q,\xi, v),$ for $v \in \QcW.$ Via the direct
sum of the above isomorphisms, we obtain the natural isomorphism
$$
\cV_{ds}(Q,\xi) \simeq V(Q, \xi).
$$
Accordingly, the map $\nj$ introduced in (\refer{e: defi of nj})
may be viewed as a linear map
$$
\nj(Q\col \xi\col \nu): \;\; \cV_{ds}(Q,\xi) \to C^{-\infty}(Q\col \xi \col \nu)^H,
$$
defined for $\nu \in \faQqdc \setminus \cup \Hyp(Q,\xi).$

To relate this map with the map $j(Q\col \xi\col \nu)$ of Theorem \refer{t: CD theorem on j}
we need  standard intertwining operators.
We recall from \bib{Wal2} and \bib{CDdvn} that for a parabolic subgroup $P \in \allparabs$
with split component equal to $A_Q,$
the standard intertwining operator
$A(Q \col P \col \xi \col \nu)$ between the representations
$\pi_{P, \xi, \nu}$ and $\pi_{Q, \xi, \nu}$
on $\Ci(K \col \xi)$ is given by an absolutely convergent integral
for $\nu \in \faQqdc$ with $\inp{\Re\nu - \rho_Q}{\ga} > 0$ for every
$\ga \in \gSr(P) \cap \gSr(\bar Q),$
and allows a meromorphic continuation in $\nu.$ Its adjoint is a continuous
linear endomorphism of $C^{-\infty}(K\col \xi),$ intertwining
the representations $\pi_{Q,\xi, \nu}^{-\infty}$ and $\pi_{P,\xi, \nu}^{-\infty}.$
It extends the standard intertwining operator $A(P \col Q \col \xi \col -\bar \nu),$
and is therefore denoted by the same symbol. Thus,
\begin{equation}
\naam{e: adjoint of A}
A(Q \col P \col \xi \col \nu)^* =
A(P \col Q \col \xi \col -\bar \nu).
\end{equation}
We also recall that
$$
A(P \col Q \col \xi \col \nu)\after A(Q\col P \col \xi \col  \nu) =
\eta(Q\col P\col \xi\col \nu) \,
I_{\Ci(K \col \xi)}
$$
with $\eta(Q\col P\col \xi \col \dotvar)$ a non-zero scalar meromorphic function on
$\faPqdc = \faQqdc.$
In particular, it follows that the standard intertwining operator is invertible
for $\nu$ in an open dense subset of $\faQqdc.$

\begin{lemma}
\naam{l: j and limit}
Let $\nu \in \faQqdc$ be such that $\Re \nu - \rho_Q$ is strictly
$\gSr(Q)$-dominant.
Then, for every $\eta \in \cV(Q,\xi)$
and $\gf \in \Ci(K\col \xi)_K,$ and for each $v \in \QcW,$ all $m\in M_Q$
and all $X \in \faQqp,$
$$
\lim_{t \to \infty}
a_t^{- \nu + \rho_Q}\;
\hinp{\gf}{\pi_{\bar Q, \xi, \bar \nu}(m a_t \, v)\,
j(\bar Q \col \xi\col \bar \nu)\eta}
=
\hinp{A(Q \col \bar Q\col \xi \col - \nu)\gf(e)}{\xi(m)\eta_v},
$$
where $a_t = \exp tX.$
\end{lemma}

\proof
The result is equivalent to Lemma 16 in \bib{D1n}.
We refer to the proof given there.
\qed

\begin{thm}
\naam{t: identification of nj}
Let $\eta \in \cV_{ds}(Q,\xi).$ Then $\nj(Q\col \xi\col \dotvar)\eta$ is
holomorphic as a function on
$\faQqdc \setminus\cup\Hyp(Q,\xi),$ with values in $C^{-\infty}(K\col \xi).$
Moreover,
$$
\nj(Q\col \xi\col \nu)\eta =
A(\bar Q \col Q  \col \xi \col \nu)^{-1} j(\bar Q \col \xi \col \nu)\eta,
$$
as an identity of meromorphic
$C^{-\infty}(K\col \xi)$-valued functions in $\nu \in  \faQqdc\setminus\cup\Hyp(Q,\xi).$
In particular, $\nj(Q\col \xi\col \dotvar)\eta$ extends to a meromorphic
$C^{-\infty}(K\col \xi)$-valued function on $\faQqdc.$
\end{thm}

For the proof of this result we need the following lemma.

\begin{lemma}
\naam{l: limit nE one}
Let $\types \subset \dK$ be a finite subset.
There exists an open dense subset $\Omega$ of the set of points $\nu \in \faQqdc$ with
$\Re \nu$ strictly $\gSr(Q)$-dominant, such that the following holds.
Let $\psi \in \cAtwoQ(\tautypes),$  $v \in \QcW,$ $m \in \spXQvp,$ $X \in \faQqp$
and put $a_t = \exp tX, (t \in \R).$
Then, for every $\nu \in \Omega,$
$$
\lim_{t \to \infty} a_t^{-\nu + \rho_Q}\,\nE(Q\col \psi \col \nu)(ma_t v)
=
(\psi)_v(m).
$$
\end{lemma}

\proof
Let $\omega$ be the set of regular points for $\nE(Q\col \dotvar),$
and $\omega_+$ the subset of $\nu \in \omega$ with
$\Re \nu$ strictly $\gSr(Q)$-dominant.

Fix a minimal parabolic group $P$ from $\allparabs,$ contained in $Q.$
Then, by \bib{BSpl1}, Prop.~\refbpln, the family $f: (\nu, x) \mapsto \nE(Q\col \nu \col x)\psi$
belongs to $\cE^\hyp_Q(\spX\col \tau).$
Moreover, for each $u \in \NKaq,$ the set of exponents
$\Exp(P,u\mid f_\nu)$ is contained in the collection
$W^{P|Q} (\nu + \gL(P|Q)) - \rho_P -\N \gD(P),$ for $\nu \in \omega.$

Fix $\nu \in \omega_+$ and let $\xi$ be an exponent in $\Exp(Q, v\mid f_\nu).$
Then it follows by application of \bib{BSanfam}, Thm.\ 3.5,
that $\xi = w(\nu + \gL)|_{\faQq} - \rho_Q - \mu,$ for certain $w \in W^{P|Q},$
$\gL \in \gL(P|Q)$ and $\mu \in \N\gDr(Q).$ It follows from the definitions
preceding \bib{BSpl1}, Prop.~\refbpln,
that
$w\gL \in - \R^+ \gD(P).$ Hence $\Re\xi(X) + \rho_Q(X) \leq  w\Re \nu(X),$
with equality if and only if $w\gL|_{\faQq} = 0$ and $\mu = 0.$
Now $\Re \nu$ is strictly $\gSr(Q)$-dominant and $X \in \faQqp,$
Hence, by a well known result on root systems,
$\Re \nu(X) \geq  \Re s \nu (X),$ for each $s \in W,$
with equality if and only if $s$ centralizes $\faQq.$
Since $W^{P|Q} \cap W_Q = \{e\},$
we conclude that
$$
\Re \xi (X) < (\Re \nu - \rho_Q)(X)
$$
for every exponent $\xi \in \Exp(Q, v\mid f_\nu),$ different from $\nu - \rho_Q.$
It follows that
\begin{equation}
\naam{e: lim f as q}
\lim_{t \to \infty} a_t^{-\nu + \rho_Q}\,f_\nu(ma_t v)
=
\lim_{t \to \infty} q_{\nu - \rho_Q}(Q,v\mid f_\nu, tX )(m),
\end{equation}
for every $\nu \in \omega_+$ for which the limit on the right-hand
side exists. It follows from \bib{BSpl1},
Def.~\refplb{}and Prop.~\refplc{},
that there exists a non-empty open subset $\Omega_0$ of
$\faQqdc$ such that
\begin{equation}
\naam{e: expression for q nu minu rho}
q_{\nu - \rho_Q}(Q,v\mid f_\nu, X )(m) = \psi_v(m)
\end{equation}
for all $\nu \in \Omega_0,$ $v \in QcW,$ $m \in \spXQvp$ and $X \in \faQq.$
On the other hand, by \bib{BSanfam}, Thm.~7.7, there exists an open dense
subset $\Omega_1$ of $\faQqdc$ such that the expression on the
left-hand side of (\refer{e: expression for q nu minu rho}) depends holomorphically on $\nu \in \Omega_1,$
for all $v, m, X$ as above. By analytic continuation it follows
that (\refer{e: expression for q nu minu rho}) is valid for $\nu \in \Omega_1.$
This implies that the limit on the right-hand side of (\refer{e: lim f as q})
has the value $\psi_v(m),$ for every $\nu \in \omega_+ \cap  \Omega_0.$
\qed

{\bf Proof of Theorem \refer{t: identification of nj}:\ }
Fix $\gf \in \Ci(K\col \xi)_K$ with $\gf(e) \neq 0$ and put $T = \eta \otimes \gf.$
Select a finite subset $\types\subset \dK$ such that $\gf \in \Ci(K\col \xi)_\types.$
Then, in the notation of ,
$\psi_T \in \cAtwo(\tautypes).$
Fix $v \in \QcW.$ Then the preimage $M_{Q,v,+}$ of $\spXQvp$ under the canonical
map $M_Q \to \spXQv$ is open dense in $M_Q.$
Fix $m \in M_Q$
and $X \in \faQqp.$ We agree to write $a_t = \exp tX.$
Let $\cA$ be the open dense subset of $\faQqdc\setminus \cup \Hyp(Q,\xi)$
consisting of points $\nu$ where both intertwining operators
$\nu \mapsto A(Q \col \bar Q \col \xi \col -\nu)^\pm$
are regular.
By Lemma \refer{l: nE as matrix coefficient} we may write, for $\nu \in \cA,$
\begin{eqnarray*}
\lefteqn{E^\circ_\types(Q\col \psi_T\col \nu)(ma_t v)(e)
=}
\\
&=&
\hinp{\gf}{\pi_{Q, \xi, \bar \nu}(m a_t v) \nj(Q\col \xi \col \bar \nu)}
\\
&=&
\hinp{A(Q\col \bar Q\col \xi \col - \nu)^{-1}\gf}
{\pi_{\bar Q, \xi, \bar \nu}(m a_t v) A(\bar Q\col Q\col \xi \col \bar \nu)
\nj(Q\col \xi \col \bar \nu)\eta}.
\end{eqnarray*}
Replacing $\cA$ if necessary, we may in addition assume that the conjugate $\bar \cA$ of
$\cA$ has empty intersection with the set $\cup\Hyp,$ where $\Hyp = \Hyp(j, \bar Q,\xi)$
is as in
Theorem \refer{t: CD theorem on j}, with $\bar Q$ in place of $\bar Q.$
By the mentioned theorem it then follows,
for $\nu \in \cA,$ that
\begin{equation}
\naam{e: A on nj equals j}
A(\bar Q\col Q\col \xi \col \bar \nu) \nj(Q\col \xi\col \bar \nu)\eta =
j(\bar Q \col \xi \col \bar \nu)\eta(\bar \nu),
\end{equation}
for $\eta(\bar \nu) \in \cV(\bar Q, \xi)$ given by
$
\eta(\bar \nu) = \ev \after A(\bar Q\col Q\col \xi \col \bar \nu)
\nj(Q\col \xi \col  \bar \nu)\eta.
$
Using Lemma \refer{l: j and limit} we now conclude that, for $\nu \in \cA$ with
$\Re \nu - \rho_Q$ strictly $\gSr(Q)$-dominant,
\begin{eqnarray}
\lefteqn{\lim_{t \to \infty} a_t^{-\nu + \rho_Q }
 \,E^\circ_\types(Q\col \psi_T\col \nu)(m a_t v)(e) }\nonumber\\
&=& \hinp{\gf(e) }{\xi(m) \ev_v \after A(\bar Q\col Q\col \xi \col \bar \nu)
\nj(Q\col \xi \col \bar \nu)\eta}
\nonumber\\
&=&
\hinp{\gf(e)}{\eta(\bar \nu)_v}.
\naam{e: limit nE two}
\end{eqnarray}
In particular, this holds for $\nu$ contained in the non-empty open set
$\cA \cap \Omega,$ with $\Omega$ as in Lemma
\refer{l: limit nE one}. For such $\nu$ it follows by the mentioned lemma
that the limit in
(\refer{e: limit nE two}) also equals $(\psi_T)_v(m) = \hinp{\gf(e)}{\xi(m) \eta_v}.$
We deduce that, for $\nu \in \bar \cA \cap \bar \Omega,$ where the bar denotes conjugation,
\begin{equation}
\naam{e: identity with gf e}
\hinp{\gf(e) }{\xi(m)\, \eta(\nu)_v}
= \hinp{\gf(e)}{\xi(m)\eta_v},
\end{equation}
for all $m \in M_{Q,v,+}.$
By continuity and density it follows that the identity
(\refer{e: identity with gf e}) holds for all $m \in M_Q.$ Since $\gf(e) \in \Hyp_\xi^\infty\setminus\{0\},$
it follows by irreducibility
of the $G$-module $\Hyp_\xi^\infty$ that $\eta(\nu)_v = \eta_v,$
for all $\nu \in \bar \cA \cap \bar \Omega.$ This identity holds for every $v \in \QcW,$
since the sets $\cA$ and $\Omega$ are independent
of the element $v \in \QcW.$ Combining this with (\refer{e: A on nj equals j})
we deduce that, for every
$\nu \in \bar \cA \cap \bar \Omega,$
\begin{equation}
\naam{e: indentity nj and A j}
\nj(Q\col \xi\col \nu) \eta = A(\bar Q \col Q\col \xi \col \nu)^{-1} j(\bar Q\col \xi \col \nu)\eta.
\end{equation}
Let $f(\nu)$ denote the expression on the left-hand side and $g(\nu)$ that on the right-hand
side of the above equation. Then $g$ is a meromorphic $C^{-\infty}(K\col\xi)$-valued
function on $\faQqdc,$ by Theorem \refer{t: CD theorem on j} and meromorphy of the intertwining
operator.
If $\gf \in \Ci(K\col \xi)_K,$
then $\nu \mapsto \hinp{f(\nu)}{\gf}$ is a holomorphic function of
$\nu \in \faQqdc\setminus\cup\Hyp(Q,\xi),$ by Lemma \refer{l: nE as matrix coefficient}.
On the other hand, $\nu \mapsto \hinp{g(\nu)}{\gf}$ is a meromorphic
function on $\faQqdc.$ By analytic continuation we deduce that
\begin{equation}
\naam{e: f equals g when tested}
\hinp{f(\nu)}{\gf} = \hinp{g(\nu)}{\gf},
\end{equation}
as an identity of meromorphic functions in the variable $\nu \in \faQqdc\setminus \cup\Hyp(Q,\xi).$
{}From the holomorphy of the function on the left-hand side it follows that the function
on the right-hand side is actually regular on $\faQqdc\setminus \cup\Hyp(Q,\xi),$
for every $\gf \in C^{-\infty}(K\col \xi)_K.$ The latter space is dense in
$C^\infty(K\col \xi)$ and $\nu \mapsto gf(\nu)$ is a meromorphic $C^{-\infty}(K\col \xi)$-valued
function. It follows that $g$ is regular on $\faQqdc\setminus \cup\Hyp(Q,\xi).$ It now follows
from (\refer{e: f equals g when tested}) that the element $f(\nu) \in C^{-\infty}(K\col \xi)$ equals
$g(\nu),$ for every $\nu \in \faQqdc\setminus \cup\Hyp(Q,\xi).$ This implies all assertions
of the theorem.
\qed

It follows from the above result that the distribution
vector $\nj(Q \col \xi \col \nu)\eta,$ defined for $\eta \in \cV_{ds}(Q,\xi),$
coincides with the similarly denoted distribution vector defined in \bib{BSft}, (3.13).

\begin{cor}
\naam{c: comparison nE}
Let $(\tau, \Vtau)$ be a finite dimensional unitary representation
of $K,$ let $P \in \allparabs$ and $\psi \in \cA_2(P, \tau).$
Then the Eisenstein integral
$\nE(P\col \psi \col \gl),$ for $\gl \in \faPqdc,$
coincides with the normalized Eisenstein integral $\nE(P,\psi, - \gl)$
defined in  \bib{CDn}, \S\ 5.1.
\end{cor}

\proof
By the functorial property of Lemma \refer{l: functoriality of nE}, which is satisfied
by both Eisenstein integrals, it suffices to prove the result
for $\tau = \tau_\types,$ with $\types\subset \dK$
a finite subset. By linearity it suffices to prove the assertion for $\psi = \psi_{\eta \otimes f},$
where $\xi \in \discserP,$ $\eta \in \cV_{ds}(P, \xi)$ and $f \in \Ci(K\col \xi)_\types.$
The associated normalized Eisenstein
integral is denoted  $\nE(P, \psi, \nu)$ in \bib{CDn}, \S 5.1.
It is represented as a matrix coefficient in \bib{CDn}, Prop.\ 4.
This representation  coincides with the one given in
Lemma \refer{l: nE as matrix coefficient}.
\qed

It follows from the  equality of the normalized Eisenstein
integrals stated above, that the Plancherel theorems formulated in
\bib{BSpl1}, \S~\refplh,
and Section \refer{s: The Plancherel decomposition},
coincide with the ones of P.\ Delorme formulated in
\bib{Dpl}, Sections 3.3 and 3.4. However, the chosen normalizations of measures
are different, resulting in different constants. We shall finish this section
by relating the various constants. The normalization of measures
for the present
paper is described in \bib{BSpl1}, \S~\refpli.
The normalization
given in \bib{Dpl}, \S 0, follows essentially the same conventions of interdependence, with
one crucial difference.
A choice of invariant measure $dx$  for $\spX,$ determines the same choice of Haar measure
$da$ for $\Aq$ in both papers. In our paper we fix the Lebesgue measure
$d\gl$ on $i \faqd$  that makes the Euclidean Fourier transform an isometry
from $L^2(\Aq, da)$ onto $L^2(i \faqd, |W|\, d\gl).$ On the other hand,
in \bib{Dpl}, \S 0, the convention is to fix the measure $d\gl':= |W| \, d\gl$ instead.

If $Q \in \allparabs,$ the same convention applies to the normalizations
of invariant measures $dx_{Q,v}$ on
$\spXQv,$ for $v \in \QcW,$ versus a choice of normalization of $da_Q$ on the group ${}^*A_{Q\iq},$
which is `the $\Aq$ of $(M_Q, M_Q \cap vHv^{-1}).$' This determines a normalization
$d\gl_Q$ of Lebesgue measure on $i\,{}^*\fa_{Q\iq}^*.$ The corresponding measure of
\bib{Dpl} is given by $d\gl_Q' = |W_Q| \, d\gl_Q.$ In both papers,
one chooses the measure on $i \faQqd$ to be the quotient of the chosen
measures on $i\faqd$ and $i\staQqd.$ This results in a choice of Lebesgue measure
$d\mu_Q$ on $i\faQqd$ in the present paper. The similar measure $d\mu_Q'$ in
\bib{Dpl} is then given by $d\mu_Q' = [W:W_Q] \,d \mu_Q.$ For the constants in the Plancherel
formula, see e.g. \bib{BSpl1}, Thm.~\refpla{}~(d), this means that
$[W:W^*_Q]$ should be replaced by
$$
[W:W_Q]^{-1}[W: W^*_Q] = [W_Q^*: W_Q]^{-1} = |W(\faQq)|^{-1}.
$$
The latter is indeed the constant occurring in, e.g., \bib{Dpl}, Thm.~3~(iii).

\def\adritem#1{\hbox{\small #1}}
\def\distance{\hbox{\hspace{3.5cm}}}
\def\apetail{@}
\def\adderik{\vbox{
\adritem{E.~P.~van den Ban}
\adritem{Mathematisch Instituut}
\adritem{Universiteit Utrecht}
\adritem{PO Box 80 010}
\adritem{3508 TA Utrecht}
\adritem{Netherlands}
\adritem{E-mail: ban{\apetail}math.uu.nl}
}
}
\def\addhenrik{\vbox{
\adritem{H.~Schlichtkrull}
\adritem{Matematisk Institut}
\adritem{K\o benhavns Universitet}
\adritem{Universitetsparken 5}
\adritem{2100 K\o benhavn \O}
\adritem{Denmark}
\adritem{E-mail: schlicht@math.ku.dk}
}
}
\vfill
\hbox{\vbox{\adderik}\vbox{\distance}\vbox{\addhenrik}}
\end{document}